\documentclass[reqno,12pt]{article}
\usepackage{amsmath,amsthm,amsfonts,amssymb,amscd}
\usepackage[dvips]{graphicx,color,epsfig}
\usepackage[T1]{fontenc}
\usepackage[latin1]{inputenc}
\newtheorem{thm}{Th\'eor\`eme}[section]
\newtheorem{cor}[thm]{Corollaire}

\newtheorem{lem}[thm]{Lemme}

\newtheorem{pro}[thm]{Proposition}

\newtheorem{defn}[thm]{D\'efinition}
\newtheorem{definition}[thm]{D\'efinition}

\newtheorem{rem}[thm]{Remarque}
\newtheorem{rems&defn}[thm]{Remarque et d\'efinition}

\newtheorem{exmp}[thm]{Exemple}

\newtheorem{rems}[thm]{Remarques}

\newcommand{\R}{\mathbb{R}}
\newcommand{\F}{\mathcal{F}}
\newcommand{\B}{\mathcal{B}}
\newcommand{\A}{{\cal A}}
\newcommand{\TT}{\mathcal{T}}

\newcommand{\h}{{\cal H}}

\newcommand{\G}{\mathcal{G}}
\newcommand{\ESS}[1]{\underset{#1}{\mathrm{ess~sup}}}
\newcommand{\INF}[1]{\underset{#1}{\mathrm{ess~inf}}}

\newcommand{\eps}{\varepsilon}
\newcommand{\1}{{\mathbf 1}}

\newcommand{\To}{\longrightarrow}

\newcommand{\biindice}[3]%.
{ #1\begin{array}[t]{c}
{\scriptstyle #2}\\
{\scriptstyle #3}\end{array}
}

\def \endproof {\quad \hfill  \rule{2mm}{2mm}\par\medskip}
\begin{document}

\title{Contr\^{o}le impulsionnel appliqu\'e \`a la gestion de
changement de  technologie dans une entreprise}
\author{Rim Amami \\
Institut de Mathématiques de Toulouse \\
Université Paul Sabatier,  Toulouse\\
rim.amami@math.univ-toulouse.fr }
\maketitle

\section{Introduction}

Le but de ce chapitre est l'étude d'un problème de contrôle
impulsionnel appliqué à la gestion du choix de technologie d'une
entreprise (voir par exemple A. Bensoussan et J.L. Lions
\cite{bensoussan}). Nous supposons que l'entreprise décide à
certains instants  de changer de technologie et  de  valeur de la firme  (par exemple une recapitalisation). 
 Les instants d'impulsion,
le choix de la nouvelle technologie et la loi des sauts sont des variables de décision, dont l'ensemble est appel\'e un contr\^ole impulsionnel. 
 Plus précisément, nous avons  la   suite  croissante  $(\tau_n)_n$ des instants d'impulsions de limite notée $\tau$,  la 
technologie $\zeta_{n+1}$ choisie à l'instant $\tau_n$ et $\Delta_n$ la taille
du saut du log de la valeur de la firme à l'instant $\tau_n$.  La  loi conditionnelle  du couple $(\zeta_{n+1},\Delta_n)$  ne dépend que de l'état du système en $\tau_n^-$.  
 On appelle contr\^ole impulsionnel la donn\'ee de tous ces param\`etres, soit la strat\'egie not\'ee $\alpha = (\tau_n,\zeta_{n+1},\Delta_n, n\geq -1).$\\
\\
En appelant $\xi_t$  le processus égal à $\zeta_{n+1}$ sur  $[\tau_n, \tau_{n+1}[$  et $Y$ est le
 processus représentant le log de la valeur de la firme.  Le bénéfice net de la firme est représenté par la  fonction $f$ et  le coût de changement de technologie est représenté par la fonction $c$.   Par  suite,  
   toute strat\'egie $\alpha$ occasionne un gain:
 \begin{equation}\label{gain}
 k(\alpha)=\int_0^{+\infty} e^{-\beta
s}f(\xi_s,Y_s)\;ds -\displaystyle\sum_{0 < \tau_n <\tau}
e^{-\beta \tau_n}c(\zeta_n,Y_{\tau_n^-},
\zeta_{n+1},Y_{\tau_n}), 
\end{equation}
où $\beta>0$ est un coefficient d'actualisation.\\
\\
 Notre  principale contribution est de prouver l'existence et fournir  une caractérisation d'une stratégie optimale 
   $\widehat{\alpha} $ qui maximise la fonction gain de la firme. \\
   \\
 Les outils math\'ematiques qui sont \`a la base d'une telle \'etude ont
 \'et\'e initi\'es par A. Bensoussan et J.L. Lions
\cite{bensoussan} et ensuite formalis\'es par d'autres auteurs.  Par exemple,  K.A. Brekke et  B. \O ksendal étudient un  problème de départ et d'arrêt. Le principal résultat de leur papier \cite {oksendal} est de trouver une suite d'arrêt optimale  du temps de  départ et d'arrêt d'un processus de production. Ils résolvent ce genre de problème en utilisant le calcul stochastique. 
\\
B. Bruder et   H. Pham considèrent un problème de contrôle   impulsionnel  en horizon fini pour les diffusions avec un décalage de décision et des délais d'exécution. Ces auteurs montrent que la fonction de valeur pour ce genre de problème   satisfait une version adaptée du principe de programmation dynamique et  ils  fournissent un algorithme pour trouver la stratégie optimale. 
\\
J.P. Lepeltier et B. Marchal  utilisent  dans \cite{leplprob} 
  une technique purement probabiliste pour la résolution du problème de contrôle impulsionnel. L'outil de base est la théorie générale du contrôle de C. Striebel qui permet d'obtenir un critère d'optimalité performant.\\
Nous mentionnons également B. Djehiche et al.  \cite{djehiche} et   M. Jeanblanc-S. Hamadène \cite {jeanblanc} qui
ont utilisé des outils purement probabilistes comme l'enveloppe de
Snell et les équations différentielles stochastiques rétrogrades
 pour
résoudre le problème optimal de changement de technologie en horizon
fini.
 \\
 Dans l'article  de P.A. Meyer \cite{meyer}, la valeur de la
 firme est mod\'elis\'ee suivant un mod\`ele canonique. La construction du
 probl\`eme  de contr\^ole  est fond\'ee sur la th\'eorie de la renaissance des
processus de Markov. Avant la premi\`ere impulsion, la loi du
syst\`eme est celle d'un processus de Markov tu\'e au moment de
cette impulsion. Puis apr\`es l'impulsion, on le fait renaitre
suivant une nouvelle loi de processus de Markov \`a
l'aide d'une probabilit\'e de transition.\\
M.H.A Davis a étudié dans \cite{Davis} un problème de contrôle
optimal déterministe. Il a introduit une simple formulation du
principe de la programmation dynamique des processus de Markov
déterministes par morceaux (PMDP) qui aident à résoudre ce type de
problème.
\\
M. Robin a abordé dans sa thèse \cite{Robin} un type de problème de 
contrôle impulsionnel avec retard déterministe, c'est à dire
qu'aucune décision ne peut être prise avant l'effet de la dernière
décision. Il a établi des résultats sur les problèmes des temps
d'arrêt optimaux essentiellement pour les processus de Markov
fellériens. Le résultat fondamental de cette thèse est la propriété
de continuité de la fonction valeur obtenue par des techniques de
pénalisation.\\
G. Mazziotto et J. Szpirglas \cite{szpirglas} ont étudié le contrôle
impulsionnel de systèmes stochastiques en information incomplète
selon des méthodes développées essentiellement dans le cadre de la
théorie du filtrage non linéaire. Leur principal résultat est un
théorème de séparation du contrôle et du filtrage dans une situation
de gestion de stock partiellement observée. Ce résultat a été obtenu
en étendant la méthode de M. Robin \cite{Robin} et en utilisant des
théorèmes de sélection \cite{parthasarathy}.
\\
Une approche diff\'erente est utilis\'ee pour la r\'esolution de ce
type de probl\`eme: le probl\`eme d'optimisation est formul\'e comme
un probl\`eme parabolique de contr\^{o}le impulsionnel avec trois
variables li\'ees \`a la fonction co\^{u}t, la technologie choisie
et la valeur de la firme (m\'ethode utilis\'ee par H. Pham, M. Mnif
et V. Vath \cite{Mnif}). Cette r\'esolution est associ\'ee au
principe de programmation dynamique des in\'equations
quasi-variationnelles de Hamilton-Jacobi-Bellman.
\\
Dans le même contexte, nous citons \cite{djehiche2,bahlali,Pham} qui
 caractérisent la fonction de valeur comme l'unique solution de
viscosité des in\'equations
quasi-variationnelles de Hamilton-Jacobi-Bellman.\\
\\
Le problème de contrôle   impulsionnel a   été étudié par J.P.  Lepeltier et B. Marchal \cite {lepl, lepeltier}. Néanmoins, leurs résultats ne  s'appliquent pas   aux situations qui nous intéressent, car nous construisons un modèle trajectoriel alors que ces deux auteurs définissent un modèle en  loi. Par ailleurs, alors que ces deux auteurs minimisent une fonction uniquement de coût (leur coût peut être défini en (\ref{gain})  par une somme au lieu d'une différence), nous cherchons à optimiser une fonction gain - coût.   Ainsi, nous ne considérons pas que la firme  pourrait disparaître: la meilleure stratégie après
la limite du temps $\tau$ est de continuer, en gardant la même technologie, au lieu d'aller à la faillite. De plus, en utilisant  des propriétés de Markov et des relations de récurrence (proposition \ref{recurence}), nous détaillons les preuves et nous présentons  une forme constructive 
   du contrôle impulsionnel.  En fait, nous utilisons un caractère markovien et homogène entre deux instants d'impulsion pour établir un critère d'optimalité (théorème  \ref{stratoptimal}). 
   \\
   \\
Notre modèle est inspiré de \cite {lepl, lepeltier} et peut être décrit comme suit. 
A chaque temps d'impulsion, on choisit une nouvelle technologie  afin d'améliorer
le profit de la firme.  En utilisant la théorie des processus de Markov, nous prouvons l'existence d'une stratégie optimale qui maximise la fonction valeur. Au lieu de la construction canonique étudié dans \cite {lepl, lepeltier}, nous choisissons
  une approche markovienne qui   permet des expressions plus explicites.  D'ailleurs,    nous retrouvons les résultats de J.P.  Lepeltier et B. Marchal à  l'aide  d'outils différents: au lieu de la théorie générale du contrôle   impulsionnel, nous donnons une forme constructive de l'évolution du système impulsé. Ensuite, nous introduisons la définition    \ref{strategadm} qui restreint l'ensemble des stratégies. En effet, nous  imposons
   deux propriétés définissant une stratégie admissible: la première assure l'intégrabilité des fonctions $f$ et $c$ et la deuxième suppose que  la loi entre deux sauts est donnée par la loi du couple initial (temps de saut, processus d'état sur ce premier intervalle).  
   \\ 
  Dans ce chapitre nous nous plaçons sous    trois   hypothèses:  des conditions de compacité sur  l'ensemble des noyaux de transition et ses sections (hypothèses 1 et 3) et des propriétés de continuité concernant  les lois  de $(\tau_0,(i,Y_.)\1_{[0,\tau_0[})$ sachant la condition initiale $(i,x),$ les fonctions de gain et de coût et les noyaux de transition  (hypothèse 2).  Nos trois hypothèses  remplacent les sept hypothèses introduites par J.P. Lepeltier et B. Marchal \cite{lepl}.  \\ 
    \\
Notre travail est organisé comme suit:  la  section 2  est
consacrée à définir le modèle correspondant au problème de
contr\^{o}le impulsionnel 
 ainsi que les filtrations
associ\'ees.  Dans  la section 3, nous
établissons un critère d'optimalit\'e. Nous \'enon\c{c}ons tout
d'abord le principe de programmation dynamique \`a l'aide du gain
maximal conditionnel apr\`es un temps d'arr\^{e}t $\theta$. Ensuite,
par des techniques markoviennes, nous \'etablissons un lien entre le
gain maximal conditionnel apr\`es
  $\theta$ et la fonction valeur de la firme ce qui nous permet de
d\'eduire un crit\`ere d'optimalit\'e dépendant  de cette fonction.  Dans la section 4,  nous définissons une stratégie
qui maximise la fonction valeur de la firme et qui réalise
l'optimalité conditionnelle    et  nous donnons un exemple qui vérifie
les hypothèses requises.  Nous terminons ce chapitre (section \ref{section5}) par  une comparaison entre les  hypothèses et définitions que nous utilisons ici   et celles introduites par J.P.  Lepeltier et B.  Marchal  \cite {lepl, lepeltier}.

  \section{Mod\'elisation }\label{section 2}
Soit $(\Omega,\F, (\F_t)_{t \geq
0}, \mathbb{P})$ un espace de probabilit\'e
 muni d'une filtration $(\F_t)_{t \geq
0}$ compl\`ete continue \`a droite  et soit un mouvement brownien $W
= (W_t)_{t \geq 0 }$. 
Nous noterons par $(\G_t)_{ t  > 0}$ la filtration
 d\'efinie par $\G_t = \displaystyle\vee_{s < t}\F_s,\; \forall\; s<t$, (elle est notée usuellement $\F_{t^-}$).\\
\\
  Un  contr\^ole impulsionnel consiste en une succession de
changements  d'état \`a $(\tau_{i})_{i \geq-1}$, où  $ \tau_{-1} = 0$, suite croissante de
 $\G\mbox{-temps d'arr\^{e}t}$ de limite not\'ee $\tau$.  Nous considérons les suites
  de temps d'arrêt vérifiant la propriété suivante: Soit  
{\small  \begin{equation}
 \label{Ptau}
 N(\omega) = \inf\{n: \tau_{n+1}(\omega) = \tau_{n}(\omega) \}\mbox{, 
 alors sur }\{N<+\infty\},~\tau_{k+1}(\omega) = \tau_{k}(\omega) \; \forall k\geq N(\omega).
 \end{equation}
 La suite $(\tau_n)$  est  strictement croissante
sur l'événement complémentaire~ $\{N=+\infty\}$.

\begin{pro}
\label{Nom} Le temps $\tau_{N}$ est un $\G$-temps d'arr\^et.
\end{pro}
\noindent\textbf{Preuve:} 
Soit
\begin{eqnarray*} 
 \{\tau_{N(\omega)} \leq t\} &=& \bigcup_k\{N(\omega) =k\} \cap  \{\tau_k(\omega) \leq t\}\nonumber\\
 & = &
\bigcup_k \{\tau_{k+1}(\omega) = \tau_{k}(\omega)  \}\cap \{ \forall j<k, \tau_{j+1}(\omega) < \tau_{j}(\omega) \}  \cap  \{\tau_k(\omega) \leq t\}. 
\end{eqnarray*} 
Du fait que les instants $\tau_k$ sont des $\G\mbox{-temps d'arr\^{e}t}$,  nous aurons $\{\tau_k(\omega) \leq t\} \in \G_t$ et  $ \{\tau_{k+1}(\omega) = \tau_{k}(\omega)  \}\cap \{ \forall j<k, \tau_{j+1}(\omega) < \tau_{j}(\omega) \}  \in \G_{\tau_k}$. \\
Par conséquent,  $\{\tau_{N(\omega)} \leq t\}  \in \G_t.$
\endproof

Soit $U$ l'ensemble fini des technologies possibles et $\mathcal{P}(U)$ la tribu triviale. 
Nous introduisons le processus
c\`adl\`ag $\xi$ \`a valeurs dans $U$:
\begin{equation}\label{xi}
   \xi_t = \xi_0 \1_{[0,\tau_0[}(t)+ \displaystyle\sum_{n \geq 0} \zeta_{n+1}
\1_{[\tau_n,\tau_{n+1}[}(t) +  \xi_{\tau^-} \1_{[\tau,+\infty[}(t).
\end{equation}
  Pour tout $n \in \mathbb{N},$  $\zeta_{n+1}$ la
technologie choisie \`a l'instant $\tau_n$, est une variable
al\'eatoire
$\G_{\tau_n}\mbox{-mesurable}$.

\begin{pro} Le processus  $\xi$ est $\G\mbox{-adapt\'e}$.
\end{pro}
\noindent\textbf{Preuve:} 
Les instants $\tau_n$ sont des $\G\mbox{-temps d'arr\^{e}t}$. Alors,
$\{\tau_n \leq t\} \in \G_t,$ et puisque
$\zeta_{n+1}$ est $\G_{\tau_n}\mbox{- mesurable}$,
$\zeta_{n+1}\1_{[\tau_n,\infty[}(t)$ est $\G_t \mbox{- mesurable}$.
De plus, par stabilit\'e des tribus par passage au compl\'ementaire
nous avons $\{\tau_{n+1} > t\} \in \G_t.$ D'o\`{u}
$\1_{\{\tau_{n+1} > t\}}$ est $\G_t \mbox{- mesurable}$.  
 Ainsi $\xi_t$ est
$\G_t\mbox{-mesurable}$ comme  somme de variables
$\G_t\mbox{-mesurables}$.  Par cons\'equent,  le processus $\xi$ est
$\G\mbox{-adapt\'e}$.
\endproof
La valeur de la firme entre deux instants d'impulsions
du syst\`eme est donn\'ee par  
 $S_t = \exp Y_t,~ t\geq 0,$
où Y est  le processus continu \`a droite d\'efini par:
\begin{equation}\label{Y}
     Y_t  =x + \displaystyle\sum_{n \geq 0} \,\Delta_n \,\1_{[\tau_n,\tau_{n+1}[}(t)+
     \int_0^t\left(b(\xi_s,Y_s)\,ds + \sigma(\xi_s,Y_s)\,dW_s\right) + Y_{\tau^-} \1_{[
     \tau,+\infty[}(t),
\end{equation}
où   $\Delta_n,$   taille du
saut du log de la valeur de la firme \`a l'instant $\tau_n$,  est
une variable al\'eatoire r\'eelle $\G_{\tau_n}\mbox{-mesurable}$ et $b:{U}\times \mathbb{R}\rightarrow \mathbb{R}$ et $\sigma:
{U}\times \mathbb{R}\rightarrow \mathbb{R}^{+}$  sont deux fonctions   
mesurables satisfaisant les conditions  de Lipschitz et    de croissance sous-linéaire:\\
\\
  - Il existe une constante   $  K \geq 0$ tels que pour tout $i  \in {U}$ et tout  $  x,y \in \mathbb{R}$: 
\begin{equation}\label{lip}
\big|b(i,x) -  b(i,y) \big|  + \big|\sigma(i,x) -  \sigma(i,y) \big|\leq K \big|x-y\big|. 
\end{equation} 
 -  Il existe une constante   $  K \geq 0$ tels que pour tout $i  \in {U}$ et tout  $  x  \in \mathbb{R}$: 
\begin{equation}\label{croissance}
\big|b(i,x)\big|^2  + \big|\sigma(i,x)   \big|^2 \leq K^2 (1+ |x |^2). 
\end{equation} 
  
  \begin{pro}\label{solforte} Sous les conditions (\ref{lip}) et (\ref{croissance}), il existe
 un unique 
 processus    $\G$-adapté  $Y$ solution forte de l'EDS (\ref{Y}).
\end{pro}

\noindent \textbf{Preuve:} 
Pour tout  $(\omega,t) \in \R^+\times\Omega,$  
il existe  un processus $Y^0$ où l'EDS associée peut être écrite sous la forme suivante:
$$dY_t^0(\omega) =   b(\xi_0(\omega),Y_t^0(\omega))dt + \sigma(\xi_0(\omega),Y_t^0(\omega))dW_t,$$
 où $b$ et $\sigma$ vérifient les conditions (\ref{lip}) et (\ref{croissance}). Donc on peut appliquer  les  théorèmes  2.5 (pp. 287 ) et 2.9 (pp. 289 ) de Karatzas et Shreve \cite{karatzas}: il existe une unique solution forte $Y_t^0$ donnée par 
$$
Y_t^0(\omega) = Y_{0}^0(\omega)  + \int_0^t\left( b(\xi_0(\omega),Y_s^0(\omega))ds + \sigma(\xi_0(\omega),Y_s^0(\omega))dW_s\right).
$$
 De plus, le processus $Y_t^0$ est $\G_t$-mesurable (voir \cite[p. 58]{lamberton} et \cite{jeanbyor}).  
En particulier, pour tout $(\omega,t) \in \R^+\times\Omega $
 et  $t = \tau_0^-$, nous obtenons:
$$Y_{\tau_0^-}^0(\omega) = Y_{0}^0(\omega)  + \int_0^{\tau_0}\left( b(\xi_0(\omega),Y_s^0(\omega))ds + \sigma(\xi_0(\omega),Y_s^0(\omega))dW_s\right),$$
et nous définissons
$$Y_{\tau_0}^0(\omega) = Y_{\tau_0^-}^0(\omega) + \Delta_0.$$
La taille du saut $\Delta_0$ est $\G_{\tau_0}$-mesurable. Par suite, le processus $Y_{\tau_0}^0$ est aussi $\G_{\tau_0}$-mesurable.
\\

Ce processus $Y^0$ restreint à l'ensemble $\{(\omega,t):~ t\leq\tau_0(\omega)\}$ est bien  $\G$-adapté et il est  l'unique solution forte de l'EDS restreinte à $\{(\omega,t):~ t\leq\tau_0(\omega)\}$.
\\

Supposons que le processus $Y$  défini  ci-dessous sur $\{(\omega,t):~ t\leq\tau_{n}(\omega)\}$ est bien $\G$-adapté et qu'il est  l'unique solution forte de l'EDS restreinte à $\{(\omega,t):~ t\leq\tau_{n}(\omega)\}$:
\begin{equation}\label{Yt}
Y_{t}(\omega) =   \sum_{\tau_k\leq t} \Delta_k     + \int_0^t\left( b(\xi_s(\omega),Y_s(\omega))ds + \sigma(\xi_s(\omega),Y_s(\omega))dW_s\right).
\end{equation}
On considère alors sur  l'ensemble $\{(\omega,t); \tau_n (\omega) \leq t < \tau_{n+1} (\omega)\},$ l'unique solution forte  $Y^n$  de l'équation différentielle stochastique :
\begin{equation}\label{EDSYn}
\left\{%
\begin{array}{ll}
    dY_t^n = b(\zeta_{n+1},Y_t^n)dt + \sigma(\zeta_{n+1},Y_t^n)dW_t^n \\
    Y_0^n = Y_{\tau_n}, \\
\end{array}%
\right.
\end{equation}
où $W^n$ est un mouvement brownien indépendant de la filtration $(\F_{\tau_n})_{n\geq 0}$.
L'unicité du processus $Y^n$ vient du fait que  les fonctions $b$ et $\sigma$ vérifient les hypothèses assurant l'existence et l'unicité (voir \cite{karatzas} et \cite{protter}).  
 Du fait de l'hypothèse de récurrence (\ref{Yt}) cette solution est également $\G$-adaptée.
Définissons alors sur l'ensemble $\{(\omega,t); \tau_n (\omega) \leq t < \tau_{n+1} (\omega)\}$ le processus
$$Y_t(\omega) := Y^n_{t- \tau_n} (\omega),$$
où $Y^n$ est solution de l'EDS (\ref{EDSYn}), ce qui montre que sur cet ensemble: 
$$
Y_{t}(\omega) =  Y_{\tau_{n}}(\omega) +  \int_{\tau_{n}}^t\left( b(\xi_0(\omega),Y_s(\omega))ds + \sigma(\xi_0(\omega),Y_s(\omega))dW_s\right),
$$
 et que $Y$ ainsi construit est également $\G$-adapté.
\\
 D'où l'hypothèse de récurrence est vérifiée pour tout $n.$
\endproof

\begin{rem}
La  filtration $\G$ \'etant une sous-filtration de $\F,$ les
processus  $\xi$ et $Y$ sont  aussi $\F\mbox{-adapt\'es}$.
\end{rem}

\begin{pro} D\'esignons par $\big(\bigcap_{s> t}\G_s, t \geq 0 \big)$ la filtration
$(\G_t, t \geq 0)$ rendue continue \`a droite. Elle v\'erifie
l'\'egalit\'e suivante:
$$\bigcap_{s> t}\G_s = \F_t.$$
\end{pro}
\noindent \textbf{ Preuve:} 
D'une part, $\G_t $ étant une sous-tribu de $\F_t$,  $\G_t \subset
\F_t$. D'o\`{u} $\bigcap_{s> t}\G_s \subset \bigcap_{s> t}\F_s$. La
filtration $\F$ \'etant continue \`a droite,
$$\bigcap_{s> t}\F_s = \F_t.$$
Ainsi, $\bigcap_{s> t}\G_s \subset \F_t.$ D'autre part,
 $$\bigcap_{s> t}\G_s = \bigcap_{s> t}\vee_{u<s}\F_u.$$
 Nous avons, $\forall s > t$, $\G_s = \vee_{u<s}\F_u$. Par
cons\'equent,
$$ \F_t \subset \vee_{u<s}\F_u = \G_s, \; \forall s >t.$$
Ce qui entra\^{\i}ne l'inclusion inverse et par suite l'\'egalit\'e
$\bigcap_{s> t}\G_s = \F_t.$
\endproof

\begin{defn} 
\label{defnstrategie}
 Un  contr\^{o}le impulsionnel 
  est la donn\'ee
d'une suite $\alpha$:
$$\alpha = (\tau_n, \zeta_{n+1}, \Delta_n, \, n \geq -1),$$
o\`{u} $(\tau_n)_{n\geq -1}$ est une suite  croissante de $\G\mbox{-temps
d'arr\^{e}t}$ de limite not\'ee $\tau$ et
v\'erifiant $\tau_{-1} = 0$,
$\zeta_{n+1}$ est la technologie choisie \`a l'instant $\tau_n$,
v.a. $\G_{\tau_n}\mbox{-mesurable}$ et $\Delta_n= Y_{\tau_n}-
Y_{\tau_n^-}$ est la taille du saut \`a l'instant $\tau_n$ et vérifiant $\Delta_{-1} =0$, v.a.
$\G_{\tau_n}\mbox{-mesurable}$, de telle sorte que la loi  du couple $(\zeta_{n+1},\Delta_n)$   sachant $\G_{\tau_n}$ ne dépend que de l'état du système en $\tau_n^-$ et que la probabilité de transition  $r$ 
définie  sur ${U}\times \mathbb{R}$ par
\begin{equation}\label{rn}
  \mathbb{P} (\zeta_{n+1}= j, Y_{\tau_n} = x +dy \,\big|\,\zeta_n = i, Y_{\tau_n^-} = x )= r (i,x;j,dy),
  \end{equation}
est ind\'ependante de $n$. 
\end{defn}

 On introduit l'ensemble $M$   de ces  noyaux  de transition 
sur $({U}\times
\mathbb{R}, \mathcal{P}({U})\otimes \mathcal{B}_{\R} )$ satisfaisant  $\forall (i,x), ~\delta_{(i,x)} \in M_{(i,x)}$    où $\delta_{(i,x)}$ est la mesure de Dirac en  $(i,x)$ et  $M_{(i,x)}$ la section en $(i,x)$  de $M$ définie comme suit: 
\begin{equation}\label{Maison}
   M_{(i,x)}= \big \{r  (i,x;.,.),  r\in M \}. 
\end{equation}
On   r\'eduit l'ensemble  des contrôles  en
supposant que  les lois de passage $r $ sont ind\'ependantes
de n. Ainsi, la famille des lois markoviennes est stationnaire.\\
Pour tout $(i,x),$ l'ensemble $M_{(i,x)}$ est muni de la topologie faible :  la suite $(r_n(i,x;.))_n$
converge vers $r(i,x;.)$ si et seulement si pour toute fonction continue bornée $g$ sur $ U \times \mathbb{R},$
\begin{equation}\label{topfaible}
    (r_n(i,x;g))_n \;\mbox{converge   vers} ~ r(i,x;g) ~\mbox{sur}~ \R.
\end{equation}

\begin{rem}\label{remsep}
1.  L'ensemble $M_{(i,x)}$ est métrisable pour la topologie faible définie ci-dessus (cf. Appendix 3 dans \cite{Billin}).\\
2.  $M_{(i,x)}$ étant inclus dans l'ensemble des probabilités sur
${U}\times \mathbb{R}$ et ${U}\times \mathbb{R}$ étant séparable, $M_{(i,x)}$ est par suite séparable (voir le  théorème 2.19 dans \cite[p. 25]{bain}).
\end{rem}

\noindent \textbf{Hypothèse 1:} \label{hypothese} Pour tout couple $(i,x) \in U \times \R$, l'ensemble $M_{(i,x)}$
est   ferm\'e  compact (donc complet) pour la topologie faible.

\begin{exmp}
Nous pouvons  prendre comme exemple l'ensemble  $M$ défini  par ses sections:
$$
   M_{(i,x)}= \left\{ r(i,x;j,dy) = p_{i,j}\otimes\frac{1}{\sqrt{2\Pi \sigma^2}} e^{\frac{-(y-x-m)^2}{2\sigma^2}}dy\cup\{\delta_{(i,x)}\}; 
 \sum_{j \in U} p_{i,j} =1
,m \in [\underline{m},\overline{m}]
    \right\}\cup\{\delta_{(i,x)}\}.$$
Soient  une fonction $g$ borélienne bornée sur ${U}\times
\mathbb{R}$ et une suite $(r_n) \in M_{(i,x)}$:
$$r_n(i,x;g) = \sum_{j\in U}p_{i,j }^n \int_{\mathbb{R}} g(j,y) \frac{1}{\sqrt{2\Pi \sigma^2}} e^{\frac{-(y-x-m_n)^2}{2\sigma^2}}\;dy,$$
où $m_n\in [\underline{m},\overline{m}]$. Il existe une sous-suite extraite $(m_{n_k})$ convergeant vers $m$ dans $\mathbb{R}$. Ensuite, la sous-suite $p_{i,j}^{n_k}$  est à valeurs dans le compact
$\{\sum_{j \in U} x_j =1\}$ et donc il existe une sous-sous suite extraite $p_{i,j}^{n_k'}$ convergente.  
 D'après le théorème de Lebesgue   de convergence majorée, la suite  $(r_{n_k'}(i,x;g))$ converge  pour toute fonction $g$
 et donc la suite $(r_{n_k'}(i,x;.))$ converge faiblement dans $ M_{(i,x)}$.\\
 Ainsi, l'ensemble $M_{(i,x)}$ est    fermé compact  pour la topologie faible. 
\end{exmp}

\begin{defn}
A chaque contrôle  $\alpha$ nous associons le gain
\begin{equation}
\label{k}
k(\alpha)=\int_0^ {+\infty} 
   e^{-\beta s}f(\xi_s,Y_s)\;ds -\displaystyle\sum_{ 0<\tau_n<\tau } e^{-\beta \tau_n}
c(\zeta_n,Y_{\tau_n^-},
\zeta_{n+1},Y_{\tau_n}),
\end{equation}
o\`{u} $\beta$ est un coefficient d'actualisation, 
$f:{U}\times \mathbb{R}\rightarrow \mathbb{R}^+$ (profit net de la firme)
et $c:{U}\times
\mathbb{R}\times{U}\times
\mathbb{R}\rightarrow \mathbb{R}^+$ (coût de changement de technologie) sont deux fonctions
bor\'eliennes positives.  La fonction $c$ 
vérifie pour tout couple  $(i,x) \in {U} \times
\mathbb{R}$, $
c(i,x,i,x) = 0,$
c'est à dire s'il n'y a pas d'impulsion, le coût est nul.
Soit le
gain moyen:
\begin{equation}\label{K}
   K(\alpha,i,x) = \mathbb{E}(k(\alpha)\big|\xi_0=i,Y_0=x).
\end{equation}
 \end{defn}

\begin{rem}\label{remN}
En utilisant le $\G$-temps d'arrêt $\tau_{N}$ (Proposition \ref{Nom}), le gain défini par (\ref{k}) peut être  écrit sous la forme suivante:
{\footnotesize
\begin{equation}\label{kN}
  k(\alpha)(\omega) =  \int_0^{\tau_{N(\omega)}} e^{-\beta s} f(\xi_s,Y_s)ds -\sum_{0<\tau_k < \tau_{N(\omega)}}e^{-\beta
\tau_k}c(\zeta_k,Y_{\tau_k^-},\zeta_{k+1},Y_{\tau_k}) +
\int_ {\tau_{N(\omega)}}^{+\infty} e^{-\beta s} f(\zeta_{N(\omega)},Y_s)ds.
\end{equation}}
\end{rem}

\begin{defn}\label{strategadm}
 La stratégie $\alpha = (\tau_n,\zeta_{n+1}, \Delta_n)_{n\geq -1} $   est dite admissible si et seulement si ces deux propriétés    sont vérifiées:
 $$1.
\left\{%
\begin{array}{ll}
   \int_0^ {\infty} 
   e^{-\beta s}f(\xi_s,Y_s) ds \in L^1\\
   \\
\sum_{ 0<  \tau_n<\tau} e^{-\beta \tau_n}
c(\zeta_n,Y_{\tau_n^-},
\zeta_{n+1},Y_{\tau_n})\in L^1.
\end{array}%
\right.
$$
2. Pour tout $n$, la loi de $(\tau_{n+1} - \tau_n,(\zeta_{n+1},Y_{\tau_n+ .})\1_{[\tau_n,\tau_{n+1}[})$ sachant $\F_{\tau_n}$ est égale à la loi de $(\tau_{n+1} - \tau_n,(\zeta_{n+1},Y_{\tau_n+ .})\1_{[\tau_n,\tau_{n+1}[})$ sachant $(\zeta_{n+1},Y_{\tau_n}),$  elle-même égale à la $\mathbb{P}_{(i,x)}$-loi de 
$(\tau_0,(i,Y_.)\1_{[0,\tau_0)})$ prise en  $ i = \zeta_{n+1} , x= Y_{\tau_n}.$ \\
On note $\A$ l'ensemble des strat\'egies
admissibles. 
\end{defn}
Nous remarquons que de fait la stratégie $\alpha$ est complètement définie par le couple $\alpha =(\tau_0, r),$ puisque la loi de $(\tau_0,(i,Y^x_.)\1_{[0,\tau_0)})$
 est déterminée par la diffusion associée à $(b(i,x), \sigma(i,x))$.

\begin{exmp}
\label{exmphyp2}
 Nous supposons que $\{0,1\}$ est l'ensemble des technologies permises telles que $0$ est l'ancienne technologie et $1$ est la nouvelle technologie.  On réduit dans cet exemple l'ensemble ${\cal A}$ en supposant que $\tau_0$ déterministe. \\
  Le profit net de la firme est donné pour tout  $(i,x) \in \{0,1\} \times  \mathbb{R}$ par $f(i,x) = \eta e^x$ 
et le coût de changement de technologie est égal à $c(i,x,1-i,y) = \exp(  x + \mu(y-x))$.  Nous supposons que $b(i,x) = b_i$ et $\sigma(i,x) = \sigma_i$, $\forall x$.\\
 L'ensemble $M$ des probabilités de  transition
est défini par  ses sections :
$$
   M_{(i,x)}= \left\{ r(i,x;j,dy) =    p_{i,j}\otimes \mathcal{N}(x+m,1); [p_{i,j}] = \left(
\begin{array}{ccc}
0 & 1 \\
1 & 0
\end{array}
\right),
m \in [\underline{m},\overline{m}]
    \right\}   \cup\{\delta_{i,x}\}.$$
Sous les hypothèses  $\beta > sup_{i \in U}( b_i + \frac{1}{2}\sigma_i^2)$ et $m < \frac{t_0}{2}\left(2 \beta - (b_0 +b_1) - \frac{1}{2}(\sigma^2_0 +\sigma^2_1)\right) - \frac{1}{2}$,
$$
\eta  \int_0^{+\infty}  e^{-\beta s} e^{Y_s} ds \in L^1. 
$$
De plus, $\tau_0,$ $ b_i$ et $ \sigma_i$ étant déterministes, nous obtenons
$$ \sum_{n\geq 0} e^{-\beta \tau_n}
\exp(  Y_{\tau_n^-} + \mu(Y_{\tau_n}-Y_{\tau_n^-})) \in L^1.
$$
  La preuve de ces deux propriétés   sera explicitée dans le chapitre 2 (en utilisant la proposition \ref{recurence}). Grâce à ces deux propriétés, nous pouvons déterminer le profit optimal de la firme et la stratégie optimale associée.
\end{exmp}

\begin{defn}
 Soient $\theta$ un $\G$-temps d'arr\^{e}t et $\alpha \in \A$ une strat\'egie admissible.
  Nous dirons qu'une
strat\'egie admissible $\mu$ se comporte comme $\alpha$ jusqu'\`a
$\theta$ exclu, que l'on note
\begin{equation}\label{deuxstratg}
    \{\mu_t =
\alpha_t, \, \forall\, t < \theta\}
\end{equation}
d\`es que les strat\'egies $\alpha$ et $\mu$ arr\^ et\'ees en
$\theta^-$ co\"{\i}ncident. C'est \`a dire que pour $\omega$ fix\'e,
il existe $n(\omega)$ tel que $\tau_n^\alpha (\omega)\leq \theta
(\omega) < \tau_{n+1}^\alpha (\omega)$ et nous avons, $\forall k
\leq n(\omega)$:
$$
\left\{%
\begin{array}{ll}
   \tau_k^\mu = \tau_k^\alpha,
   \\
   \zeta_{k}^\mu = \zeta_{k}^\alpha,
   \\
   r^\mu = r^\alpha =  r.
\end{array}%
\right.
$$
De m\^{e}me, nous dirons qu'une strat\'egie admissible $\mu$ se
comporte comme $\alpha$ jusqu'\`a $\theta$ inclus, que l'on note:
 $$ \{\mu_t =
\alpha_t, \, \forall\, t \leq \theta\}$$ d\`es que les strat\'egies
$\alpha$ et $\mu$ arr\^ et\'ees en $\theta$ co\"{\i}ncident. C'est à
dire que pour $\omega$ fix\'e, il existe $n(\omega)$ tel que
$\tau_n^\alpha (\omega)\leq \theta (\omega) < \tau_{n+1}^\alpha
(\omega)$, et par suite, $\forall k \leq n (\omega)$:
$$
\left\{%
\begin{array}{ll}
   \tau_k^\mu = \tau_k^\alpha,
   \\
   \zeta_{k+1}^\mu = \zeta_{k+1}^\alpha,
   \\
   r^\mu = r^\alpha=  r.
\end{array}%
\right.
$$

Le cas particulier $\theta = \tau_n^\alpha$ les strat\'egies
$\alpha$ et $\mu$ arr\^ et\'ees en $\tau_n^-$ se traduit par le fait
que  $\forall\,\omega$,  $n (\omega) = n$, et $\forall k \leq n$:
$$
\left\{%
\begin{array}{ll}
   \tau_k^\mu = \tau_k^\alpha,
   \\
   \zeta_{k}^\mu = \zeta_{k}^\alpha,
   \\
   r^\mu = r^\alpha=  r.
\end{array}%
\right.
$$

Le cas particulier $\theta = \tau_n^\alpha$ les strat\'egies
$\alpha$ et $\mu$ arr\^ et\'ees en $\tau_n$ se traduit par le fait
que $\forall\,\omega$, $n (\omega) = n$, et $\forall k \leq n $:
$$
\left\{%
\begin{array}{ll}
   \tau_k^\mu = \tau_k^\alpha,
   \\
   \zeta_{k+1}^\mu = \zeta_{k+1}^\alpha,
   \\
   r^\mu = r^\alpha  =  r.
\end{array}%
\right.
$$
\end{defn}

\section{Crit\`eres d'optimalit\'e}\label{section 3}
Le probl\`eme d'optimalit\'e pos\'e consiste \`a prouver l'existence
d'une strat\'egie admissible  $\widehat{\alpha}$ qui maximise la
fonction gain  $ K(\alpha,i,x)$ d\'efinie par l'expression (\ref{K}),
c'est \`a dire trouver une strat\'egie $\widehat{\alpha}$ telle que
\begin{equation}\label{optimal}
    K(\widehat{\alpha},i,x) = \ESS {\alpha \in \A} K(\alpha,i,x).
\end{equation}
La strat\'egie $\widehat{\alpha}$ est dite optimale.

\subsection{Gains maximaux conditionnels}
Nous introduisons, tout d'abord, deux notions que nous utiliserons
fr\'equemment (voir  \cite[p. 87-92]{Karoui}):
\begin{defn} \label{defnsysteme} Soit une filtration $\h$ et $\mathcal{T}$ une sous-famille de $\h$-temps d'arr\^et. Une famille de variables al\'eatoires
$\{X^\alpha_\theta, \theta \in\mathcal{T} , \alpha \in
\A\} $ est appel\'ee un
$(\h,\mathcal{T} ,\A)\mbox{-syst\`eme}$ si, pour tout
$ \theta \in\mathcal{T} $ et tout $\alpha \in \A$,
l'ensemble des strat\'egies admissibles, nous avons:
\begin{itemize}
    \item [i/] $\forall\, \gamma \in \mathcal{T} ,$
    sur l'ensemble $\{\theta = \gamma\}$, $X^\alpha_\theta=
    X^\alpha_\gamma \quad \mathbb{P}\;\mbox{p.s.}$
    \item [ii/] Les variables aléatoires  $X^\alpha_\theta$\, sont ${\h}_{\theta}$-mesurables.
 \item [iii/] Si $\mu \in \A$, $\mu_t = \alpha_t ~ \forall t< \theta$ et $\mu = \alpha$ sur $D  \in {\h}_{\theta}$ (i.e.  $ ~ \forall \omega\in D$, $\mu(t,\omega) = \alpha(t,\omega) $  pour tout $t$) alors
  $X^\alpha_\theta=
    X^\mu_\theta$    sur $D$        $\mathbb{P}\;\mbox{p.s.}$ 
\end{itemize}
\end{defn}

\begin{defn}
Un
$(\h,\mathcal{T},\A)$sur-martingale-syst\`eme 
(resp. martingale-syst\`eme, sous-martingale-syst\`eme) est un
$(\h,\mathcal{T},\A)$-syst\`eme tel que:
\begin{itemize}
    \item [i/] Pour tout $\theta \in\mathcal{T} $ et  tout $\alpha \in \A$, $X^\alpha_\theta$ est
    $\mathbb{P}\mbox{-int\'egrable}$.
    \item [ii/] Si $\gamma\;et\;\theta$ sont deux \'el\'ements de
    $\mathcal{T} $ tels que $\gamma \leq  \theta$, alors:
    $$\mathbb{E}\left(X^\alpha_\theta \big| \h_\gamma\right) \leq X^\alpha_\gamma \qquad \mbox{p.s. (resp. } \; = , \geq).$$
    \item [iii/]  Si $\mu \in \A$, $\mu_t = \alpha_t ~ \forall t< \theta$, alors
  $X^\alpha_\theta=
    X^\mu_\theta \quad  \mathbb{P}\;\mbox{p.s.}$ 
 \end{itemize}
\end{defn}

La m\'ethode de r\'esolution des probl\`emes de contr\^{o}le
stochastique de type (\ref{optimal}) repose sur le principe de
Bellman (voir \cite[p. 95]{Karoui}).
Plus pr\'ecis\'ement, si on conna\^{\i}t une stratégie optimale
$\widehat{\alpha}$
 jusqu'à un temps d'observation $T$ , et une autre optimale $\tilde{\alpha}$ de $T$
 \`a $T+h$, il reste optimal entre $0$ et $T+h$
 de  garder $\widehat{\alpha}$ jusqu'\`a $T$  et de la prolonger apr\`es par
 $\tilde{\alpha}$. Ce principe est en fait un critère nécessairement et suffisant d'optimalité. C'est la raison pour laquelle  nous introduisons les
 gains maximaux
conditionnels  suivants:
\begin{defn} Soit $\mathcal{T} $  la famille des $\G$-temps d'arr\^{e}t
et $\mathcal{T}^* $ la famille des $\G$-temps d'arr\^{e}t   strictement positifs. 
Pour toute
strat\'egie $\alpha \in \A$, nous appelons gain
maximal conditionnel  la famille d\'efinie  par
$$F^\alpha_\theta= 
 \ESS{\mu_t=\alpha_t,t<\theta} \{\mathbb{E}[k(\mu)\,|\G_\theta]\},~
\theta \in\mathcal{T}^*.$$ 
 Respectivement, 
$$ F^{\alpha^+}_{\theta}= \ESS{\mu_t=\alpha_t,t\leq\theta} \{ \mathbb{E}[k(\mu)\,|\F_\theta]\} , \theta\in
\mathcal{T}.$$
\end{defn}

 La variable al\'eatoire $F_{\theta}^{\alpha}$ (resp.
$F_{\theta}^{\alpha ^ +}$) repr\'esente le gain maximal conditionnel
connaissant l'\'evolution du processus impuls\'e jusqu'\`a l'instant
$\theta$ en excluant (resp. incluant) la renaissance \'eventuelle en
cet instant.\\
\\
Pour toute stratégie $\alpha \in \A$ et tout temps d'arrêt $\theta \in \mathcal{T}^*,$ le gain apr\`es $\theta $ est   donn\'e p.s. par
l'expression suivante:
\begin{equation}\label{gaintheta}
k_\theta(\alpha) = \int_\theta ^{+\infty} e^{-\beta s}f(\xi_s,Y_s)ds -\displaystyle\sum_{\theta \leq \tau_n <\tau}e^{-\beta
\tau_n}
    c(\zeta_n,Y_{\tau_n^-},
    \zeta_{n+1},Y_{\tau_n}). 
    \end{equation}
     Respectivement, le gain strictement  apr\`es $\theta$ est donn\'e p.s. par
\begin{equation}
\label{gaintheta+}
   k_{\theta^+} (\alpha) = \int_\theta ^{+\infty} e^{-\beta s}f(\xi_s,Y_s ) ds -
\displaystyle\sum_{\theta < \tau_n  < \tau } e^{-\beta
\tau_n }
    c(\zeta_n ,Y_{ \tau_n ^-}, \zeta_{n+1} ,Y_{\tau_n }).
\end{equation}

D'où la remarque suivante:
\begin{rem} \label{remF0+}
1.  Les égalités (\ref{k})  et (\ref{gaintheta+}) montrent  que 
$ k_{0^+}(\alpha) =  k(\alpha)$. On en déduit que,  pour toute stratégie $\alpha \in \A$, 
$$F^{\alpha^+}_0=  \sup_{\mu\in
    \A}\mathbb{E}[k(\mu)].$$ 
2. L'expression $k(\alpha) - k_{\theta}(\alpha)$
$\left(\mbox{resp.}\,k(\alpha) - k_{\theta^+}(\alpha)\right)$
\'etant $\G_\theta $ (resp. $ \F_\theta $)-mesurable, nous avons $\forall \theta \in \mathcal{T}^* $ $\forall \alpha \in \A,$
    $$ F_{\theta}^{\alpha} = k(\alpha) - k_\theta(\alpha)
    +\displaystyle \ESS{\{\mu_t = \alpha_t,\,t <\theta\}}\{
\mathbb{E}\left(k_\theta(\mu)|\,\G_{\theta}\right)\}.$$
 Respectivement, $\forall \theta \in \mathcal{T}$
$$ F_{\theta}^{\alpha^+} = k(\alpha) - k_{\theta^+}(\alpha)
+ \displaystyle \ESS{\{\mu_t = \alpha_t,\,t \leq \theta\}}\{
\mathbb{E} \left(k_{\theta^+}(\mu)|\,\F_\theta\right)\}~p.s.
$$

\end{rem}
\begin{pro} La famille de  gains maximaux conditionnels $(F^\alpha_\theta,~\alpha \in \A,\theta~\in \mathcal{T}^*)$ est un
$(\G,\mathcal{T}^* ,\A)$-syst\`eme. Respectivement,
$(F^{\alpha^+}_\theta,~\alpha \in \A,\theta~\in \mathcal{T}) $ est un
$(\F,\mathcal{T} ,\A)$-syst\`eme.
\end{pro}

\noindent\textbf{Preuve:} 
Il s'agit de montrer que la  famille de  variables aléatoires $(F^\alpha_\theta)$ v\'erifie les propri\'et\'es de la d\'efinition \ref{defnsysteme}. \\
\\
1. $\forall\,\gamma \in\mathcal{T}^* $, nous avons:
$$
F_{\gamma}^{\alpha} = k(\alpha) - k_\gamma(\alpha)
    +\displaystyle \ESS{\{\mu_t = \alpha_t,\,t <\gamma\}}
\mathbb{E}\left(k_\gamma(\mu)|\,\G_{\gamma}\right).
$$
Sur l'ensemble $\{\theta = \gamma\},$ nous avons $k_\theta(\mu)
=k_\gamma(\mu)$. Ensuite, nous sommes ramen\'es \`a montrer
l'\'egalit\'e
$$\1_{\{\theta = \gamma\}}\,\mathbb{E}\left(k_\gamma(\mu)|\,\G_{\theta}\right) =
\1_{\{\theta =
\gamma\}}\,\mathbb{E}\left(k_\gamma(\mu)|\,\G_{\gamma}\right).$$ De
prime abord, $\{\theta = \gamma\} \in \G_\gamma \cap \G_\theta.$ 
D'une part, l'expression $$\1_{\{\theta =
\gamma\}}\,\mathbb{E}\left(k_\gamma(\mu)|\,\G_{\gamma}\right)\,\1_{\{\gamma
\leq t\}}$$ est $\G_t$-mesurable. Elle est \'egale \`a l'expression
$\1_{\{\theta =
\gamma\}}\,\mathbb{E}\left(k_\gamma(\mu)|\,\G_{\gamma}\right)\,\1_{\{\theta
\leq t\}}$ qui est $\G_t$-mesurable et  donc l'expression
$\1_{\{\theta=\gamma\}}\,\mathbb{E}\left(k_\gamma(\mu)|\,\G_{\gamma}\right)\}$
est
$\G_\theta$-mesurable.\\
D'autre part, l'expression $\1_{\{\theta =
\gamma\}}\,\mathbb{E}\left(k_\gamma(\mu)|\,\G_{\theta}\right)$ est $\G_\theta$-mesurable.\\
Puis, soit $X$ une v.a. $\G_\theta$-mesurable. Nous avons $X
\,\1_{\{\theta = \gamma\}} \in \G_{\gamma}$. Par suite,
$$\mathbb{E}\left(X \,\1_{\{\theta = \gamma\}}  \,\mathbb{E}\left[ k_\gamma(\mu)|\,\G_{\gamma}\right]\right)=
\mathbb{E}\left(X \,\1_{\{\theta = \gamma\}}  \,
k_\gamma(\mu)\right).$$ De plus,
$$\mathbb{E}\left(X \,\1_{\{\theta = \gamma\}}  \,
k_\gamma(\mu)\right) = \mathbb{E}\left(X \,\1_{\{\theta = \gamma\}}
\,\mathbb{E}\left[ k_\gamma(\mu)|\,\G_{\theta}\right]\right).$$ Il
en suit,
$$\mathbb{E}\left(X \,\1_{\{\theta = \gamma\}}  \,\mathbb{E}\left[
k_\gamma(\mu)|\,\G_{\gamma}\right]\right)=\mathbb{E}\left(X
\,\1_{\{\theta = \gamma\}}  \,\mathbb{E}\left[
k_\gamma(\mu)|\,\G_{\theta}\right]\right).$$
 Ainsi,
$$F^\alpha_\theta =F^\alpha_\gamma \quad \mathbb{P}~ \mbox{p.s.  sur l'ensemble } \{\theta = \gamma\}. $$
2.  L'expression $k(\alpha) - k_\theta(\alpha)$ est
$\G_\theta$-mesurable. De plus, par d\'efinition de l'ess-sup d'une
famille mesurable, l'expression $ \ESS{\{\mu_t =
\alpha_t,\,t <\gamma\}}
\mathbb{E}\left(k_\gamma(\mu)|\,\G_{\gamma}\right)$ est aussi
$\G_\theta$-mesurable. Ainsi, $F^\alpha_{\theta}$ est
$\G_\theta$-mesurable. \\
\\
3. 
 Si $\mu \in \A$, $\mu_t = \alpha_t ~ \forall t< \theta$ et $\mu = \alpha$ sur $D  \in {\G}_{\theta}$ 
 nous avons:
 $$
F_{\theta}^{\mu} = (k(\mu) - k_\theta(\mu))
    +\displaystyle \ESS{\{\mu_t = \alpha_t,\,t <\theta\}}
\mathbb{E}\left(k_\theta(\alpha)|\,\G_{\theta}\right).
$$
Ainsi, sur $D$, nous obtenons $F^\alpha_\theta =F^\mu_\theta \quad\mathbb{P}\; \mbox{p.s.} $
 $F_\theta^\alpha$ est bien  un
$(\G,\mathcal{T}^* ,\A)$-syst\`eme.\\
\\
On proc\`ede d'une mani\`ere analogue pour montrer que
$(F_\theta^{\alpha^+})$ est   un
$(\F,\mathcal{T} ,\A)$-syst\`eme en constatant que
$F_\theta^{\alpha^+}$ est $\F_\theta$-mesurable puisque c'est la
somme de deux expressions $\F_\theta$-mesurables.
\endproof

\begin{defn}
\label{cond opt}  Pour toute stratégie $\alpha \in \A$ et tout $\G$-temps d'arr\^{e}t $\theta \in \mathcal{T},$ une strat\'egie $\widehat{\alpha}$ est dite
$(\theta,\alpha)\mbox{-conditionnellement optimale}$ si elle est admissible, $\widehat{\alpha}_t = \alpha_t \; \forall t\leq \theta$ et elle v\'erifie  p.s.
$$ F_{\theta}^{\widehat{\alpha}^+} =
\mathbb{E}(k(\widehat{\alpha})|\,\F_{\theta}) \quad\mathbb{P}\; \mbox{p.s.} $$
\end{defn}

 Le principe d'optimalit\'e est bas\'e
sur la propri\'et\'e de commutation de l'ess-sup et de l'esp\'erance
conditionnelle, qui est en particulier satisfaite si la famille
consid\'er\'ee est un ensemble filtrant croissant (voir 
\cite[p. 230]{Karoui}). D'où le lemme suivant:

\begin{lem}\label{ens decr}
Pour tout $\G\mbox{-temps d'arr\^{e}t } ~   \theta \in \mathcal{T}$  et toute
strat\'egie $\alpha \in \A$, l'ensemble
$$\{\mathbb{E}(k(\mu) | \G_\theta ); \forall t
< \theta, \, \mu_t = \alpha_t\}$$ est un  ensemble filtrant
croissant. Egalement on a la m\^{e}me propri\'et\'e pour les
ensembles suivants:
$$ \{\mathbb{E}(k_\theta(\mu) |
\G_\theta ); \forall t < \theta, \, \mu_t = \alpha_t\}\quad
et\quad \{\mathbb{E}(k_{\theta^+}(\mu) | \F_\theta);
\forall t \leq \theta, \, \mu_t = \alpha_t\}.$$
\end{lem}

\noindent \textbf{Preuve:} 
 Consid\'erons deux strat\'egies admissibles
$\mu^1$ et $\mu^2$  v\'erifiant
$$\{\mu^1_t = \alpha_t, t < \theta\} \quad et \quad \{\mu^2_t = \alpha_t, t < \theta\}.$$
Introduisons, pour simplifier l'\'ecriture, les variables
al\'eatoires $\G_\theta\mbox{-mesurables}$ suivantes:
$$
F^1 = \mathbb{E}(k(\mu^1)|\,\G_{\theta})\quad et \quad
 F^2 =
\mathbb{E}(k(\mu^2)|\,\G_{\theta}).
$$
Nous définissons  la strat\'egie  admissible $\mu =
(\tau_n,\zeta_{n+1},\Delta_n)$ par
$$\mu = \left\{%
\begin{array}{ll}
    \mu^2& \hbox{sur} \; \{F^1 \leq F^2\}\\
    \mu^1 & \hbox{sur} \;\{F^1 > F^2\} \\
\end{array}%
\right.   $$ Elle v\'erifie $\{\mu_t = \alpha_t, t < \theta\}$.
\begin{eqnarray*}
\mathbb{E}\left(k(\mu)|\,\G_{\theta}\right) & = &
\mathbb{E}\left[1_{\{F^1 \leq F^2\}}k(\mu) |\,\G_{\theta}\right] +
\mathbb{E}\left[1_{\{F^1 > F^2\}}k (\mu)|\,\G_{\theta}\right] \nonumber\\
&=& \mathbb{E}\left(k (\mu^2)|\,\G_{\theta}\right)\; 1_{\{F^1 \leq
F^2\}} + \mathbb{E}\left(k (\mu^1)|\,\G_{\theta}\right) \;1_{\{F^1
> F^2\}}  \nonumber\\
&=&  F^1 \vee F^2.
\end{eqnarray*}
Pour l'ensemble $ \{\mathbb{E}(k_\theta(\mu) | \G_\theta)~;~
\forall t < \theta,~\mu_t = \alpha_t\}$, au lieu des v.a. $F^1$
et $F^2$, il suffit de prendre les variables al\'eatoires
$\G_\theta\mbox{-mesurables}$ suivantes:
$$
F^1_1 = \mathbb{E}(k_\theta(\mu^1)|\,\G_{\theta})\quad et \quad
 F^2_2 =
\mathbb{E}(k_\theta(\mu^2)|\,\G_{\theta}),
$$
et la d\'emarche sera la m\^{e}me que pr\'ec\'edemment. Tandis que
pour l'ensemble $\{\mathbb{E}(k_{\theta^+}(\mu) | \F_\theta
)~;~
\forall t \leq \theta,~\mu_t = \alpha_t\}$, nous
proc\`edons d'une mani\`ere analogue en consid\'erant les
strat\'egies $\mu^3$ et $\mu^4$ v\'erifiant
$$\{\mu^3_t = \alpha_t, t \leq \theta\} \quad et \quad \{\mu^4_t = \alpha_t, t \leq \theta\},$$
et  les v.a. $\F_\theta\mbox{-mesurables}$ suivantes:
$$
F^3 = \mathbb{E}(k_{\theta^+}(\mu^3)|\,\F_{\theta})\quad et \quad
 F^4 =
\mathbb{E}(k_{\theta^+}(\mu^4)|\,\F_{\theta}).
$$
\endproof

\begin{cor}\label{commutation}
 Pour tout $\G\mbox{-temps d'arr\^{e}t} ~  \theta \in \mathcal{T}$  et toute
strat\'egie $\alpha \in \A$, l'ensemble
$$\{\mathbb{E}(k(\mu) | \G_\theta ); \forall t
< \theta, \, \mu_t = \alpha_t\}$$ étant  filtrant
croissant, nous pouvons commuter l'ess-sup et la  $\G \mbox{ (resp.}~\F)-$espérance conditionnelle. 
Il en est de même pour les ensembles:
$$ \{\mathbb{E}(k_\theta(\mu) |
\G_\theta ); \forall t < \theta, \, \mu_t = \alpha_t\}\quad
et\quad \{\mathbb{E}(k_{\theta^+}(\mu) | \F_\theta);
\forall t \leq \theta, \, \mu_t = \alpha_t\}.$$
\end{cor}
Pour plus de détails, voir la proposition A.$2$ de N. El Karoui \cite[p. 230]{Karoui}.

\begin{pro}
\label{sous-marting}
 Le gain maximal conditionnel $(F_{\theta}^{\alpha})$
(resp. $(F_{\theta}^{\alpha ^ +})$) forme un $ (\G,\mathcal{T}^*,
\A)$ $(\mbox{resp.}\,(\F,\mathcal{T} ,
\A))$-sur-martingale-syst\`eme positif.
\end{pro}

\noindent \textbf{Preuve:}\\ 
1. Soient $\theta$ et $\gamma$ deux
$\G\mbox{-temps d'arr\^{e}t}$ avec $\gamma \leq \theta$. Alors, p.s.
\begin{eqnarray*}
\mathbb{E}(F_{\theta}^{\alpha}|\,\G_{\gamma}) & = & \mathbb{E}\left[
 \ESS{\mu_t=\alpha_t,t<\theta}E\left(k(\mu)\,\big|\G_\theta\right)\big|\,\G_{\gamma}\right]\nonumber\\
&\leq &  \mathbb{E}\left[
 \ESS{\mu_t=\alpha_t,t<\gamma}E\left(k(\mu)\,\big|\G_\theta\right)\big|\,\G_{\gamma}\right],
\end{eqnarray*}
le sup \'etant pris sur un ensemble plus vaste  dans la deuxi\`eme
in\'egalit\'e. Nous pouvons commuter l'ess-sup et la
$\G_{\gamma}\mbox{-esp\'erance conditionnelle}$, gr\^{a}ce  au corollaire \ref{commutation}, et obtenir p.s.
$$\mathbb{E}(F_{\theta}^{\alpha}|\,\G_{\gamma}) \leq F_{\gamma}^{\alpha}.$$
2. On applique le m\^{e}me raisonnement  pour le gain maximal
conditionnel $F_{\theta}^{\alpha^+}$:
\begin{eqnarray*}
\mathbb{E}(F_{\theta}^{\alpha^+}|\,\F_{\gamma}) & = &
\mathbb{E}\left[
 \ESS{\mu_t=\alpha_t,t\leq\theta}E\left(k(\mu)\,\big|\F_\theta\right)\big|\,\F_{\gamma}\right]\nonumber\\
&\leq &  \mathbb{E}\left[
 \ESS{\mu_t=\alpha_t,t\leq\gamma}E\left(k(\mu)\,\big|\F_\theta\right)\big|\,\F_{\gamma}\right],\nonumber\\
\end{eqnarray*}
le sup \'etant pris sur un ensemble plus vaste  dans la deuxi\`eme
in\'egalit\'e. Nous pouvons commuter l'ess-sup et la
$\F_{\gamma}\mbox{-esp\'erance conditionnelle}$, gr\^{a}ce   au corollaire \ref{commutation}, et obtenir p.s.
$$\mathbb{E}(F_{\theta}^{\alpha^+}|\,\F_{\gamma}) \leq F_{\gamma}^{\alpha^+}.$$
Il reste à prouver que les variables aléatoires $F_{\theta}^{\alpha^+}$
et $F_{\theta}^{\alpha}$
sont positives pour tout  $\G$-temps d'arrêt $\theta$.  Il suffit de
choisir  une stratégie $\alpha$ appartenant à $\A$ dont le
gain est positif. Soit une stratégie $\alpha \in \A$ avec
$\xi_t = i,~\forall t$ et $\tau_0=+\infty$. Le
gain associé à $\alpha$ est donné par
$$
k(\alpha) = \int_0^{+\infty} e^{-\beta s} f(i,Y_s)\,ds.
$$
Or la fonction $f$ est positive et donc l'expression précédente est
positive. Ce qui entraîne que son $\F_{\theta}$ (respectivement $\G_{\theta}$)-espérance conditionnelle est aussi positive. En passant à l'essentiel
supremum, on déduit que $F_{\theta}^{\alpha^+}$ et $F_{\theta}^{\alpha}$ sont
positives.  
\endproof

Une cons\'equence imm\'ediate de cette proposition est le premier
crit\`ere d'optimalit\'e qui permet de réduire consid\'erablement la classe des
stratégies susceptibles d'être  optimales:
\begin{cor}
\label{martingale} Une condition n\'ecessaire et suffisante pour que
la strat\'egie $\widehat{\alpha}$ soit optimale est que le gain
maximal conditionnel $F_.^{\widehat{\alpha}^+}$ soit  un
$(\F,\mathcal{T} ,\A)$-martingale-syst\`eme, c'est à dire
$\forall \, \theta, \gamma$ deux $\G$-temps d'arrêt, $\gamma \leq
\theta,$ nous avons
$$\mathbb{E}(F_{\theta}^{\alpha^+}|\,\F_{\gamma}) =
F_{\gamma}^{\alpha^+} \quad \mathbb{P} \;\mbox{p.s}.$$
\end{cor}

\noindent \textbf{Preuve:}  La  strat\'egie admissible
$\widehat{\alpha}$ est  optimale, alors elle vérifie  pour tout $\theta\in \mathcal{T}$:
$$\mathbb{E}(k (\widehat{\alpha})) = \sup_{\alpha \in
\A}\mathbb{E}(k (\alpha)) \geq
\sup_{\{\alpha:\alpha_t=\widehat{\alpha}_t,\, t \leq \theta \}}
\mathbb{E}(k (\alpha))  \geq \mathbb{E}(k (\widehat{\alpha})).$$
D'o\`{u} l'\'egalit\'e. De plus, la commutation de l'ess-sup et
l'esp\'erance conditionnelle et l'\'egalit\'e pr\'ec\'edente
entra\^{\i}nent
\begin{equation}\label{pr}
    \mathbb{E} (F^{\widehat{\alpha}^+}_{\theta}) = \mathbb{E}\left (
    \displaystyle \ESS{\{\alpha_t = \widehat{\alpha}_t,\,t \leq\theta\}}
\mathbb{E}(k(\alpha)|\,\F_{\theta})\right)=\sup_{\alpha \in
\A}\mathbb{E}(k (\alpha))
=\mathbb{E}(k(\widehat{\alpha})) = F_0^{\widehat{\alpha}^+}.
\end{equation}
La derni\`ere \'egalit\'e provient de la d\'efinition \ref{cond opt}
appliqu\'ee au temps $0$  et montre que l'on a bien un $(\F,\mathcal{T} ,\A)$-martingale-syst\`eme.\\

Inversement, supposons que $(F^{\widehat{\alpha}^+}_{\theta})$ est
 un $(\F,\mathcal{T} ,\A)$-martingale-syst\`eme, c'est \`a dire
$$p.s.~~F^{\widehat{\alpha}^+}_{\gamma} =  \mathbb{E} (F^{\widehat{\alpha}^+}_{\theta}|\,\F_{\gamma})\quad \mbox{ pour tout } \gamma\leq\theta, \gamma,\theta\in \mathcal{T}.$$
Citons le th\'eor\`eme $1.17$ d'El Karoui \cite{Karoui} : \label{thm
karoui}\\
 "\textit{Une condition n\'ecessaire et suffisante pour
qu'un contr\^{o}le $\widehat{\alpha}$ soit optimal est que, pour
tout temps d'observation $\gamma$, il soit
$(\widehat{\alpha},\gamma)\mbox{-conditionnellement optimal}$, ou ce
qui est \'equivalent, que le gain maximal conditionnel par
rapport \`a $\widehat{\alpha}$, soit un
$(\F,\mathcal{T} , \A)$-martingale-syst\`eme, c'est \`a dire que 
si $\gamma$ et $\theta$ sont deux temps d'observation avec $\gamma
\leq \theta$: $$
\mathbb{E}[F_{\theta}^{\widehat{\alpha}^+}|\,\F_\gamma] =
F_{\gamma}^{\widehat{\alpha}^+}\quad \mbox{p.s}."$$}

On peut donc dire que la strat\'egie $\widehat{\alpha}$ est
conditionnellement optimale  donc optimale.
\endproof

Nous  introduisons une nouvelle notion de gains maximaux:

\begin{defn}
\label{coutcond} Pour toute strat\'egie admissible $\alpha$ et 
tout $\G\mbox{-temps d'arrêt} ~ \theta$, nous appellons gain maximal
conditionnel après $\theta  \in \mathcal{T}^*$  la variable aléatoire d\'efinie
 p.s. par:
$$
    W_{\theta}^{\alpha} :=   \ESS{\{\mu_t = \alpha _t,\; \forall \,t <\theta\}}
    \mathbb{E} \left[ k_\theta(\mu)
    \,|\G_{\theta}\right].
$$
Respectivement, le gain maximal
conditionnel strictement  après $\theta  \in \mathcal{T}$ est la variable aléatoire  définie p.s. par:
\begin{eqnarray*}
W_{\theta}^{\alpha^+} &:= & \displaystyle \ESS{\{\mu_t = \alpha
_t,\; \forall \,t \leq
    \theta\}}
    \mathbb{E} \left[k_{\theta^+}(\mu)
    \,|\F_\theta\right].
\end{eqnarray*}
\end{defn}

\begin{pro} Pour toute stratégie $\alpha \in \A$ et tout
$\G$-temps d'arrêt $\theta$, le gain maximal  conditionnel
$F^\alpha_\theta$ peut \^{e}tre \'ecrit, p.s., sous la forme
suivante:
\begin{equation}\label{rel FW}
    F^\alpha_\theta = \left(k (\alpha)- k_\theta(\alpha)\right) +
     W^\alpha_\theta,~ \mbox{ où } \theta>0~p.s.
\end{equation}
De même, le gain maximal  conditionnel $F^{\alpha^+}_\theta$ est
donn\'e par l'expression
\begin{equation}\label{rel FW+}
    F^{\alpha^+}_\theta = \left(k(\alpha) - k_{\theta^+}(\alpha)\right) +
     W^{\alpha^+}_\theta.
\end{equation}
\end{pro}

\noindent \textbf{Preuve: } 
De la d\'efinition du gain maximal  conditionnel, nous avons
  pour  $\theta\in\TT^*$:
$$F^\alpha_\theta=\left(k (\alpha)- k_\theta(\alpha)\right) +
\ESS{\{\mu_t = \alpha _t,\; \forall \,t < \theta\}}  \mathbb{E}
\left( k_\theta(\mu) |\, \G_\theta\right),$$ 
o\`u l'on reconnait
$W_\theta^\alpha$ dans le deuxi\`eme terme, soit (\ref{rel FW}).\\

De m\^{e}me  pour tout $\theta\in\TT$:
\begin{eqnarray*}
 F_{\theta}^{\alpha^+} &=& \displaystyle \ESS{\{\mu_t = \alpha
_t,\; \forall \,t \leq
    \theta\}}
    \mathbb{E} \left[k(\mu)
    \,|\F_\theta\right]\nonumber\\
&=&\left(k(\alpha) - k_{\theta^+}(\alpha)\right) + \displaystyle
\ESS{\{\mu_t = \alpha_t,\,t \leq \theta\}} \mathbb{E}
\left(k_{\theta^+}(\mu)|\,\F_\theta\right).
\end{eqnarray*}
o\`u l'on reconnait $W_\theta^{\alpha^+}$ dans le deuxi\`eme terme.
D'o\`{u} l'\'egalit\'e (\ref{rel FW+}) est v\'erifi\'ee.
\endproof

\begin{rem}
 On d\'eduit imm\'ediatement des
\'egalit\'es (\ref{rel FW}) et (\ref{rel FW+}) et du fait que les
gains maximaux conditionnels $F^\alpha_\theta$ et
$F^{\alpha+}_\theta$ sont des $(\G,\mathcal{T}^*,\A)$ (respectivement $(\F,\mathcal{T} ,\A)$)-syst\`emes que  $(W_{\theta}^{\alpha})$ d\'efinit un ($\G,
\mathcal{T} ^*,\A$)-syst\`eme et que
$(W_{\theta}^{\alpha^+})$ d\'efinit un ($\F,
\mathcal{T} ,\A$)-syst\`eme.
\end{rem}

\begin{lem}\label{admcout}
Le gain maximal conditionnel $W_{\tau_n}^{\alpha^+}$ converge vers
$\int_{\tau}^{+\infty}e^{-\beta s} f(\xi_s,Y_s) ds$ p.s.  lorsque $n$ tend vers l'infini.
\end{lem}

\noindent \textbf{Preuve}: 
De la proposition \ref{sous-marting}, la famille
$(F_{\theta}^{\alpha^+},~\theta\in\mathcal{T} )$ est un $(\F,
\mathcal{T} , \A)$-sur-martingale-syst\`eme positif,
donc
 $(F_{\tau_n}^{\alpha^+},~n\geq 0)$ est une
sur-martingale discrète positive pour la filtration $\F_{\tau_n}$. 
Gr\^{a}ce à la convergence des sur-martingales positives, il existe
une limite p.s. positive de
 $F_{\tau_n}^{\alpha^+}$, lorsque n tend vers l'infini, not\'ee
$F_{\infty}^{\alpha^+}$ et qui vérifie:
 \begin{equation}\label{surmartingale}
 F_{\tau_n}^{\alpha^+}\geq
\mathbb{E}\left[F_{\infty}^{\alpha^+}|\F_{\tau_n}\,\right]
\quad\mbox{p.s.}
\end{equation}
 Cette limite est $\F_{\tau}$-mesurable. En effet,
nous avons, $\forall\, B$ ouvert:
$$\{\omega \in \Omega:\,F_{\infty}^{\alpha^+} \in B\} = \bigcup_N \bigcap_{n \geq N}\{F_{\tau_n}^{\alpha^+} \in B\}
\in \vee_n \F_{\tau_n}\subset\F_{\tau}.$$

 D'une part, la commutation
de l'essentiel sup et l'espérance conditionnelle par rapport à la
filtration $\F_{\tau_n}$ permet d'établir l'expression suivante:
$$\mathbb{E}(F_{\tau_n}^{\alpha^+}) = \sup _{\mu_t = \alpha_t, \, t \leq \tau_n} \mathbb{E}\left[\,\mathbb{E}\left(k(\mu)|\,\F_{\tau_n}\right)\,\right] =
\sup _{\mu_t = \alpha_t, \, t \leq \tau_n}
\mathbb{E}\left(k(\mu)\right),$$
la suite $(F_{\tau_n}^{\alpha^+})$
étant une sur-martingale discrète positive  convergeant p.s. et dans
$L^1$ vers $F_{\infty}^{\alpha^+}$, la suite
$\mathbb{E}(F_{\tau_n}^{\alpha^+})$ décro\^{\i}t vers
$E[F_{\infty}^{\alpha^+}]$ et $\sup _{\mu_t = \alpha_t,
\, t \leq \tau_n} \mathbb{E}\left(k(\mu)\right)$
décroit vers $\mathbb{E}\left(k(\alpha)\right)$. 
 Par suite, nous avons $\mathbb{E}[F_{\infty}^{\alpha^+}]=
 \mathbb{E}\left(k(\alpha)\right)$. \\
 En utilisant  $F_{\tau_n}^{\alpha^+} - \mathbb{E}[k(\alpha) |\, \F_{\tau_n}]\geq 0,$
 nous obtenons $F_{\infty}^{\alpha^+}- k(\alpha)\geq 0$. Ainsi, nous avons $F_{\infty}^{\alpha^+}= k(\alpha)$.

D'autre part, en utilisant la remarque  \ref{remN},  $k_{\tau_n^+}(\alpha)$ converge vers  $\int_{\tau}^{+\infty}e^{-\beta s} f(\xi_s,Y_s) ds$ p.s. 
 lorsque
n tend vers l'infini. De l'expression (\ref{rel FW+}), il vient :
$$\lim_{n\rightarrow +\infty}  W_{\tau_n}^{\alpha^+} = \lim_{n\rightarrow +\infty} ( F_{\tau_n}^{\alpha^+} +k_{\tau_n^+}(\alpha)- k(\alpha))
= F_{\infty}^{\alpha^+}- k(\alpha)  + \int_{\tau}^{+\infty}e^{-\beta s} f(\xi_s,Y_s) ds.$$
Ce qui entra\^{\i}ne que
$W_{\tau_n}^{\alpha^+}$ converge vers  $ \int_{\tau}^{+\infty}e^{-\beta s} f(\xi_s,Y_s) ds$ p.s.
\endproof

 Le principe de la programmation
dynamique est un principe fondamental pour la théorie du contrôle
stochastique. Il a été initié dans les années cinquante par Bellman
et il s'énonce ainsi dans notre cas :

\begin{pro}
\label{progdyn}
Pour toute stratégie $\alpha \in \A$ et tout couple
$(\gamma,\theta)$ de $\G$-temps d'arrêt, $\gamma\in\TT^*,$  $\theta\in\TT$, nous
avons p.s.
\begin{eqnarray}
\label{inegcond}
W_\gamma^{\alpha} & \geq &
 \mathbb{E} \left[ \int_{\gamma  }^{\theta  }
e^{-\beta s}f(\xi_s,Y_s)\,ds -\displaystyle\sum_{\gamma\leq \tau_n < \theta} e^{-\beta \tau_n}
    c(\zeta_n,Y_{\tau_n^-}, \zeta_{n+1},Y_{\tau_n})  \,|\,\G_\gamma\right]  \nonumber\\
& + &
    \mathbb{E} \left(W_\theta^{\alpha}|\,\G_\gamma\right).
\end{eqnarray}
 Respectivement pour tout couple
$(\gamma,\theta)$, $0\leq \gamma\leq \theta$, nous avons p.s.
\begin{eqnarray}
\label{inegcond+}
 W_\gamma^{\alpha^+} & \geq &
 \mathbb{E} \left[ \int_{\gamma  }^{\theta }
e^{-\beta s}f(\xi_s,Y_s)\,ds -
\displaystyle\sum_{   \gamma< \tau_n \leq
\theta }  e^{-\beta \tau_n}
    c(\zeta_n,Y_{\tau_n^-}, \zeta_{n+1},Y_{\tau_n})\,|\, \F_\gamma\right]  \nonumber\\
& + &
    \mathbb{E}(W_\theta^{\alpha^+}|\,\F_\gamma).
    \end{eqnarray}
De plus, $\widehat{\alpha}$ est optimale si et seulement si
l'\'egalit\'e (\ref{inegcond+}) a lieu p.s. pour tout couple
$(\theta,\gamma)$.
\end{pro}
\noindent \textbf{Preuve  } \\
1. D'apr\`es la proposition \ref{sous-marting}, la variable aléatoire  $
F^\alpha_\theta$ est  un $(\G,\mathcal{T}^* ,\A)$
sur-martingale-syst\`eme, c'est \`a dire que
 pour $\gamma \leq\theta$
\begin{equation}\label{ingssmart}
p.s.~~F^\alpha_\gamma \geq \mathbb{E}
(F^\alpha_\theta\,|\G_{\gamma}).
\end{equation}
En \'ecrivant le gain maximal  conditionnel $F^\alpha_.$
sous la forme (\ref{rel FW}) dans l'in\'egalit\'e (\ref{ingssmart}),
nous obtenons
$$k (\alpha)- k_\gamma(\alpha)+ W^\alpha_\gamma \geq \mathbb{E} (k (\alpha)- k_\theta(\alpha)+ W^\alpha_\theta\,|\G_{\gamma}).$$
D'o\`{u}, 
$$ W^\alpha_\gamma \geq \mathbb{E} ( k_\gamma(\alpha)- k_\theta(\alpha)+ W^\alpha_\theta\,|\G_{\gamma}).$$
Nous obtenons ainsi l'in\'egalit\'e (\ref{inegcond}).\\

Nous appliquons le m\^{e}me raisonnement pour d\'emontrer
l'in\'egalit\'e (\ref{inegcond+}): D'apr\`es la
proposition \ref{sous-marting}, la v.a $ F^{\alpha^+}_\theta$ est un
$(\F,\mathcal{T} , \A)$ sur-martingale-syst\`eme,
c'est \`a dire que
 pour $\gamma \leq\theta$
\begin{equation}\label{ingssmart+}
p.s.~~F^{\alpha^+}_\gamma \geq \mathbb{E}
(F^{\alpha^+}_\theta\,|\F_{\gamma}).
\end{equation}
En \'ecrivant le gain maximal conditionnel $F^{\alpha^+}_.$ sous la forme (\ref{rel FW+})
dans l'in\'egalit\'e (\ref{ingssmart+}), nous obtenons
$$k (\alpha)- k_{\gamma^+}(\alpha)+ W^{\alpha^+}_\gamma \geq \mathbb{E} (k (\alpha)- k_{\theta^+}(\alpha)+ W^{\alpha^+}_\theta\,|\F_{\gamma}).$$
D'o\`{u},
$$ W^\alpha_{\gamma^+} \geq \mathbb{E} ( k_{\gamma^+}(\alpha)- k_{\theta^+}(\alpha)+ W^{\alpha^+}_\theta\,|\F_{\gamma}).$$
2. Supposons que la strat\'egie $\widehat{\alpha}$ est optimale. Par
suite, d'apr\`es le corollaire \ref{martingale},
$F^{\widehat{\alpha}^+}_.$ est  un $(\F,\mathcal{T} ,
\A)$-martingale-syst\`eme et on a l'\'egalit\'e
$$F^{\widehat{\alpha}^+}_\gamma = \mathbb{E}
(F^{\widehat{\alpha}^+}_\theta\,|\,\F_{\gamma}).$$
 R\'e\'ecrivons cette
\'egalit\'e en remplaçant $F^{\widehat{\alpha}^+}_.$ par
l'expression (\ref{rel FW+})  prise en $\theta$ et $\gamma$:
$$k (\widehat{\alpha})- k_{\gamma^+}(\widehat{\alpha})+ W^{\widehat{\alpha}^+}_\gamma = \mathbb{E} (k (\widehat{\alpha})-
k_{\theta^+}(\widehat{\alpha})+
W^{\widehat{\alpha}^+}_\theta\,|\F_{\gamma}).$$
 Ainsi  $k
(\widehat{\alpha})- k_{\gamma^+}(\widehat{\alpha})$ étant  $\F_\gamma$-
mesurable,  il passe  sous l'espérance conditionnelle et
$$W^{\widehat{\alpha}^+}_\gamma =\mathbb{E} \left(k_{\gamma^+}(\widehat{\alpha})-
k_{\theta^+}(\widehat{\alpha})+W_\theta^{\widehat{\alpha}^+}\,\big|\,\F_\gamma\right),$$
 l'\'egalit\'e (\ref{inegcond+}) est v\'erifi\'ee. \\
 \\
3. 
Réciproquement,  supposons que la stratégie $\alpha$ v\'erifie l'\'egalit\'e
(\ref{inegcond+}) pour tout $\gamma,\theta$, en particulier avec $\theta=\tau$:
$$ 
W_\gamma^{\alpha^+} =
 \mathbb{E} \left[ \int_{\gamma  }^{\tau}
e^{-\beta s}f(\xi_s,Y_s)\,ds -
\sum_{   \gamma< \tau_n <
\tau }  e^{-\beta \tau_n}
    c(\zeta_n,Y_{\tau_n^-}, \zeta_{n+1},Y_{\tau_n})+
   W_\tau^{\alpha^+}|\,\F_\gamma\right],$$ 
 où l'on remplace  $W_\gamma^{\alpha^+}$ et $W_\tau^{\alpha^+}$ par (\ref{rel FW+})  soit
 {\footnotesize{
 $$F^{\alpha^+}_\gamma = 
 \mathbb{E} \left[ k(\alpha) -k_{\gamma^+}(\alpha)+\int_{\gamma  }^{\tau}
e^{-\beta s}f(\xi_s,Y_s)\,ds -
\sum_{   \gamma< \tau_n <
\tau }  e^{-\beta \tau_n}
    c(\zeta_n,Y_{\tau_n^-}, \zeta_{n+1},Y_{\tau_n})+
   F^{\alpha^+}_\tau -
    k(\alpha) + k_{\tau^+}(\alpha)|\,\F_\gamma\right],$$ 
    }}
 et l'on simplifie utilisant les définitions et en particulier
 $k_{\tau^+} (\alpha) = \int_{\tau}^{+\infty}  e^{-\beta s} f(\xi_s^{\alpha},Y_s^{\alpha}) ~ ds $ :
 $$F^{\alpha^+}_\gamma = 
 \mathbb{E} \left[ 
   F^{\alpha^+}_\tau |\,\F_\gamma\right],$$ 
 $(F^{\alpha^+}_\gamma )$  est  un $(\F,\mathcal{T} ,\A)$-martingale système et le corollaire \ref{martingale} permet de conclure.
\endproof

Le second crit\`ere d'optimalit\'e permet d'examiner ce qui se passe
entre deux instants de changement de technologie et d'\'etat:
\begin{thm}
\label{thm optimalit}
 Pour toute strat\'egie admissible $\alpha$, 
tout couple $(i,x)\in {U}\times \mathbb{R}$ et pour $(\tau,\xi,Y)$   relatif \`a la strat\'egie
$\alpha$, 
nous avons  les in\'egalit\'es  p.s.  suivantes:
\begin{eqnarray}
\mbox{Pour tout }n\geq 0,\nonumber\\
 W^\alpha_{\tau_n}& \geq & -
e^{-\beta {\tau_n}}\int_{{U}\times
\mathbb{R}}c(\zeta_n,Y_{\tau_n^-},i,x)r(\zeta_n,Y_{\tau_n^-},i,dx)+
 \mathbb{E}(W^{\alpha^+}_{\tau_n}|\;\G_{\tau_n}).\nonumber\\\label{egalite2}\\
 \mbox{Pour tout }n\geq -1,\nonumber\\
W^{\alpha^+}_{\tau_n} & \geq & \mathbb{E}
\left(\int_{\tau_n}^{\tau_{n+1}}e^{-\beta
s}f(\xi_s,Y_s)\;ds|\;\F_{\tau_n}\right)+
\mathbb{E}(W^\alpha_{\tau_{n+1}}|\;\F_{\tau_n}). \label{egalite3}
\end{eqnarray}
On déduit de (\ref{egalite2}) prise en $n=0$ et  (\ref{egalite3}) prise en $n=-1$ :
 
\begin{equation}
W^{\alpha^+}_0  \geq  \mathbb{E}\left(\int_0^{\tau_0} e^{-\beta
s}f(\xi_s,Y_s)\;ds -  e^{-\beta
\tau_0}c(\xi_0,Y_{\tau_0^-},\zeta_1,Y_{\tau_0})|\;\F_0\right)
+
\mathbb{E}(W^{\alpha^+}_{\tau_0}|\;\F_0),
\label{egalite1}
\end{equation}
o\`{u}  $\F_0= \sigma (\xi_0, Y_{0})$.  
De plus, la strat\'egie
$\widehat{\alpha}$ est optimale si et seulement si l'\'egalit\'e a
lieu simultan\'ement dans  (\ref{egalite2}) et
(\ref{egalite3}).
\end{thm}
\noindent \textbf{Preuve}\\
1. L'essentiel supremum étant pris sur un ensemble plus restreint, nous
pouvons \'ecrire p.s.
\begin{eqnarray*}
W^\alpha_{\tau_n} &\geq& \ESS{\{\mu_t = \alpha_t,\, t \leq \tau_n^\alpha\}} \mathbb{E}\big(
\int_ {\tau_n^\mu}^{\tau_{N(\omega)}^\mu} e^{-\beta s}f(\xi_s^\mu,Y_s^\mu)\;ds -
\sum_{  n \leq k  < N }e^{-\beta\tau_k^\mu}
c(\zeta_k^\mu,Y_{(\tau_k^\mu)^-},
\zeta_{k+1}^\mu,Y_{\tau_k^\mu})\nonumber\\
 & + &    \int_ {\tau_{N(\omega)}^\mu}^{+\infty} e^{-\beta s}f(\zeta_{N(\omega)}^\mu,Y_s^\mu)\;ds   \big|\;\G_{\tau_n^\mu}\big).
\end{eqnarray*}
Puisque
$\mu_t = \alpha_t, \forall t \leq \tau_n $, nous avons
$\tau_n^\mu = \tau_n^\alpha$ et
$\G_{\tau_n^\mu} = \G_{\tau_n^\alpha}$.  Par suite, nous pouvons
\'ecrire l'\'egalit\'e:

$$\mathbb{E}(e^{-\beta \tau_n^\mu} c(\zeta_n^\mu,Y_{(\tau_n^\mu)^-}, \zeta_{n+1}^\mu,Y_{\tau_n}^\mu)|\;\G_{\tau_n^\mu}) = \int_{{U}\times
\mathbb{R}}e^{-\beta \tau_n^\alpha}c(\zeta_n^\alpha,
Y_{(\tau_n^\alpha)^-},i,x)r(\zeta_n^\alpha,
Y_{(\tau_n^\alpha)^-},i,dx).$$ Cette expression
 est $\G_{\tau_n}$-mesurable et ne d\'epend pas de
la strat\'egie $\mu$.  Elle peut  donc  \^{e}tre sortie
 de l'ess-sup pour obtenir p.s.

$$W^\alpha_{\tau_n} \geq -e^{-\beta \tau_n}\int_{{U}\times
\mathbb{R}}c(\zeta_n^\alpha,
Y_{(\tau_n^\alpha)^-},i,x)r(.,i,dx) +
 \ESS{\{\mu_t = \alpha_t,\, t \leq
\tau_n\}}\mathbb{E}
\left(k_{\tau_n^+}(\mu)|\;\G_{\tau_n^\mu}\right).$$

 Du fait que  $\F_{\tau_n^\mu} =
\F_{\tau_n^\alpha}$, nous avons:
\begin{eqnarray*}
W^\alpha_{\tau_n} &\geq&- e^{-\beta \tau_n}\int_{{U}\times
\mathbb{R}}c(\zeta_n^\alpha,
Y_{(\tau_n^\alpha)^-},i,x)r(.,i,dx)+
\mathbb{E}\left(\ESS{\{\mu_t = \alpha_t,\, t \leq \tau_n\}}
\mathbb{E}
(k_{\tau_n^+}(\mu)|\;\F_{\tau_n^\mu})|\;\G_{\tau_n^\alpha}\right).
\end{eqnarray*}
D'o\`{u} l'in\'egalit\'e (\ref{egalite2}).\\
\\
2. Par définition du gain maximal conditionnel apr\`es $\theta$,
\begin{eqnarray*}
W^{\alpha^+}_{\tau_n}  &=& \ESS{\{\mu_t = \alpha_t,\, t \leq
\tau_n\}} \mathbb{E}\big( \int_ {\tau_n^\mu}^{\tau_{N(\omega)}^\mu} e^{-\beta s}f(\xi_s^\mu,Y_s^\mu)\;ds - \sum_{  n < k < N} e^{-\beta
\tau^\mu_k}c(\zeta_k,Y^\mu_{\tau_k^-},
\zeta_{k+1},Y^\mu_{\tau_k^\mu}) \nonumber\\
 & + &    \int_ {\tau_{N(\omega)}^\mu}^{+\infty} e^{-\beta s}f(\zeta_{N(\omega)}^\mu,Y_s^\mu)\;ds  \big|\;\F_{\tau_n}\big)
\end{eqnarray*}
 Soit  dans
la famille des $\mu $  qui v\'erifient
  $\{ \mu_t = \alpha_t, \forall t \leq
\tau_{n}^\alpha\}$ celles qui coincident avec $\alpha$  jusqu'à
$(\tau_{n+1}^\alpha)^-,$ donc $\tau_{n+1}^\alpha = \tau_{n+1}^\mu$
et $ \int_ {\tau_n}^{\tau_{n+1}^\mu} e^{-\beta
s}f(\xi_s,Y^\mu_s)\;ds$ ne dépend pas de $\mu$. Ainsi,
 $$W^{\alpha^+}_{\tau_n}  \geq \mathbb{E}\left(\int_
{\tau_n^\alpha}^{\tau_{n+1}^\alpha} e^{-\beta s}f(\xi_s,Y_s)\;ds
|\;\F_{\tau_n}\right) + \mathbb{E}\left( k_{\tau_{n+1}}(\mu)
|\;\F_{\tau_n}\right)$$ Comme $\F_{\tau_n^\mu}\subset
\G_{\tau_{n+1}^\mu} $ la minoration peut se récrire :
$$
W^{\alpha^+}_{\tau_n} \geq  \mathbb{E}\left(\int_
{\tau_n^\alpha}^{\tau_{n+1}^\alpha} e^{-\beta s}f(\xi_s,Y_s)\;ds
|\;\F_{\tau_n}\right) +\mathbb{E}\left( \mathbb{E}\left(
k_{\tau_{n+1}} (\mu)|\;\G_{\tau_{n+1}}\right)|\,\F_{\tau_n}\right).
$$
Puisque $\mu$ coincide avec $\alpha$  jusqu'à
$(\tau_{n+1}^\alpha)^-,$
 par essentiel sup il vient:
$$W^{\alpha^+}_{\tau_n} \geq  \mathbb{E}\left(\int_
{\tau_n^\alpha}^{\tau_{n+1}^\alpha} e^{-\beta s}f(\xi_s,Y_s)\;ds +
\ESS{\{\mu_t = \alpha_t,\, t < \tau_{n+1}\}} \mathbb{E}\left(
k_{\tau_{n+1}}(\mu) |\;\G_{\tau_{n+1}}\right) |\;\F_{\tau_n}\right)
.$$

Enfin, l'essentiel sup commute avec la  $\F_{\tau_n}$-espérance
conditionnelle ce qui implique
l'in\'egalit\'e (\ref{egalite3}).\\
\\
3. Supposons que la strat\'egie $\widehat{\alpha}$ est optimale.
 D'une
part,  les in\'egalit\'es (\ref{egalite2}) et (\ref{egalite3})
impliquent
\begin{eqnarray*}
 W^{\widehat{\alpha}}_{\tau_n}& \geq & -e^{-\beta \tau_n}\int_{{U}\times
\mathbb{R}}c(\zeta_n,Y_{\tau_n^-},i,x)r(.,i,dx)+
 \mathbb{E}(W^{\widehat{\alpha}^+}_{\tau_n}|\;\G_{\tau_n})\nonumber\\
& \geq & -e^{-\beta \tau_n}\int_{{U}\times
\mathbb{R}}c(\zeta_n,Y_{\tau_n^-},i,x)r(.,i,dx)\nonumber
\\ &+&
 \mathbb{E}\left(\int_{\tau_n}^{\tau_{n+1}}e^{-\beta
 s}f(\xi_s,Y_s)\;ds +
W^{\widehat{\alpha}}_{\tau_{n+1}}|\;\G_{\tau_n}\right),\nonumber\\
\end{eqnarray*}
o\`{u} $(\tau,\xi,Y)$ est relatif \`a la strat\'egie
$\widehat{\alpha}$.\\
 D'autre part, l'\'egalit\'e (\ref{inegcond}) appliqu\'ee aux temps
d'arr\^{e}t $\tau_n$ et $\tau_{n+1}$ entra\^{\i}ne l'\'egalit\'e
\begin{eqnarray*}
W^{\widehat{\alpha}}_{\tau_n} &= & \mathbb{E}\left[ \int_
{\tau_n}^{\tau_{n+1}} e^{-\beta s}f(\xi_s,Y_s)\;ds -e^{-\beta
\tau_n}c(\zeta_n,Y_{\tau_n^-}, \zeta_{n+1},Y_{\tau_n})
|\;\G_{\tau_n}\right]+ \mathbb{E}\left[
W^{\widehat{\alpha}}_{\tau_{n+1}}|\;\G_{\tau_n}\right].
\end{eqnarray*}
D'o\`{u} (\ref{egalite2}) est une \'egalit\'e.\\
\\
4. Ci-dessous, le triplet  $(\tau,\xi,Y)$ est relatif \`a la
strat\'egie $\widehat{\alpha}$. En rempla\c{c}ant dans
l'in\'egalit\'e (\ref{egalite3})  le gain maximal conditionnel
$W_{\tau_{n+1}}^{\widehat{\alpha}}$ par l'\'egalit\'e
(\ref{egalite2}) appliqu\'ee au temps $\tau_{n+1}$, nous obtenons
\begin{eqnarray*}
W^{\widehat{\alpha}^+}_{\tau_n} &\geq& \mathbb{E}\left[ \int_
{\tau_n}^{\tau_{n+1}} e^{-\beta s}f(\xi_s,Y_s)\;ds +
W^{\widehat{\alpha}}_{\tau_{n+1}}
|\;\F_{\tau_n}\right]\nonumber\\&\geq& \mathbb{E}\left[ \int_
{\tau_n}^{\tau_{n+1}} e^{-\beta s}f(\xi_s,Y_s)\;ds
 +W^{\widehat{\alpha}^+}_{\tau_{n+1}} - e^{-\beta
\tau_{n+1}}c(\zeta_{n+1}, Y_{\tau_{n+1}^-},
\zeta_{n+2},Y_{\tau_{n+1}})|\;\F_{\tau_n}\right].\nonumber\\
\end{eqnarray*}
De plus, l'\'egalit\'e (\ref{inegcond+}) appliqu\'ee aux temps
d'arr\^{e}t $\tau_n$ et $\tau_{n+1}$ entra\^{\i}ne l'\'egalit\'e
\begin{eqnarray*}
W^{\widehat{\alpha}^+}_{\tau_n} &= & \mathbb{E}\left[ \int_
{\tau_n}^{\tau_{n+1}} e^{-\beta s}f(\xi_s,Y_s)\;ds - e^{-\beta
\tau_{n+1}}c(\zeta_{n+1}, Y_{\tau_{n+1}^-},
\zeta_{n+2},Y_{\tau_{n+1}}) |\;\F_{\tau_n}\right]\nonumber\\& +&
\mathbb{E}\left[
W^{\widehat{\alpha}^+}_{\tau_{n+1}}|\;\F_{\tau_n}\right].
\end{eqnarray*}
 Ce qui entra\^{\i}ne que (\ref{egalite3}) est une \'egalit\'e.\\
\\
 5. Supposons qu'il existe une strat\'egie admissible
$\alpha$ pour laquelle il y ait \'egalit\'e dans  (\ref{egalite1}),
(\ref{egalite2})
 et (\ref{egalite3}) pour tout $n \geq 0$ (respectivement
 pour tout $n \geq -1$). Nous avons, d'apr\`es
 l'\'egalit\'e (\ref{egalite1})
p.s.
$$W_0^{\alpha^+} = \mathbb{E}\left[
\int_ {0}^{\tau_{0}^\alpha} e^{-\beta
s}f(\xi_s^\alpha,Y_s^\alpha)\;ds -  e^{-\beta
\tau_0}c(\xi_0,Y_{\tau_0^-},\zeta_1,Y_{\tau_0}) \big|\;\F_0\right] +
\mathbb{E}\left[ W^{\alpha^+}_{\tau_{0}}\big|\;\F_0\right]. $$
 Nous posons l'hypoth\`ese
de r\'ecurrence jusqu'au rang $n$ suivante:
\begin{equation}\label{recurce}
   W^{\alpha^+}_{0} =  \mathbb{E}\left [\int_0^{\tau_{n}^\alpha} e^{-\beta s}f(\xi_s^\alpha,Y_s^\alpha)\;ds
-\sum_{k =0 }^{n} e^{-\beta \tau_k^\alpha}c(\zeta_k^\alpha,
Y_{(\tau_k^\alpha)^-}, \zeta_{k+1}^\alpha,Y_{\tau_k^\alpha})+
W^{\alpha^+}_{\tau_{n}}\big|\;\F_0\right].
\end{equation}
Montrons qu'elle est vraie au rang $n+1$. Les \'egalit\'es
(\ref{egalite3}) et (\ref{egalite2}) (prises en $(n+1)$) impliquent :
\begin{equation}\label{recurence n}
W_{\tau_n}^{\alpha^+} =  \mathbb{E}\left [
\int_{\tau_n^\alpha}^{\tau_{n+1}^\alpha} e^{-\beta
s}f(\xi_s,Y_s)\;ds -e^{-\beta
\tau_{n+1}^\alpha}c(\zeta_{n+1}^\alpha, Y_{(\tau_{n+1}^\alpha)^-},
\zeta_{n+2}^\alpha,Y_{\tau_{n+1}^\alpha}) +
W^{\alpha^+}_{\tau_{n+1}}\big|\;\F_{\tau_{n}}\right].
\end{equation}
En rempla\c{c}ant 
$W_{\tau_n}^{\alpha^+}$ par son expression (\ref{recurence n}),
l'\'egalit\'e (\ref{recurce}) devient
$$ W^{\alpha^+}_{0} =  \mathbb{E}\left [\int_0^{\tau_{n+1}^\alpha} e^{-\beta
s}f(\xi_s^\alpha,Y_s^\alpha)\;ds  -\sum_{k =0 }^{n+1} e^{-\beta
\tau_k^\alpha}c(\zeta_k^\alpha, Y_{(\tau_k^\alpha)^-},
\zeta_{k+1}^\alpha,Y_{\tau_k^\alpha})+
W^{\alpha^+}_{\tau_{n+1}}\big|\;\F_0\right].$$
 Ainsi l'hypoth\`ese
de r\'ecurrence est v\'erifi\'ee au rang $n+1$. Enfin, lorsque n
tend vers l'infini, $W_{\tau_n}^{\alpha ^+}$ converge vers $\int_{\tau}^{+\infty} f(\xi_s,Y_s) ds$
 p.s. (lemme \ref{admcout}) et par suite nous avons
$$W_0^{\alpha^+} = \mathbb{E} (k(\alpha) |\, \F_0).$$
En appliquant la propriété (\ref{rel FW+}) en $\theta=0$, nous obtenons  $F_0^{\alpha^+}=  W_0^{\alpha^+},$ soit
 $$F_0^{\alpha^+}= \mathbb{E} (k(\alpha)\big|\;\F_0)=E_{(i,x)}(k(\alpha)).$$
 Par conséquent, la strat\'egie $\alpha$ est optimale.
\endproof

\subsection{Propriétés markoviennes}
Le critère d'optimalité donné par le théorème \ref{thm optimalit}
 est insuffisant pour aider à la construction d'une stratégie
 optimale car les variables aléatoires $W^\alpha_{\tau_n}$ et $W^{\alpha^+}_{\tau_n}$ qui interviennent
 d\'ependent de la strat\'egie admissible $\alpha$. En 
utilisant le caract\`ere markovien et homog\`ene entre deux instants
d'impulsion, ainsi que la forme markovienne de chaque renaissance,
nous pouvons
 esp\'erer obtenir que les gains maximaux conditionnels ne d\'ependent
que de l'\'etat de syst\`eme \`a l'instant du conditionnement. \\
D'apr\`es le th\'eor\`eme de Doob, $\forall\, Y$ int\'egrable il
existe une fonction mesurable $g$ telle que:
$$\mathbb{E}\left[Y \big|\, \sigma(\xi_t, Y_t)
\right]:= g(\xi_t, Y_t).$$
 Ainsi, prenons en compte la notation:
$$
    \mathbb{E}_{\{i,x\}}(Y):= g(i,x).
$$
\begin{pro}
\label{mesrbrho}
 1. Introduisons les fonctions   $\rho$ et
$\rho^+$ sur   l'espace de probabilité  $({U}\times
\mathbb{R}, \mathcal{P}({U})\otimes\B_{
\mathbb{R}}, \mathbb{P})$:
$$\rho(i,x) = \ESS{\mu \in \A}~ \mathbb{E}_{\{i,x\}}(k(\mu)) \quad et\quad \rho^+(i,x) =
\ESS{\{\mu \in \A, \,\tau_0^\mu>0\}}
~\mathbb{E}_{\{i,x\}}(k(\mu)),$$
 o\`{u} l'ess-sup est pris au sens de la mesure de 
 Lebesgue. Ces deux fonctions sont $\mathcal{P}({U})\otimes\B_{
\mathbb{R}}$ mesurables.\\
\\ 
2.  De plus, 
\begin{equation}\label{comprho}
\rho(i,x) \geq \rho^+(i,x) \geq \mathbb{E}_{i,x}\left[ \int_0^{+\infty} e^{-\beta s} f(\xi_s,Y_s) ds \right] > 0.
\end{equation}
\end{pro}

\noindent \textbf{Preuve} \\
1. La mesurabilit\'e de ces fonctions,  sur  l'ensemble  de probabilité
$({U}\times \mathbb{R},
\mathcal{P}({U})\otimes\B_{ \mathbb{R}},
\mathbb{P})$,  est assur\'ee comme cons\'equence de  la proposition
6.1.1 de J. Neveu \cite{neveu1}:\\
 \textit{"Pour toute famille F de
fonctions r\'eelles mesurables  $ f: \Omega \longrightarrow
\mathbb{R}$ d\'efinies sur un espace de probabilit\'e
$(\Omega,\mathcal{ A}, \mathbb{P})$ il existe une, et \`a une
\'equivalence pr\`es une seule, fonction mesurable  $g: \Omega
\longrightarrow \mathbb{R}$ telle que
\begin{itemize}
    \item  $g \geq f$ p.s. pour tout $f\in F$,
    \item si $h$ est une fonction mesurable telle que $h \geq f$
    p.s. pour tout $f\in F$, alors $h \geq g$ p.s.
\end{itemize}
Cette fonction $g$, qui est la borne sup\'erieure de la famille $F$
au sens de l'in\'egalit\'e p.s., est not\'e $ \ESS~(F)$. En outre,
il existe au moins une suite $(f_n, n\in \mathbb{N})$ extraite  de F
telle
que $ \ESS~(F) = \sup f_n$ p.s. \\
Si la famille est filtrante croissante, la suite $(f_n, n\in
\mathbb{N})$ peut \^{e}tre choisie p.s. croissante et alors
$$ \ESS~(F) = \lim_n \uparrow f_n \quad p.s."$$ }
2. D'une part, puisque l'ess-sup est pris sur un ensemble plus restreint, il est trivial de prouver que 
$\rho(i,x) \geq \rho^+(i,x). $ \\
D'autre part, pour  la stratégie $\mu$ qui vérifie $\tau_0 =+ \infty$ et $\xi_t = i ~\forall t,$ nous obtenons:
$$\rho^+(i,x) \geq
\mathbb{E}_{\{i,x\}}\left[ \int_0^{+\infty} e^{-\beta s} f(\xi_s,Y_s) ds\right]. $$
D'où, la propriété  (\ref{comprho}) est vérifiée.

\endproof

\begin{pro}\label{recurence}
Pour toute stratégie $\mu$ et  toutes fonctions $f$ et $c$  telles que 
$$(\omega,t) \To e^{-\beta
(t-\tau_n)}f(\xi_t,Y_t) \1_{[\tau_k, \tau_{k+1}[}  \mbox{ et }\omega  \To e^{-\beta(\tau_k-\tau_n)}c(\zeta_k,Y_{\tau_k^-},\zeta_{k+1},Y_{\tau_k})  $$
sont intégrables $\forall n, \forall k\geq n$, il existe des fonctions
 mesurables  $(F_i)_{1 \leq i \leq 4}$ sur $\mathbb{N}\times{U}\times \mathbb{R}$ telles que:

\begin{eqnarray}
\mathbb{E}\left[A_{n,k+1}\int_{\tau_k}^{\tau_{k+1}}e^{-\beta
(t-\tau_n)}f(\xi_t,Y_t)dt\,\big|\,\F_{\tau_n}\right]&=&F_1(k-n,\zeta_{n+1},Y_{\tau_n}),\;\forall
k \geq n \geq -1\nonumber\\
\label{rec1}
\\
\mathbb{E}\left[ A_{n,k+1} \int_{\tau_k}^{\tau_{k+1}}e^{-\beta
(t-\tau_n)}f(\xi_t,Y_t)dt\,\big|\,\G_{\tau_n}\right]&=&F_2(k-n,\zeta_{n},Y_{\tau_n^-}),\;
\forall k \geq n\geq 0\nonumber\\ \label{rec2}
\\
\mathbb{E}\left[A_{n,k+1} e^{-\beta(\tau_k-\tau_n)}c(\zeta_k,Y_{\tau_k^-},\zeta_{k+1},Y_{\tau_k})\,\big|\,\F_{\tau_n}\right]&=&F_3(k-n,\zeta_{n+1},Y_{\tau_n}),\;
\forall k >n\geq -1 \nonumber\\
 \label{rec3}
\\
\mathbb{E}\left[A_{n,k+1}  e^{-\beta(\tau_k-\tau_n)}c(\zeta_k,Y_{\tau_k^-},\zeta_{k+1},Y_{\tau_k})\,\big|\,\G_{\tau_n}\right]&=&F_4(k-n,\zeta_{n},Y_{\tau_n^-}),\;
\forall k \geq n\geq 0, \nonumber\\ \label{rec4}
\end{eqnarray}
 où $A_{n,k+1} = \1_{\{\tau_n< \tau_{n+1}<\ldots<\tau_{k+1}\}}$.
\end{pro}

\noindent\textbf{Preuve:}  Pour montrer les assertions
pr\'ec\'edentes, il suffit de  proc\'eder par r\'ecurrence.
 Notons d'abord par $(\F_t^n , t \geq 0)$ la filtration
$(\F_{t+ \tau_n}, t \geq 0)$. Ensuite,
 mentionnons  que  les processus $Y_.^n$ et
$Y_{.+\tau_n}$ ont la m\^{e}me loi  $\F_{\tau_n}$-conditionnelle
sur l'intervalle $[0, \tau_{n+1}-\tau_n[$   où $Y_.^n$ est un
processus de Markov homog\`ene partant de $ Y_{\tau_n}$ et
v\'erifiant:
$$
\left\{%
\begin{array}{ll}
    dY_t^n = b(\zeta_{n+1},Y_t^n)dt + \sigma(\zeta_{n+1},Y_t^n)dW_t^n \\
    Y_0^n = Y_{\tau_n}, \\
\end{array}%
\right.
$$
o\`{u} $W^n$ est un mouvement brownien ind\'ependant de la filtration $\F_{\tau_n}$. 
De plus, le processus $(t, Y_.^n)$ est un processus de Markov non
homog\`ene pour la filtration $(\F_t^n , t \geq 0)$.\\
\\
1.\; Commen\c{c}ons par l'assertion (\ref{rec1}). Pour $l=k-n=0,\,
n\geq-1,$ nous avons:
$$\mathbb{E}\left[ A_{n,n+1}\int_{\tau_n}^{\tau_{n+1}}e^{-\beta
(t-\tau_n)}f(\xi_t,Y_t)dt\,\big|\,\F_{\tau_n}\right]=$$
\begin{eqnarray*}
&&
\mathbb{E}\left[ A_{n,n+1}\int_{\tau_n}^{\tau_{n+1} }e^{-\beta
 (t-\tau_n)}f\big(\zeta_{n+1},
Y^n_{t- \tau_n})\,dt\,\big|\,\F_{\tau_n}\right]\nonumber\\
&=& \mathbb{E}\left[ A_{n,n+1} \int_{0}^{\tau_{n+1} -  \tau_n }e^{-\beta
 t}f\big(\zeta_{n+1},
Y^n_t)\,dt\,\big|\,\F_{\tau_n}\right]
\end{eqnarray*}
En utilisant 2. de la définition \ref{strategadm}, nous obtenons:
$$\mathbb{E}\left[ A_{n,n+1}\int_{0}^{\tau_{n+1} -  \tau_n }e^{-\beta
 t}f\big(\zeta_{n+1},
Y^n_t)\,dt\,\big|\,\F_{\tau_n}\right]  $$
 $$= \mathbb{E}_{\{i,x\}}\left[  A_{-1,0} \int_{0}^{\tau_0}e^{-\beta
 t}f\big(\zeta_{n+1},Y_t^{-1})dt \right]_{\big|i=\zeta_{n+1},x=Y_{\tau_n}},$$
où $A_{-1,0} = \1_{\{\tau_0>0 \}}.$   On obtient donc l'existence d'une fonction $F_1(0,.,.)$ telle
que
$$
\mathbb{E}\left[ A_{n,n+1}\int_{\tau_n}^{\tau_{n+1}}e^{-\beta
(t-\tau_n)}f(\xi_t,Y_t)dt\,\big|\,\F_{\tau_n}\right]=F_1(0,\zeta_{n+1},Y_{\tau_n}),
$$
o\`{u}
$$
F_1(0,i,x) = \mathbb{E}_{\{i,x\}}\left[\1_{\{\tau_0>0 \}}  \int_{0}^{\tau_0}e^{-\beta
 t}f\big(i,
Y^x_t)\,dt\right].
$$
D'o\`{u}, l'assertion (\ref{rec1}) est
v\'erifi\'ee pour $l = 0$.\\
\\
$\G_{\tau_n}$ \'etant une sous-tribu de $\F_{\tau_n}$, on
conditionne ce dernier r\'esultat par $\G_{\tau_n}$, pour tout $n\geq0$. Ensuite,  en
se servant du th\'eor\`eme de Fubini et en utilisant la loi
conditionnelle $r $ de passage du couple 
 $(\zeta_n,Y_{\tau_n^-})$ au couple $(\zeta_{n+1},Y_{\tau_n})$ sachant $\G_{\tau_n}$ qui est ind\'ependante de $n$,
nous obtenons:
\begin{eqnarray*}
\mathbb{E}\left[F_1(0,\zeta_{n+1},
Y_{\tau_n})\,\big|\,\G_{\tau_n}\right]
 &=& \int_{{U}\times
\mathbb{R}}
F_1(0,j,y)\,r(\zeta_n,Y_{\tau_n^-},j,dy).
\end{eqnarray*}
Cette derni\`ere expression est effectivement une fonction mesurable
de $\zeta_n$ et $Y_{\tau_n^-}$. D'o\`{u} l'existence d'une fonction
$F_2(0,.,.)$  d\'efinie par
$$F_2(0,i,x)=\int_{{U}\times
\mathbb{R}} F_1(0,j,y)\,r(i,x,j,dy).$$
 Par cons\'equent, les assertions
(\ref{rec1}) et (\ref{rec2}) sont v\'erifi\'ees pour $l = k-n=0$.
Supposons qu'elles sont vraies jusqu' \`a $l = k-n$ et montrons qu'elles sont vraies au rang $l+1$.\\
\\
2.\; D'une part, du fait que $\F_{\tau_n} \subset \G_{\tau_{n+1}}$ et $A_{n,n+1}\in\G_{\tau_{n+1}},$
nous avons:
$$\mathbb{E}\left[A_{n,k+2} \int_{\tau_{k+1}}^{\tau_{k+2}}e^{-\beta (t-\tau_n)}f(\xi_t,Y_t)dt\,\big|\,\F_{\tau_n}\right]
=$$
$$
\mathbb{E}\left[A_{n,n+1}\mathbb{E}\left(A_{n+1,k+2}\int_{\tau_{k+1}}^{\tau_{k+2}}e^{-\beta
(t-\tau_n)}f(\xi_t,Y_t)dt\,\big|\,\G_{\tau_{n+1}}\right)\;\big|\,\F_{\tau_n}\right].$$
D'apr\`es  l'hypoth\`ese de r\'ecurrence (\ref{rec2}) appliqu\'ee
\`a $l = (k+1)-(n+1)$, nous obtenons:
$$\mathbb{E}\left[A_{n,n+1}\mathbb{E}\left(A_{n+1,k+2}\int_{\tau_{k+1}}^{\tau_{k+2}}e^{-\beta
(t-\tau_n)}f(\xi_t,Y_t)dt\,\big|\,\G_{\tau_{n+1}}\right)\;\big|\,\F_{\tau_n}\right] = $$
\begin{eqnarray*}
&& \mathbb{E}\left[ A_{n,n+1}  e^{-\beta(\tau_{n+1}-\tau_n)}
F_2(k-n,\zeta_{n+1},Y_{\tau_{n+1}^-})\,\big|\,\F_{\tau_n}\right]\nonumber\\
&=& \mathbb{E}_{\{i,x\}}\left[ A_{-1,0}  e^{-\beta\,\tau_0}F_2(k-n,\zeta_{n+1},Y^0_{\tau_{n+1} -
\tau_n})\right]_{\big| i=\zeta_{n+1},x=Y_{\tau_n}}.
\end{eqnarray*}
 La dernière égalité est donnée grâce à   la définition \ref{strategadm}. 
D'o\`{u} l'assertion (\ref{rec1}) est
v\'erifi\'ee pour $l+1$ avec
$$F_1(l+1,i,x)=\mathbb{E}_{\{i,x\}}\left[ \1_{\{\tau_0>0\}}  e^{-\beta\,\tau_0}\;F_2(l,i,Y^x_{\tau_0^-})\right].$$
D'autre part,  $\G_{\tau_n} \subset \F_{\tau_n}$ et par suite en
conditionnant l'assertion (\ref{rec1}) par $\G_{\tau_n}$ il vient
pour $k +1 =(l+1) +n$,
\begin{eqnarray*}
\mathbb{E}\left[F_1(l+1,\zeta_{n+1},Y_{\tau_n})\;\big|\,\G_{\tau_n}\right]&=&
\int_{{U}\times \mathbb{R}} F_1(l+1,j,y)
r (\zeta_n, Y_{\tau_n^-},j,dy).
\end{eqnarray*}
L'expression pr\'ec\'edente est bien une  fonction  de $l+1$ et
mesurable pour les v.a. $\zeta_n$ et $Y_{\tau_n^-}$ et l'assertion
(\ref{rec2}) est vraie $\forall k\geq n$ avec
$$F_2(l+1, i,x)=\int_{{U}\times \mathbb{R}} F_1(l+1,j,y)
r (i, x,j,dy).$$
\\
3.\;On calcule pour $l =k- n = 0,\, n\geq0$:
\begin{eqnarray*}
\mathbb{E}\left[A_{n,n+1} c(\zeta_n,Y_{\tau_n^-},\zeta_{n+1},Y_{\tau_n})\,\big|\,\G_{\tau_n}\right]
&=& \mathbb{E}\left[ c(\zeta_n,Y_{\tau_n^-},\zeta_{n+1},Y_{\tau_n}) \mathbb{E} \left(A_{n,n+1}  \big|\,\F_{\tau_n}\right)\,\big|\,\G_{\tau_n}\right]\nonumber\\
&=&
\mathbb{E}\left[c(\zeta_n,Y_{\tau_n^-},\zeta_{n+1},\Delta_n + Y_{\tau_n-}) \mathbb{E} \left(A_{-1,0} \,\big|\,\F_{\tau_n}\right)\,\big|\,\G_{\tau_n}\right]\nonumber\\
&=& \int_{{U}\times \mathbb{R}}
\mathbb{P}_{j,y}(\tau_0>0) c(\zeta_n,Y_{\tau_n-},j,y)\,r(\zeta_n,Y_{\tau_n-},j,dy).
\end{eqnarray*}
La derni\`ere égalité provient de la d\'efinition de la loi
conditionnelle $r$ de passage de
$(\zeta_n,Y_{\tau_n^-})$ \`a $(\zeta_{n+1},Y_{\tau_n})$ sachant
$\G_{\tau_n}$. Ainsi, on a l'existence de $F_4(0,.,.)$ d\'efinie par
$$F_4(0,i,x)=\int_{{U}\times \mathbb{R}}
\mathbb{P}_{j,y}(\tau_0>0) c(i,x,j,y)\,r(i,x,j,dy).$$
 Ensuite montrons que
l'assertion (\ref{rec3}) est vraie pour $k = n+1,\, n \geq-1$. Du
fait que $\F_{\tau_n} \subset \G_{\tau_{n+1}}$,  on multiplie l'assertion (\ref{rec4}), prise au rang $n+1$,  par $A_{n,n+1}e^{-\beta (\tau_{n+1}- \tau_n)}$, puis on la conditionne   par  $\F_{\tau_n}$
et l' on obtient 

$$\mathbb{E}\left[A_{n,n+1}   e^{-\beta (\tau_{n+1}- \tau_n)}F_4(0,\zeta_{n+1},Y_{\tau_{n+1}^-})\,\big|\,\F_{\tau_n}\right]=$$
$$
=\mathbb{E}_{\{i,x\}}\left[A_{-1,0} e^{-\beta \tau_0} F_4(0,\zeta_{n+1},Y^n_{\tau_{n+1}-\tau_n}\right]_{\big|i=\zeta_{n+1},x=Y_{\tau_n}},
$$

 la dernière égalité étant donnée par 2. de la définition \ref{strategadm}. Par suite, il existe $F_3(1,.,.)$
d\'efinie par
$$F_3(1, i,x)=\mathbb{E}_{\{i,x\}}\left[\1_{\{\tau_0>0\}} e^{-\beta \tau_0} F_4(0,i,Y^x_{\tau_0^-})\right] .$$
Ainsi l'assertion (\ref{rec3}) est v\'erifi\'ee pour $ k = n+1$.\\
\\
Par suite on pose les hypoth\`eses de r\'ecurrence (\ref{rec3}) et
(\ref{rec4}) jusqu'au rang $l =k-n$ et on montre qu'elles sont
vraies pour $l+1$.\\
\\
4. Pour tout $n \geq -1,$ $\F_{\tau_n}\in\G_{\tau_{n+1}}$  et $A_{n,n+1} \in \G_{\tau_{n+1}}$, nous obtenons:
$$
\mathbb{E}\left[A_{n,k+2} e^{-\beta(\tau_{k+1}-\tau_n)}c(\zeta_{k+1},Y_{\tau_{k+1}^-},\zeta_{k+2},Y_{\tau_{k+1}})\,\big|\,\F_{\tau_n}\right]$$
=$$ \mathbb{E}\left[A_{n,n+1}  \mathbb{E}\left
(A_{n+1,k+2} e^{-\beta(\tau_{k+1}-\tau_n)}c(\zeta_{k+1},Y_{\tau_{k+1}^-},\zeta_{k+2},Y_{\tau_{k+1}})\,\big|\G_{\tau_{n+1}}\right)\,\big|\,\F_{\tau_n}\right].
$$
D'apr\`es l'hypoth\`ese de r\'ecurrence (\ref{rec4}), appliqu\'ee
\`a $l =(k+1)-(n+1)$, l'\'egalit\'e pr\'ec\'edente devient:
$$\mathbb{E}\left[A_{n,k+2}e^{-\beta(\tau_{k+1}-\tau_n)}c(\zeta_{k+1},Y_{\tau_{k+1}^-},\zeta_{k+2},Y_{\tau_{k+1}})\,\big|\,\F_{\tau_n}\right]=$$
\begin{eqnarray*}
&&\mathbb{E}\left[A_{n,n+1}e^{-\beta(\tau_{n+1}-\tau_n)}F_4(k-n,\zeta_{n+1},Y_{\tau_{n+1}^-})\,\big|\,\F_{\tau_n}\right]\nonumber\\
&=& \mathbb{E}\left[A_{-1,0}e^{-\beta(\tau_{n+1}-\tau_n)}F_4(k-n,\zeta_{n+1},Y^0_{\tau_{n+1}-\tau_n})\right]_{\big|\zeta_{n+1} =i,Y_{\tau_n} =x}.
\end{eqnarray*}
 En utilisant la définition \ref{strategadm}, l'égalit\'e (\ref{rec3}) est vérifiée  $ \forall \, k > n
$ avec  l'existence de $F_3(l+1,.,.)$ d\'efinie par
$$F_3(l+1,i,x)=\mathbb{E}_{\{i,x\}}\left[\1_{\{\tau_0>0\}}e^{-\beta \tau_0}
F_4(l,i,Y^x_{\tau_{0}^-})\right].$$
Ensuite,  $\G_{\tau_n} \subset \F_{\tau_n}$. Par suite, pour
tout $n \geq 0,$ en conditionnant l'assertion  (\ref{rec3}) en $k
+1> n $ par la filtration $\G_{\tau_n}$, nous avons:
$$\mathbb{E}\left[A_{n,k+2}e^{-\beta(\tau_{k+1}-\tau_n)}c(\zeta_{k+1},Y_{\tau_{k+1}^-},\zeta_{k+2},Y_{\tau_{k+1}})\,\big|\, \G_{\tau_n}\right]$$
\begin{eqnarray*}
&=& \mathbb{E}\left[A_{n,n+1}   \mathbb{E}\left
(A_{n,k+2} e^{-\beta(\tau_{k+1}-\tau_n)}c(\zeta_{k+1},Y_{\tau_{k+1}^-},\zeta_{k+2},Y_{\tau_{k+1}})\,\big|\F_{\tau_{n}}\right)\,\big|\,\G_{\tau_n}\right]\nonumber\\
&=&
\mathbb{E}\left[A_{n,n+1} F_3(k+1-n,\zeta_{n+1},Y_{\tau_n})\,\big|\,\G_{\tau_n}\right].
\end{eqnarray*}
De la d\'efinition de la loi conditionnelle $r $ de
passage de $(\zeta_n,Y_{\tau_n^-})$ \`a $(\zeta_{n+1},Y_{\tau_n})$
sachant $\G_{\tau_n}$, il vient:
$$
\mathbb{E}\left[A_{n,n+1} F_3(k+1-n,\zeta_{n+1},Y_{\tau_n})\,\big|\,\G_{\tau_n}\right]=
\int_{{U}\times \mathbb{R}} \mathbb{P}_{j,y}(\tau_0>0)
F_3(k+1-n,j,y)r(\zeta_n,Y_{\tau_n-},j,dy),
$$
 On d\'eduit
l'existence de $F_4(l+1,.,.)$  d\'efinie par
$$F_4(l+1,i,x)=\int_{{U}\times
\mathbb{R}}\mathbb{P}_{j,y}(\tau_0>0)
F_3(l+1,j,y)r(i,x,j,dy).$$
\endproof

Au cours de cette preuve, nous avons démontré les propriétés suivantes:
\begin{cor}\label{relationrec}
Les fonctions  $(F_n(.,i,x))_{1 \leq n\leq 4}$  satisfont les relations de récurrence suivantes:

 \begin{eqnarray*}  
F_1(0,i,x)& =& \mathbb{E}_{\{i,x\}}\left[\1_{\{\tau_0>0\}} \int_{0}^{\tau_0}e^{-\beta
 t}f\big(i,
Y^x_t)\,dt\right]\nonumber\\ 
F_2(l, i,x)&=&\int_{{U}\times \mathbb{R}} F_1(l,j,y)
r (i, x,j,dy)\quad \hbox{si} \quad  l \geqslant 0\nonumber\\
F_1(l,i,x)&=&\mathbb{E}_{\{i,x\}}\left[ \1_{\{\tau_0>0\}}  e^{-\beta\,\tau_0}\, F_2(l-1,i,Y^x_{\tau_0^-})\right] \quad \hbox{si} \quad  l>0\nonumber\\
F_4(0,i,x)&=&\int_{{U}\times \mathbb{R}}\mathbb{P}_{j,y}(\tau_0>0)
c(i,x,j,y)\,r (i,x,j,dy) \nonumber\\
F_3(l,i,x)&=& \mathbb{E}_{\{i,x\}}\left[\1_{\{\tau_0>0\}} e^{-\beta \tau_0}
F_4(l-1,i,Y^x_{\tau_{0}^-})\right]\quad \hbox{si} \quad  l \geq 1 \nonumber\\ 
F_4(l,i,x)&=&\int_{{U}\times
\mathbb{R}}\mathbb{P}_{j,y}(\tau_0>0)
F_3(l,j,y)r (i,x,j,dy)\quad \hbox{si} \quad  l > 0.\nonumber\\
\end{eqnarray*}
\end{cor}

D'où le corollaire suivant:

\begin{cor}\label{coutrho}
Nous avons  les \'egalit\'es suivantes:
\begin{equation}\label{W rho}
   W_{\tau_n}^\alpha = e^{-\beta \tau_n^\alpha}\rho
(\zeta_n^\alpha,Y_{(\tau_n^\alpha)^-}),\; \forall\; n\geq 0 \qquad
\mbox{p.s.}
\end{equation}
\begin{equation}\label{W+rho}
W_{\tau_n}^{\alpha +} = e^{-\beta \tau_n^\alpha}\rho^+
(\zeta_{n+1}^\alpha,Y_{\tau_n^\alpha}), \; \forall\; n\geq -1\qquad
\mbox{p.s.}
\end{equation}
 En particulier, $W_0^{\alpha +} =\rho^+(i,x).$
\end{cor}

\noindent \textbf{Preuve: }
 On considère la trajectoire $(\xi,Y)$ à partir de $\tau^\alpha_n$ notée $(\tilde\xi,\tilde Y)$:
$$\rho(\zeta^\alpha_n,Y_{(\tau^\alpha_n)^-}) =   \ESS{\mu \in \A}\,
    \mathbb{E}\left[k(\mu)
    \,\big|\tilde\xi_0 = \zeta^\alpha_n, \tilde Y_0=Y_{(\tau^\alpha_n)^-}\right].$$
    Puisque l'on considère les trajectoires partant de $\tau_n^\alpha$,
    de fait l'essentiel sup est pris sur l'ensemble
 $\{\mu_t = \alpha _t,\; \forall \,t < \tau_n^\alpha\}$. En utilisant l'expression (\ref{kN}), pour de telles  stratégies,  $\tau_n^\mu =\tau_n^\alpha$  (not\'e $\tau_n$ ci-dessous),  il vient:
 
 \begin{eqnarray*} 
 k_{\tau_n}(\mu)(\omega)& =& \sum_{n \leq k < N}A_{n,k+1 } \int_{\tau_k}^{\tau_{k+1}} e^{-\beta s}f(\xi_s,Y_s)\,
ds -   \sum_{n \leq k<N}A_{n,k+1 } e^{-\beta \tau_k}  c(\zeta_k, Y_{\tau_k^-},
 \zeta_{k +1}, Y_{\tau_{k}}) \nonumber\\
 & +& \1_{\{\tau_n \leq \tau_N\}} \int_{\tau_{N(\omega)}}^{+\infty}  e^{-\beta s} f(\zeta_{N(\omega)},Y_s)ds.
\end{eqnarray*}
 Par ailleurs, nous pouvons utiliser  les assertions (\ref{rec2}) et (\ref{rec4}) de la
proposition \ref{recurence} pour obtenir:
  {\small
 \begin{eqnarray*}
 \mathbb{E}\left[k(\mu)(\omega)
    \,\big|\tilde\xi_0 = \zeta_n, \tilde Y_0=Y_{\tau_n^-}\right] 
 &=& \sum_{j \geq 0}  e^{-\beta \tau_n} F_2(j, \zeta_n,Y_{\tau_n^-}) -\sum_{j \geq 0}  e^{-\beta \tau_n} F_4(j,\zeta_n,Y_{\tau_n^-})\\
 &+&\1_{\{\tau_n \leq \tau_N\}}\mathbb{E} \left[\int_ {\tau_{N(\omega)}}^{+\infty}  e^{-\beta s} f(\zeta_{N(\omega)},Y_s)ds\big|\tilde\xi_0 = \zeta_n, \tilde Y_0=Y_{\tau_n^-}\right]. 
  \end{eqnarray*}
  }
Ensuite, en  prenant $j= k-n,$ la dernière expression peut être écrite sous la forme suivante:
 \begin{eqnarray*}
  &&\sum_{k \geq n} F_2(k-n, \zeta_n,Y_{\tau_n^-}) -\sum_{k \geq n} F_4(k-n, \zeta_n,Y_{\tau_n^-})\nonumber\\
  & + &  \1_{\{\tau_n \leq \tau_N\}} \mathbb{E} \left[\int_ {\tau_{N(\omega)}}^{+\infty}  e^{-\beta s} f(\zeta_{N(\omega)},Y_s)ds\big|\tilde\xi_0 = \zeta_n, \tilde Y_0=Y_{\tau_n^-}\right]= \mathbb{E}\left[e^{\beta \tau_n}k_{\tau_n}(\mu)
    \,\big|\G_{\tau_n}\right] ,
  \end{eqnarray*}
où l'intégrand est détaillé comme suit
{\small
\begin{eqnarray*} 
  e^{\beta \tau_n}k_{\tau_n}(\mu)& =& \sum_{n \leq k <  N}A_{n,k+1 } \int_{\tau_k}^{\tau_{k+1}} e^{-\beta
(s-\tau_n)}f(\xi_s,Y_s)\,
ds -  \sum_{n \leq k<N}A_{n,k+1 }  e^{-\beta
(\tau_k-\tau_n)} c(\zeta_k, Y_{\tau_k^-},
 \zeta_{k +1}, Y_{\tau_{k}}) \nonumber\\
 & +&\1_{\{\tau_n \leq \tau_N\}}  \int_{\tau_{N(\omega)}}^{+\infty}  e^{-\beta(s-\tau_n)} f(\zeta_{N(\omega)},Y_s)ds. 
\end{eqnarray*}
}
L'essentiel sup étant pris sur l'ensemble $\{\mu_t = \alpha _t,\; \forall \,t < \tau_n^\alpha\}$, on obtient 
 $$\rho(\zeta_n,Y_{\tau_n^-}) = \displaystyle \ESS{\{\mu_t = \alpha _t,\; \forall \,t <
    \tau_n\}}
    \mathbb{E}\left[e^{\beta \tau_n}k_{\tau_n}(\mu)
    \,\big|\G_{\tau_n}\right],$$
    où l'on reconnait $ e^{\beta \tau_n}\; W_{\tau_n}^{\alpha}$
    (définition \ref{coutcond}). 
Ainsi, l'\'egalit\'e (\ref{W rho}) est v\'erifi\'ee.\\
 
   De mani\`ere analogue,
 nous v\'erifions l'\'egalit\'e
(\ref{W+rho}), pour tout $n \geq -1$: 
On considère  la trajectoire $(\xi,Y)$ qui part à droite de $\tau_n^\alpha$ notée $(\tilde\xi,\tilde Y)$:
$$\rho^+(\zeta_{n+1},Y_{\tau_n}) = \ESS{\mu \in \A,\,\tau_0^\mu>0}
    \mathbb{E}\left[k(\mu)
    \,\big|\tilde\xi_0 = \zeta_{n+1}^\alpha, \tilde
    Y_0=Y_{\tau_n^\alpha}\right].$$
Puisque l'on considère les trajectoires partant à droite de
$\tau_n^\alpha$,
    de fait l'essentiel sup est pris sur l'ensemble
 $\{\mu_t = \alpha _t,\; \forall \,t \leq \tau_n^\alpha\}$. Pour de telles
 stratégies,  $\tau_n^\mu =\tau_n^\alpha$ (not\'e $\tau_n$) et  il vient:
 {\small 
 \begin{eqnarray*}
  \mathbb{E}\left[k(\mu)
    \,\big|\tilde\xi_0 = \zeta_{n+1}, \tilde Y_0=Y_{\tau_n}\right] 
    &=& \sum_{j \geq 0} e^{-\beta \tau_n}\; F_1(j,\zeta_{n+1},Y_{\tau_n} ) -\sum_{j > 0} e^{-\beta \tau_n}\;F_3(j,\zeta_{n+1},Y_{\tau_n})\nonumber\\
   &+&\1_{\{\tau_n \leq \tau_N\}}\mathbb{E} \left[\int_ {\tau_{N(\omega)}}^{+\infty}  e^{-\beta s} f(\zeta_{N(\omega)},Y_s)ds\big|\tilde\xi_0 = \zeta_{n+1}, \tilde Y_0=Y_{\tau_n}\right]. 
  \end{eqnarray*}}
La dernière égalité est obtenue en utilisant les assertions (\ref{rec2}) et (\ref{rec4})  de la proposition \ref{recurence}. 
 Ensuite, en  prenant $j= k-n,$ la dernière expression peut être écrite sous la forme suivante:
  {\small
  \begin{eqnarray*}
&&\sum_{k \geq n} F_1(k-n, \zeta_{n+1},Y_{\tau_n}) -\sum_{k > n} F_3(k-n,\zeta_{n+1},Y_{\tau_n}) \nonumber\\
&+& \1_{\{\tau_n \leq \tau_N\}}\mathbb{E} \left[\int_ {\tau_{N(\omega)}}^{+\infty}  e^{-\beta s} f(\zeta_{N(\omega)},Y_s)ds\big|\tilde\xi_0 = \zeta_{n+1}, \tilde Y_0=Y_{\tau_n}\right] =\mathbb{E}\left[e^{\beta \tau_n}k_{\tau_n^+}(\mu)\big| \F_{\tau_n}\right],
  \end{eqnarray*}}
 où l'intégrand est détaillé comme suit
{\small
\begin{eqnarray*} 
  e^{\beta \tau_n}k_{\tau_n^+}(\mu)& =& \sum_{n \leq k < N}A_{n,k+1 }\int_{\tau_k}^{\tau_{k+1}} e^{-\beta
(s-\tau_n)}f(\xi_s,Y_s)\,
ds -  \sum_{n <k<N}A_{n,k+1 }  e^{-\beta
(\tau_k-\tau_n)} c(\zeta_k, Y_{\tau_k^-},
 \zeta_{k +1}, Y_{\tau_{k}}) \nonumber\\
 & +& \1_{\{\tau_n \leq \tau_N\}} \int_{\tau_{N(\omega)}}^{+\infty}  e^{-\beta(s-\tau_n)} f(\zeta_{N(\omega)},Y_s)ds. 
\end{eqnarray*}}
 L'essentiel sup étant pris sur l'ensemble $\{\mu_t = \alpha _t,\; \forall \,t \leq \tau_n^\alpha\}$, on obtient 
$$\rho^+(\zeta_{n+1},Y_{\tau_n}) = \displaystyle \ESS{\{\mu_t = \alpha _t,\; \forall \,t \leq
    \tau_n\}}
   \mathbb{E}\left[e^{\beta \tau_n}k_{\tau_n^+}(\mu)
    \,\big|\F_{\tau_n}\right],$$
       où l'on reconnait $e^{\beta \tau_n}\; W_{\tau_n}^{\alpha^+}$
    (définition \ref{coutcond}).
\endproof
D'où la remarque suivante:
\begin{rem}
La fonction $\rho$ présente le gain moyen de la firme avant le changement de technologie, tandis que la fonction $\rho^+$ est le gain moyen de la firme juste après le saut.
\end{rem}

 A la suite de J.P. Lepeltier et B. Marchal \cite{lepl}, 
 nous introduisons la définition suivante:

\begin{definition}\label{applimrho}
Pour tout $(i,x) \in {U}\times\mathbb{R}$,
 nous définissons l'application suivante:
$$m \rho^+(i,x) :(i,x)\longmapsto \ESS{\nu \in M_{(i,x)} }\int_{{U}\times\mathbb{R}}\nu
(i,x;j,dy)\, \left(-c(i,x,j,y) + \rho^+(j,y)\right).$$
Posons, de plus, pour tout ensemble
$M^*_{(i,x)}=M_{(i,x)}-\delta_{(i,x)}$
$$ m^\ast\rho^+(i,x) : (i,x)\longmapsto \ESS{\nu \in M^*_{(i,x)}}\int_{{U}\times\mathbb{R}}\nu (i,x;j,dy)
  (-c(i,x,j,y) + \rho^+(j,y)) .$$
  \end{definition}

\begin{rems}
\label{rhomesurable}
1. Les applications $m\rho^+$ et $m^*
\rho^+$ sont mesurables comme essentiel supremum des fonctions mesurables 
$$ (i,x)\longmapsto  \int_{{U}\times\mathbb{R}}\nu
(i,x;j,dy)\, \left(-c(i,x,j,y) + \rho^+(j,y)\right).$$
  2. Les applications $\rho^+$, $\rho$, $m\rho^+$, $m^\ast\rho^+$, $m\rho$, $m^\ast\rho$ sont li\'ees par les
  relations suivantes:
$$\left\{%
\begin{array}{ll}
    m\rho^+(i,x)=\rho^+(i,x)\vee m^\ast\rho^+(i,x) \\
   m\rho(i,x)=\rho(i,x)\vee   m^\ast\rho(i,x). 
\end{array}%
\right.    $$
\end{rems}

\noindent \textbf{Hypothèse 2:} 
  \\
  (i)  la $\mathbb{P}_{(i,x)}$-loi du couple  $(\tau_0, (\xi_.,Y_.) \1_{[0,\tau_0)})$  est  faiblement  continue en $x$.
  \\
  (ii) les fonctions $f$ et $c$ sont  continues bornées.
  \\
(iii) Pour toute fonction $f$ borélienne bornée,  $\forall r\in M, ~   (i,x)\longmapsto r(i,x;f)$ est continue.

  \begin{pro}
\label{prnoyaubor}   Sous  les hypothèses   1 et  2 (iii), il existe un noyau
bor\'elien $r^*$ tel que $r^*(i,x,.,.)\in M_{(i,x)}$ vérifiant,
\begin{equation}\label{noyaubor}
   m\rho^+(i,x) = \int_{{U}\times \mathbb{R}}r^\ast(i,x;j,dy)(-c(i,x,j,y)+\rho^+(j,y)).
\end{equation}
\end{pro}

\noindent \textbf{Preuve }\\
1. D'après la définition de la fonction $m\rho^+$ et les propriétés de
l'essentiel supremum il existe, $\forall (i,x)$  $\forall n$, un noyau
bor\'elien $r^n_{(i,x)} \in M$
 tel que,
$$
    -\frac{1}{n}+  m \rho^+ (i,x) \leq \int_{{U}\times \mathbb{ R}} r^n_{(i,x)} (i,x;j,dy)\left(-c(i,x;j,y) +\rho^+(j,y)\right)
    \leq  m \rho^+ (i,x).
$$
 L'ensemble $M_{(i,x)}$  étant    compact  fermé pour la topologie faible (hypothèse 1), il existe une
suite extraite $(r^{n_j}_{(i,x)})_{j\geq 0}$   et une mesure 
$  r^\ast_{(i,x)}\in  M_{(i,x)}$ telles que
 $$  r^\ast_{(i,x)}(i,x,..)\mbox{ est la limite faible de la suite }
 (r^{n_j}_{(i,x)}(i,x,..))_j.$$
  $r^\ast_{(i,x)}(i,x,..) $ est une probabilité obtenue comme limite faible 
d'une suite extraite $(r^{n_j}_{(i,x)})_{j\geq 0}$ de
probabilités de l'espace compact $M_{(i,x)}$.\\
\\
2. Ensuite,  appliquons le  théorème $3.38$ de C. Castaing \cite[p. 85]{castaing}:\\

(i) On a la multi-application
\begin{eqnarray*}
\Sigma : {U}\times\mathbb{ R} &\To& \mathcal{P}({U}\times\mathbb{ R})\nonumber\\
(i,x) &\To& \{r(i,x,.,.), ~ r\in M\}
\end{eqnarray*}
 où $\mathcal{P}({U}\times\mathbb{ R})$ est l'ensemble des probabilités sur ${U}\times\mathbb{ R}$ muni de la tribu $\mathcal{T}$ des boréliens de la topologie faible. C'est un ensemble complet, métrisable et séparable (et donc un espace Polonais).
Pour prouver que $graph ~ \Sigma= \{((i,x), r(i,x,.,.))\}$ est mesurable,  nous utilisons le lemme suivant que nous citons par souci de complétude: \\

 {\it {\bf Lemme 2  } \cite[p. 135]{castaing1}:  Soient $(\Omega,\A, \mu)$ un espace mesuré avec $\mu$  positive finie et $\A^\ast$ le prolongement de Lebesgue de $\A.$ Soient $E$ un espace Polonais et  $\Gamma$ une multi-application de $\Omega$ à valeurs dans les fermés non vides de $E$. Alors les conditions suivantes sont équivalentes:\\
  - Pour tout $x$ fixé dans $E$, la fonction $\omega \To d(x,\Gamma(\omega))$ est $\A^\ast$ -mesurable.\\
  - $\Gamma $ est de graphe mesurable, c'est  à dire son graphe appartient à $\A^\ast \otimes\mathcal{B}$, où $\mathcal{B}$ est la tribu borélienne de $E$.} 
  \\
 Dans notre cas,    $(\Omega, \F,(\F_{t}),\mathbb{P})$ est  un espace de probabilité complet filtré et $\mathcal{P}({U}\times\mathbb{ R})$   est  un espace  Polonais.  De plus, l'hypothèse 1 donne que $\Sigma$ est à valeurs dans les fermés. Par ailleurs, d'après le théorème 2.19 \cite[p. 25]{bain}, il existe un sous-ensemble dénombrable de fonctions continues  $(\varphi_j)$ sur $\mathcal{P}({U}\times\mathbb{ R})$ avec $\| \varphi_j\| = 1\; \forall j,$ tel que 
  \begin{eqnarray*}
d : \mathcal{P}({U}\times\mathbb{ R}) \times \mathcal{P}({U}\times\mathbb{ R}) &\To& [0,+\infty[\nonumber\\
(P,Q) &\To& \sum_{j =1}^{+\infty}\frac{|P \varphi_j - Q \varphi_j |}{2^j}
\end{eqnarray*}
  définit une métrique  sur $\mathcal{P}({U}\times\mathbb{ R})$ qui engendre la topologie faible. 
 Pour tout $P \in \mathcal{P}({U}\times\mathbb{ R})$, définissons l'application:
 $$(i,x)\longmapsto d(P, \Sigma_{i,x}) =  \mbox{ess}\inf\{d(P,Q), Q \in \Sigma_{i,x} \} = \INF{r\in M} \;d(P,r(i,x,.)),$$
 où 
  $$d(P,r(i,x,.)) =  \sum_{j =1}^{+\infty}\frac{|P \varphi_j - r(i,x,\varphi_j) |}{2^j}.$$
 Sous l'hypothèse 2 (iii), $\forall j ~  x \longmapsto r(i,x,\varphi_j) $ est mesurable. Par conséquent, $\forall r \in M$ $(i,x)\longmapsto  d(P,r(i,x,.))$ est  mesurable comme limite croissante de fonctions mesurables.    Enfin,    $x \longmapsto \INF{r\in M} ~d(P,r(i,x,.))$ est  mesurable comme un essentiel  infimum de fonctions mesurable.    Par conséquent, d'après le lemme 2 \cite[p. 135]{castaing1}, le graphe de $\Sigma$ est mesurable. \\
\\
(ii) Nous utilisons le résultat suivant que nous citons par souci de complétude:\\
\\
 {\it {\bf Proposition 12.4} \cite[p. 74]{kechris} : Soient $(X,\A )$ un espace mesurable, Y espace séparable métrique et $f: X\to Y$ une fonction mesurable. Alors,  $graph(f) \subset X\times Y $ est aussi mesurable.}  \\
Dans notre cas, $({U}\times\mathbb{ R}, \mathcal{P}({U})\otimes \mathcal{B}_{\mathbb{ R}})$ 
est un espace mesurable, $\R$ un espace séparable métrique muni de la distance usuelle et   l'application $m\rho^+: {U}\times\mathbb{ R} \To\R$  définie dans la définition \ref{applimrho}    est mesurable (proposition \ref{rhomesurable}). Cette proposition 12.4  \cite[p. 74]{kechris}
montre alors que $graph(m\rho^+)$ est mesurable.\\
\\
 (iii)  Puisque $(-c + \rho^+)$ est une fonctions bornée mesurable sur  $(U\times \R)^2$,   la fonction $g: {U}\times\mathbb{ R}\times \mathcal{P}({U}\times\mathbb{ R}) \To \mathbb{ R}$ 
 définie par: $$g(i,x,\mathbb{ P}) =  \int_{{U}\times\mathbb{ R}} \left(-c(i,x;j,y) +\rho^+(j,y)\right)\mathbb{ P}(j,dy)$$
 est mesurable.  \\
\\
(iv) D'après le point 1 ci-dessus, $\forall (i,x)$, il existe $r^*_{(i,x)}(i,x;.) \in M_{(i,x)}$ tel que 
$$m \rho^+ (i,x) = \int_{{U}\times\mathbb{ R}} r^*_{(i,x)} (i,x;j,dy)\left(-c(i,x;j,y) +\rho^+(j,y)\right).$$
Ainsi,  $m\rho^+(i,x) \in g(i,x,\Sigma(i,x))$ et   $ \forall (i,x)$ $g(i,x,\Sigma(i,x))\cap m\rho^+(i,x) \neq \emptyset$. Par suite, d'après le théorème  $3.38$ \cite[p. 85]{castaing},  il existe une sélection mesurable de $\Sigma$ notée  $r^\ast$ telle que $g(i,x,r^\ast(i,x;.,.)) = m\rho^+(i,x)$  
et le noyau $r^*$ vérifie donc l'égalité (\ref{noyaubor}). 
\endproof

Pour la suite, on va utiliser une topologie plus forte que celle introduite en  (\ref{topfaible}):  \\

\noindent\textbf{Hypothèse 3:}  L'ensemble 
 $A=\{r(.,-c+\rho^+),~r\in M\} $
est fermé et compact pour la   topologie 
\begin{equation}\label{topHyp3}
(r_n(.,-c+\rho^+))_n ~\mbox{converge uniformément  vers   }   r(.,-c+\rho^+)   ~\mbox{ sur tout compact de } \R.
\end{equation}

\begin{pro}\label{opcont}
 Sous l'hypothèse 3,
 \\
 1. l'application 
 $$x\longmapsto \int_{U \times\mathbb{R}}r^*
(i,x;j,dy)\, \left(-c(i,x,j,y) + \rho^+ (j,y)\right)$$
 coincide avec l'application $m\rho^+,$
\\

2.    $x\longmapsto m\rho ^ +(i,x)$ est continue sur $\R.$

\end{pro}
\noindent\textbf{Preuve:} 
\\
1.  Pour tout  $ x\in D,$ où  $D$ est un ensemble dénombrable dense,  d'après  les propriétés de
l'essentiel supremum il existe
  une suite  $(r_{n,x}) \in M$ telle que  $ r_{n,x}(i,x,-c(i,x,j,y) + \rho^+ (j,y))$ converge vers  $r^\ast(i,x,-c(i,x,j,y) + \rho^+ (j,y)).$\\
 L'ensemble  $A$ est fermé et compact pour la   topologie définie dans l'hypothèse  3, d'où il existe une sous-suite $(r_{n_k})$ de limite $\widehat{r} \in M$ telle que $r_{n_k}(.,-c + \rho^+ )$ converge uniformément vers $\widehat{r} (.,-c  + \rho^+  )$ sur tout compact de  $\R.$  

Nous obtenons grâce à l'unicité de la limite, $ \forall x\in D, ~ \forall r \in M$:
  {\small $$r^\ast (i,x,-c(i,x,j,y) + \rho^+ (j,y)) = \widehat{r} (i,x,-c(i,x,j,y) + \rho^+ (j,y)) \geq r  (i,x,-c(i,x,j,y) + \rho^+ (j,y)). $$}
  Soit $K$ un compact de $\R$.  Pour tout $x \in K$ et $r \in M$,  introduisons $\eta_{\widehat{r} }$ et   $\eta_{r },$ les modules d'uniforme continuité de $\widehat{r} $ et $r$ sur $K$. Ensuite, il existe  $x_n \in D \cap B(x, \eta_{\widehat{r} })\cap B(x, \eta_{r})$ tel que, $ \forall \eps >0:$
  
 \begin{eqnarray*}
 \widehat{r} (i,x,-c(i,x,j,y) + \rho^+ (j,y))& \geq&  \widehat{r} (i,x_n,-c(i,x,j,y) + \rho^+ (j,y)) - \eps \nonumber\\
 &\geq & r  (i,x_n,-c(i,x,j,y) + \rho^+ (j,y)) -\eps \nonumber\\
 &\geq & r  (i,x,-c(i,x,j,y) + \rho^+ (j,y)) -2 \eps.
 \end{eqnarray*} 
  Par suite, $\forall \eps >0:$
  $$ \widehat{r} (i,x,-c(i,x,j,y) + \rho^+ (j,y)) \geq r  (i,x,-c(i,x,j,y) + \rho^+ (j,y)) -2 \eps.$$
   Ce qui implique que, $\forall x \in K$ et $\forall  K \subset \R:$  
  $$\widehat{r} (i,x,-c(i,x,j,y) + \rho^+ (j,y)) \geq r  (i,x,-c(i,x,j,y) + \rho^+ (j,y)).$$
  Ainsi,  $\forall x\in \R,$
   $$\widehat{r} (i,x,-c(i,x,j,y) + \rho^+ (j,y)) = \sup_{r \in M} r (i,x,-c(i,x,j,y) + \rho^+ (j,y))$$ 
  où l'on reconnait $m\rho^+(i,x)$ ce qui montre 1.\\
   
2. Les hypothèses 2 (ii) (iii) montrent la continuité de
   $$ x\mapsto \widehat{r} (i,x,-c(i,x,j,y) + \rho^+ (j,y)),$$
   d'où celle de
   $$ x\mapsto  r^\ast (i,x,-c(i,x,j,y) + \rho^+ (j,y)).$$

 \endproof

\begin{pro}\label{mrhorho}
Pour toute strat\'egie admissible $\alpha$ et tout $n\geq 0$, on a
\begin{equation}\label{promrhw}
    W_{\tau_n}^\alpha = e^{-\beta \tau_n^\alpha} \, m\rho^+ (\zeta_n^\alpha,Y_{(\tau_n^\alpha)^-}) \quad \mbox{p.s.}
\end{equation}
De plus, pour toute stratégie  $\alpha \in \A$ et tout $n\geq 0,$ on a
 $$    m\rho^+(\zeta_n^\alpha, Y_{(\tau_n^\alpha)^-}) = \rho(\zeta_n^\alpha, Y_{(\tau_n^\alpha)^-}) \quad \mbox{p.s}.$$
\end{pro}

\noindent\textbf{Preuve}\\
1. Il  s'agit d'\'etablir $\forall \alpha \in \A$
l'\'egalit\'e
\begin{equation}\label{mrhW}
   e^{-\beta \tau_n^\alpha} \, m\rho^+ (\zeta_n^\alpha,Y_{(\tau_n^\alpha)^-}) = \displaystyle \ESS{\{\mu_t = \alpha _t,\; \forall \,t <
    \tau_n^\alpha\}}
    \mathbb{E}(k_{ \tau_n^\alpha}(\mu)|\,\G_{\tau_n^\alpha})\quad \mbox{ p.s.}
\end{equation}
 Soit une strat\'egie $\mu$ v\'erifiant $\{\mu_t = \alpha _t,\;
\forall \,t < \tau_n\}$, alors   $\tau_n^\mu = \tau_n^\alpha$ et
 $r(\zeta_n^\alpha,Y_{(\tau_n^\alpha)^-},.,.)\in
M_{\{\zeta^\alpha_n,Y_{(\tau^\alpha_n)^-}\}}$, donc l'op\'erateur
$m\rho^+ (\zeta_n^\alpha,Y_{(\tau_n^\alpha)^-})$ v\'erifie:
$$
e^{-\beta \tau_n^\alpha} \, m\rho^+
(\zeta_n^\alpha,Y_{(\tau_n^\alpha)^-}) \geq  e^{-\beta
\tau_n^\alpha} \int_{{U}\times
\mathbb{R}}r(\zeta_n^\alpha,Y_{(\tau_n^\alpha)^-},i,dx)(-c(\zeta_n^\alpha,Y_{(\tau_n^\alpha)^-},i,x)+\rho^+(i,x)).$$
 De plus, pour toutes ces stratégies $\mu$ on a:
\begin{eqnarray*}
&& e^{-\beta \tau_n^\alpha}\int_{{U}\times
\mathbb{R}}r(\zeta_n^\alpha,Y_{(\tau_n^\alpha)^-},i,dx)(-c(\zeta_n^\alpha,Y_{(\tau_n^\alpha)^-},i,x)+\rho^+(i,x))\nonumber\\
&=& e^{-\beta \tau_n^\alpha}
\mathbb{E}\big(-c(\zeta_n^\mu,Y_{(\tau_n^\mu)^-},\zeta_{n+1}^\mu,Y_{\tau_n^\mu})
+ \rho^+(\zeta_{n+1}^\mu,Y_{\tau_n^\mu})|\,\G_{\tau_n^\alpha}\big).
\end{eqnarray*}
Gr\^{a}ce \`a l'\'egalit\'e (\ref{W+rho}), on peut remplacer
$e^{-\beta \tau_n^\alpha}\rho^+(\zeta_{n+1}^\mu,Y_{\tau_n^\mu})$ :
$$e^{-\beta \tau_n^\alpha} \, m\rho^+ (\zeta_n^\alpha,Y_{(\tau_n^\alpha)^-}) \geq \mathbb{E}
\left(-e^{-\beta
\tau_n}c(\zeta_n^\mu,Y_{(\tau_n^\mu)^-},\zeta_{n+1}^\mu,Y_{\tau_n^\mu})+
W_{\tau_n}^{\mu+}|\,\G_{\tau_n^\alpha}\right).$$
 Soit encore, puisque  $W_{\tau_n}^{\mu+}\geq \mathbb{E}(k_{\tau_n^+}(\mu)|\,
\F_{\tau_n^\mu})$, pour  cette strat\'egie $\mu$ qui vérifie
$\{\mu_t = \alpha _t,\; \forall \,t < \tau_n\}$ :
$$
e^{-\beta \tau_n^\alpha} \, m\rho^+
(\zeta_n^\alpha,Y_{(\tau_n^\alpha)^-}) \geq
 \mathbb{E} \left[-e^{-\beta
\tau_n}c(\zeta_n^\mu,Y_{(\tau_n^\mu)^-},\zeta_{n+1}^\mu,Y_{\tau_n^\mu})+
\mathbb{E}(k_{\tau_n^+}(\mu)|\,
\F_{\tau_n^\mu})\,\big|\,\G_{\tau_n^\alpha}\right].
$$
La tribu $\G_{\tau_n}$ \'etant une sous-tribu de $\F_{\tau_n}$,
$\forall\mu\in \{\mu_t = \alpha _t,\; \forall \,t < \tau_n\}$:
$$
e^{-\beta \tau_n^\alpha} \, m\rho^+
(\zeta_n^\alpha,Y_{(\tau_n^\alpha)^-})
 \geq  \mathbb{E}\left[-e^{-\beta
\tau_n}c(\zeta_n^\mu,Y_{(\tau_n^\mu)^-},\zeta_{n+1}^\mu,Y_{\tau_n^\mu})+
k_{\tau_n^+}(\mu)|\,\G_{\tau_n^\alpha}\right].$$
 Par suite, on a:
$$
e^{-\beta \tau_n^\alpha} \, m\rho^+
(\zeta_n^\alpha,Y_{(\tau_n^\alpha)^-})
 \geq \ESS{\{\mu_t = \alpha _t,\; \forall \,t < \tau_n\}} \mathbb{E}\left[-e^{-\beta
\tau_n}c(\zeta_n^\mu,Y_{(\tau_n^\mu)^-},\zeta_{n+1}^\mu,Y_{\tau_n^\mu})+
k_{\tau_n^+}(\mu)|\,\G_{\tau_n^\alpha}\right].$$
 On en déduit
donc l'in\'egalit\'e
$$e^{- \beta \tau_n^\alpha} \, m\rho^+ (\zeta_n^\alpha,Y_{(\tau_n^\alpha)^-}) \geq  W_{\tau_n}^\alpha.$$
Inversement, soit le  noyau
bor\'elien $r^* $ introduit dans la proposition \ref{prnoyaubor}. On applique l'expression (\ref{noyaubor})  au point $(i,x)= (\zeta_n^\alpha,Y_{(\tau_n^\alpha)^-}),$ il vient p.s.
 \begin{eqnarray*}
e^{-\beta \tau_n^\alpha} \,    m\rho^+
(\zeta_n^\alpha,Y_{(\tau_n^\alpha)^-}) &=& e^{-\beta \tau_n^\alpha}
\int_{{U}\times
\mathbb{R}}r ^*(\zeta_n^\alpha,Y_{(\tau_n^\alpha)^-};dy)\left(-c(\zeta_n^\alpha,Y_{(\tau_n^\alpha)^-},j,y)+\rho^+(j,y)\right).\nonumber\\
 \end{eqnarray*}
Le noyau bor\'elien $r^*$ est une loi
 conditionnelle de passage de la technologie $\zeta_n$ \`a $\zeta_{n+1}$ et de $Y_{\tau_n^-}$ \`a
 $Y_{\tau_n}$. Ainsi, il lui est associé une strat\'egie $\mu^{*}$ qui v\'erifie
$\{\mu_t^{*}= \alpha _t,\; \forall \,t < \tau_n\}$ telle que p.s.
\begin{eqnarray*}
e^{-\beta \tau_n^\alpha} \,   m\rho^+
(\zeta_n^\alpha,Y_{(\tau_n^\alpha)^-}) &= &  e^{-\beta \tau_n} \,
\mathbb{E}\left(-c(\zeta_n^{\mu^*},Y_{(\tau_n^{\mu^*})^-},\zeta_{n+1}^{\mu^*},Y_{\tau_n^{\mu^*}})+
    \rho^+(\zeta_{n+1}^{\mu^*},Y_{\tau_n^{\mu^*}})|\,\G_{\tau_n^{\mu^*}}\right).
\end{eqnarray*}
De plus, d'apr\`es l'\'egalit\'e (\ref{W+rho}),
 on a  p.s.:
$$
e^{-\beta \tau_n} \,  m\rho^+ (\zeta_n,Y_{\tau_n^-}) =
\mathbb{E}\left(-e^{-\beta \tau_n}
c(\zeta_n^{\mu^*},Y_{(\tau_n^{\mu^*})^-},\zeta_{n+1}^{\mu^*},Y_{\tau_n^{\mu^*}})
+ W_{\tau_n}^{\mu_*^+} |\,\G_{\tau_n^{\mu^*}}\right).$$
 En
appliquant l'in\'egalit\'e (\ref{egalite2}), on obtient l'in\'egalit\'e:
$$
e^{-\beta \tau_n} \,  m\rho^+ (\zeta_n,Y_{\tau_n^-}) \leq
W_{\tau_n}^{\mu^*} = \ESS{\{\nu_t = \mu^*_t, t < \tau_n \}}
\mathbb{E}(k_{\tau_n}(\nu)|\,\G_{\tau_n^\nu}).
 $$
 Parce que $\mu_t^* = \alpha _t,\; \forall \,t < \tau_n,$ les deux
 ensembles $\{\nu_t = \mu^*_t, t <
\tau_n \}\;et\;\{ \nu_t = \alpha_t, t < \tau_n \}$ co\"{\i}ncident.
Ainsi $$\ESS{\{\nu_t = \mu^*_t, t < \tau_n \}}
\mathbb{E}(k_{\tau_n}(\nu)|\,\G_{\tau_n^\nu})= \ESS{\{ \nu_t =
\alpha_t, t < \tau_n\}}
\mathbb{E}(k_{\tau_n}(\nu)|\,\G_{\tau_n^{\nu}}).$$ Par suite,
$$
e^{-\beta \tau_n} \, m\rho^+ (\zeta_n,Y_{\tau_n^-})) \leq
W_{\tau_n}^{\mu^*}= W_{\tau_n}^{\alpha}.
 $$
  Ainsi  on  obtient 
l'in\'egalit\'e inverse et donc l'\'egalit\'e (\ref{mrhW}).
\\
\\
2.  
Pour tout $n\geq 0,$ des \'egalit\'es (\ref{W rho}) et
    (\ref{promrhw}), il vient:
    $$m\rho^+(\zeta_n^\alpha, Y_{(\tau_n^\alpha)^-}) = \rho(\zeta_n^\alpha, Y_{(\tau_n^\alpha)^-})\quad \mbox{p.s}.$$

%En particulier, pour $n = 0$ et toute strat\'egie $\alpha$ qui
%d\'emarre avec $\xi_0 = i$  telle que $\tau_0=t$, $t>0,$ on a 
%p.s.
%$$m\rho^+(i, Y_{t-}) = \rho(i, Y_{t-}).$$
%On obtient   donc pour tout $ i\in {U}$ et tout $x \in E$:
%$$\rho(i,x) =m\rho^+(i,x).$$
\endproof

\begin{rem}\label{ensE}  Soit l'ensemble 
$ E = \{  (i,x) \in U \times \R:   ~ \rho(i,x) =m\rho^+(i,x)  \}.$
Alors cet ensemble contient l'ensemble $\{  (i,x) \in U \times \R:   ~   \exists ~  \alpha   ~  \exists ~ (\omega,t), Y_t (\omega) = x  \}.$
On verra plus tard que, suite à des propriétés topologiques, de fait $E=U \times \R.$
\end{rem}

\begin{pro}
\label{egalrho}
L'application $\rho^+$ satisfait à l'\'egalit\'e suivante:
\begin{equation}\label{egalrho1}
\rho^+(i,x) = \ESS{T>0, T \in
\underline{R}_{-1}}\mathbb{E}_{\{i,x\}}\left(\int_0^{T((i,x),.)}
e^{-\beta s}f(i,Y_s)\, ds + e^{-\beta T((i,x),.)}
m\rho^+(i,Y_{T^-((i,x),.)})\right),
\end{equation}
 o\`u  $\underline{R}_{-1}$ est l'ensemble des
applications mesurables~ $T$ de $(U  \times \R
\times \Omega, \mathcal{P}(U) \otimes
\mathcal{B}(\R ) \otimes\F)$ dans
$(\mathbb{R}_+,\mathcal{B}_{\mathbb{R}_+})$,
tel que pour  $(i,x)\in {U}\times\mathbb{R},$
$T((i,x),.)$ est  un $\G\mbox{-temps d'arrêt}$.
\end{pro}
\noindent\textbf{Preuve}
\\
1. On rappelle que
$$
\rho^+(i,x)= \ESS{\mu \in \A,
\,\tau_0>0}\mathbb{E}_{\{i,x\}}(k(\mu)).$$ 
 Soit  $T((i,x),.)\in \underline{R}_{-1}$ tel que $T((i,x),.)>0 $ p.s.  Pour toute strat\'egie $\mu$
qui d\'emarre avec la technologie $\xi_0 = i$ et telle que
$\tau_0^\mu= T((i,x),.)> 0$ p.s., on a la suite
d'\'egalit\'es:
{\small
\begin{eqnarray*}
k(\mu)& = &\mathbb{E}_{\{i,x\}}\Big[\int_0^{T((i,x),.)} e^{-\beta
s}f(i,Y_s^\mu)\, ds + \int_{T((i,x),.)}^{+\infty} e^{-\beta
s}f(i,Y_s^\mu)\, ds  \nonumber \\
&-& e^{-\beta T((i,x),.)} c(i,Y^\mu_{T^-((i,x),.)}, \zeta_1^\mu,
Y^\mu_{T((i,x),.)}) - \sum_{   0< T((i,x),.)<\tau_n <\tau} 
\,e^{-\beta \tau_n}
c(\zeta_n^\mu,Y^\mu_{(\tau_n)^-}, \zeta_{n+1}^\mu, Y^\mu_{\tau_n})\Big]\nonumber\\
 &=&  \mathbb{E}_{\{i,x\}}\big[ -e^{-\beta T((i,x),.)}
c(i,Y^\mu_{T^-((i,x),.)},
\zeta_1^\mu, Y^\mu_{T((i,x),.)})+\int_0^{T((i,x),.)} e^{-\beta s}f(i,Y^\mu_s)\, ds \nonumber\\
&+& \mathbb{E}\left(k_{T^+((i,x),.)} (\mu)|\;\F_
{T((i,x),.)}\right)\big].
\end{eqnarray*}
}
 D'o\`{u}, en prenant l'essentiel sup sur le dernier terme
{\small
$$
\rho^+(i,x) \geq \mathbb{E}_{\{i,x\}}\left(-  e^{-\beta T((i,x),.)}
c(i,Y^\mu_{T^-((i,x),.)}, \zeta_1^\mu, Y^\mu_{T((i,x),.)})
+\int_0^{T((i,x),.)} e^{-\beta s}f(i,Y^\mu_s)\, ds +
W^{\mu^+}_{T((i,x),.)}\right).$$}
 D'apr\`es
l'\'egalit\'e (\ref{W+rho}) prise en $T((i,x),.)$, nous avons:
{\small
\begin{eqnarray*}
 \rho^+(i,x) &\geq&
\mathbb{E}_{\{i,x\}}\big(- e^{-\beta T((i,x),.)}
c(i,Y^\mu_{T^-((i,x),.)}, \zeta_1^\mu, Y^\mu_{T((i,x),.)})
+\int_0^{T((i,x),.)} e^{-\beta s}f(i,Y^\mu_s)\, ds \nonumber\\&+&
e^{-\beta T((i,x),.)} \rho^+ (\zeta_1^\mu, Y^\mu_{T((i,x),.)})\big).
\end{eqnarray*}}
 Le processus $Y$ est $\G\mbox{-adapt\'e}$ et
$T((i,x),.)$ est  un $\G\mbox{-temps d'arr\^{e}t}$. D'o\`{u} en
conditionnant par la  tribu $\G_{T((i,x),.)}$, nous pouvons
\'ecrire:
\begin{eqnarray*} \rho^+(i,x)&\geq& \mathbb{E}_{\{i,x\}}\big[e^{-\beta T((i,x),.)}
\mathbb{E}\left[-c(i,Y_{T^-((i,x),.)}, \zeta_1^\mu,
Y^{\mu}_{T((i,x),.)} )+ \rho^+ (i,
Y^{\mu}_{T((i,x),.)})|\,\G_{T((i,x),.)}\right]\nonumber\\& +&
\int_0^{T((i,x),.)} e^{-\beta s}f(i,Y_s)\, ds \big].
\end{eqnarray*}
 On considère  la
strat\'egie $\mu$ qui d\'emarre avec $\xi_0=i$ et
$\tau_0=T((i,x),.)>0 $ p.s. avec  $T((i,x),.)\in \underline{R}_{-1}$ et de  noyau $r\in M,$ nous
avons:
$$\mathbb{E}\left[-c(i,Y^\mu_{T^-((i,x),.)}, \zeta_1^\mu, Y^{\mu}_{T((i,x),.)} )+ \rho^+ (i,
Y^{\mu}_{T((i,x),.)})|\,\G_{T((i,x),.)}\right]=$$
$$
 \int_{{U}\times \mathbb{R}} r(i,Y^{\mu}_{T^-((i,x),.)}, j,dy)(-c(i,Y^{\mu}_{T^-((i,x),.)}, j, y )+\rho^+ (j,y))
$$
o\`u, prenant l'essentiel sup pour $r\in M,$ on reconnait  la d\'efinition de 
$m\rho^+(i,Y_{T^-((i,x),.)}).$ Puis en passant \`a l'ess sup sur les $T((i,x),.)
\in \underline{R}_{-1}$  tels que $T((i,x),.)>0 $ p.s. il vient:
$$
 \rho^+(i,x)\geq   \ESS{T>0, \,T\in
\underline{R}_{-1}}\mathbb{E}_{\{i,x\}}\left(\int_0^{T((i,x),.)}
e^{-\beta s}f(i,Y_s)\, ds  + e^{-\beta T((i,x),.)}
m\rho^+(i,Y_{T^-((i,x),.)})\right).
$$
2.\, D'apr\`es l'\'egalit\'e (\ref{promrhw})
 prise en $n=0$ avec $\alpha=\mu$
déjà utilisée dans le 1. qui vérifie
$\xi_0=i,~\tau_0^\mu=T((i,x),.)>0 $ p.s.:
$$W^\mu_{T((i,x),.)} =e^{-\beta
T((i,x),.)} m\rho^+(i,Y_{T^-((i,x),.)}),$$
d'où,
$$
\mathbb{E}_{\{i,x\}}\left(\int_0^{T((i,x),.)} e^{-\beta s}f(i,Y_s)\,
ds + e^{-\beta T((i,x),.)} m\rho^+(i,Y_{T^-((i,x),.)})\right) =$$
$$
\mathbb{E}_{\{i,x\}}\left(\int_0^{T((i,x),.)} e^{-\beta s}f(i,Y_s)\,
ds + W^\mu_{T((i,x),.)} \right).$$
 Par
définition du gain maximal conditionnel, nous avons, pour la
strat\'egie $\mu$ qui d\'emarre avec $\xi_0=i$ et $\tau_0=T((i,x),.)>
0$ p.s.:
$$
\mathbb{E}_{\{i,x\}}\left[\int_0^{T((i,x),.)} e^{-\beta s}f(i,Y_s)\,
ds + W^\mu_{T((i,x),.)} \right)\geq \mathbb{E}_{\{i,x\}} \left(k
(\mu)- k_{T((i,x),.)}(\mu)+ \mathbb{E}
[k_{T(i,.)}(\mu)|\;\G_{T((i,x),.)} ]\right]
$$
L'expression  $k (\mu)- k_{T((i,x),.)}(\mu)$ étant
$\G_{T((i,x),.)}\mbox{-mesurable}$,  nous avons
{\small
$$
\mathbb{E}_{\{i,x\}} \left[k (\mu)- k_{T((i,x),.)}(\mu)+ \mathbb{E}
[k_{T((i,x),.)}(\mu)|\;\G_{T((i,x),.)} ]\right]=
\mathbb{E}_{\{i,x\}}(k (\mu)). $$ }
D'o\`{u}, pour toute strat\'egie
$\mu$ qui d\'emarre avec $\xi_0=i$ et  $\tau_0=T((i,x),.)> 0$ p.s. :
$$
\mathbb{E}_{\{i,x\}}\left(\int_0^{T((i,x),.)} e^{-\beta s}f(i,Y_s)\,
ds + W^\mu_{T((i,x),.)} \right)\geq\mathbb{E}_{\{i,x\}}(k (\mu)).
$$
Par d\'efinition de l'ess sup (cf. \cite{neveu1}), il en suit:
$$
\mathbb{E}_{\{i,x\}}\left(\int_0^{T((i,x),.)} e^{-\beta s}f(i,Y_s)\,
ds + W^\mu_{T((i,x),.)} \right)\geq \ESS{\mu \in \A,
\tau_0^\mu>0}\mathbb{E}_{\{i,x\}}(k (\mu)) = \rho^+(i,x).
$$
  D'o\`{u} l'\'egalit\'e (\ref{egalrho1}).
\endproof
Par suite, nous pouvons exprimer le crit\`ere d'optimalit\'e
\'etabli dans le th\'eor\`eme \ref{thm optimalit} \`a l'aide de
l'application $\rho^+$ qui, grâce à la proposition \ref{egalrho}, est ind\'ependante de
toute strat\'egie admissible.
\begin{thm}\label{stratoptimal}
Pour toute strat\'egie admissible $\alpha =(\tau_0,r)$, nous avons
les in\'egalit\'es suivantes p.s.:
\begin{eqnarray}
\mbox{Pour tout}\;  n \geq 0,\nonumber\\
m\rho^+(\zeta_n,Y_{\tau_n^-}) & \geq & \int_{{U}\times
\mathbb{R}}r(\zeta_n,Y_{\tau_n^-},j,dy)\left(-c(\zeta_n,Y_{\tau_n^-},j,y)+
\rho^+(j,y)\right).\label{egalite5}\\
\mbox{Pour tout}\;  n \geq -1,\nonumber\\
e^{-\beta \tau_n}\,\rho^+(\zeta_{n+1},Y_{\tau_n}) & \geq &
\mathbb{E}_{\{\zeta_{n+1},Y_{\tau_n}\}}\left(\int_{\tau_n}^{\tau_{n+1}
}e^{-\beta s}f(\zeta_{n+1},Y_s)\;ds+ e^{-\beta \tau_{n+1}
} m\rho^+(\zeta_{n+1},Y_{\tau_{n+1}^-}) \right).\nonumber\\
\label{egalite6}
\end{eqnarray}
En particulier, pour $n=-1$ il vient:
\begin{equation}
\rho^+(i,x)  \geq 
\mathbb{E}_{\{i,x\}}\left(\int_0^{\tau_0}e^{-\beta s}f(i,Y_s)\;ds+
e^{-\beta \tau_0}
m\rho^+(i,Y_{\tau_0^-})\right).\label{egalite4}
\end{equation}
 De plus, $\widehat{\alpha}$ est optimale si
et seulement si l'\'egalit\'e a lieu simultan\'ement dans
 (\ref{egalite5}) et (\ref{egalite6}).
\end{thm}

\noindent\textbf{Preuve:} 
%Soit une strat\'egie $\alpha \in \A$ .\\
% Gr\^{a}ce \`a la proposition \ref{egalrho} appliqu\'ee
%au temps  $T((i,x),.) = \tau_0^\alpha \in \underline{R}_{-1}$, on
%obtient:
%$$\rho^+(i,x) \geq
%\mathbb{E}_{\{i,x\}}\left(\int_0^{\tau_0^\alpha} e^{-\beta
%s}f(i,Y_s)\;ds+ e^{-\beta \tau_0^\alpha}
%m\rho^+(i,Y_{\tau_0^-})\right).$$
% Ainsi l'in\'egalit\'e
%(\ref{egalite4}) est v\'erifi\'ee.\\
%\\
%Quant \`a l
L'in\'egalit\'e (\ref{egalite5}), pour tout $n\geq0$, 
est tir\'ee
 de l'in\'egalit\'e (\ref{egalite2}) du
th\'eor\`eme \ref{thm optimalit} en rempla\c{c}ant
$W_{\tau_n}^\alpha$ par $e^{-\beta \tau_n} \, m\rho^+
(\zeta_n,Y_{\tau_n^-})$ (\'egalit\'e (\ref{promrhw})) et
$W_{\tau_n}^{\alpha^+}$ par $e^{-\beta \tau_n}\rho^+
(\zeta_{n+1},Y_{\tau_n})$
(\'egalit\'e (\ref{W+rho})).\\
\\
Puis, pour tout $n\geq-1$, l'in\'egalit\'e (\ref{egalite6})
s'obtient imm\'ediatement \`a partir de l'in\'egalit\'e
(\ref{egalite3}) en rempla\c{c}ant $W_{\tau_n}^{\alpha^+}$ par
$e^{-\beta \tau_n}\rho^+ (\zeta_{n+1},Y_{\tau_n})$ (\'egalit\'e
(\ref{W+rho})) et $W_{\tau_{n+1}}^\alpha$ par $e^{-\beta \tau_{n+1}}
\, m\rho^+
(\zeta_{n+1},Y_{\tau_{n+1}^-})$ (\'egalit\'e (\ref{promrhw})).\\
\\
Ensuite, supposons que la strat\'egie $\widehat{\alpha}$ est
optimale : dans ce cas,
 les \'egalit\'es (\ref{egalite1}),(\ref{egalite2}) et (\ref{egalite3})  du th\'eor\`eme \ref{thm
optimalit} sont v\'erifi\'ees . Ce qui entra\^{\i}ne que
(\ref{egalite5}) et (\ref{egalite6}) sont des \'egalit\'es.
 L'égalité (\ref{egalite4}) c'est exactement l'égalité
(\ref{egalite6}) en $n=-1$.\\
\\
 Inversement, supposons que les in\'egalit\'es
(\ref{egalite5}) et (\ref{egalite6}) sont des \'egalit\'es.\\
L'égalité (\ref{egalite2}) s'obtient, pour tout $n\geq 0,$ de
l'égalité (\ref{egalite5}) en rempla\c{c}ant $e^{-\beta \tau_n} \,
m\rho^+ (\zeta_n,Y_{\tau_n^-})$ par $W_{\tau_n}^\alpha$ (\'egalit\'e
(\ref{promrhw})) et $e^{-\beta \tau_n}\rho^+
(\zeta_{n+1},Y_{\tau_n})$ par $W_{\tau_n}^{\alpha^+}$
(\'egalit\'e (\ref{W+rho})).\\
\\
Quant à l'égalité (\ref{egalite3}), elle est tirée, pour tout $n
\geq-1$, de l'égalité (\ref{egalite6}) en rempla\c{c}ant $e^{-\beta
\tau_n}\rho^+ (\zeta_{n+1},Y_{\tau_n})$  par
$W_{\tau_n}^{\alpha^+}$(\'egalit\'e (\ref{W+rho})) et
 $e^{-\beta \tau_{n+1}} \, m\rho^+
(\zeta_{n+1},Y_{\tau_{n+1}^-})$ par $W_{\tau_{n+1}}^\alpha$ (\'egalit\'e (\ref{promrhw})).\\
\\
 Par
suite les \'egalit\'es  (\ref{egalite2}) et (\ref{egalite3})  donc
(\ref{egalite1}) sont v\'erifi\'ees. Ainsi,  d'apr\`es le
th\'eor\`eme \ref{thm optimalit}, la strat\'egie $\alpha$ est
optimale.

\endproof

\section{ Une stratégie optimale} \label{section 4}

De la remarque \ref{rhomesurable}.2, il vient:
$$m\rho^+(i,x)  = \rho^+(i,x) \vee m^\ast\rho^+(i, x) $$
et sur l'ensemble  
$$ E = \{(i,x)\in U\times\R,~\rho(i,x)  = m\rho^+(i,x)\},$$ 
 nous obtenons:
$$\rho(i,x)  = \rho^+(i,x)  \vee m^\ast\rho^+(i,x) .$$
On en déduit que $\rho(i,x)  \geq m^\ast \rho^+(i,x) $, pour tout
$(i,x)   \in  E.$\\

%Notons que $E = C \cup I,$ où $C$ et $I$ sont deux ensembles définis comme suit:}

\begin{defn}
Nous appelons ensemble optimal de continuation le sous-ensemble de
${U}\times\mathbb{R}$:
$$ C = \big\{(i,x) \in E: \rho(i,x) > m^\ast \rho^+(i,x)\big\}.$$
\end{defn}
Supposons que $(i,x) \in C$, il n'est pas intéressant de changer
de stratégie à cet instant. En effet $\rho(i,x)$, le meilleur gain
en partant de $x$ est strictement supérieur au meilleur gain quand
on réalise au moins un changement. Il faut donc laisser le
syst\`eme \'evoluer librement, on voit que le premier instant de
changement de technologie est celui où le syst\`eme atteint la fronti\`ere de $C$.\\

Nous noterons par I son compl\'ementaire  dans $E$ qui est \textbf{l'ensemble
d'impulsion}. Il est défini comme suit:
 $$ I = \{(i,x)\in E: \rho(i,x) = m^\ast \rho^+(i,x)\}.$$
Notons que $E = C \cup I.$ 
\\

Pour construire la suite de temps d'impulsion, on introduit le temps suivant:
\begin{defn} 
Nous définissons le temps:
 $$T^\ast ((i,x),.) = \left\{%
\begin{array}{ll}
   \inf\{ t > 0:  (\xi_t,Y_t) \in I\}.
 \\
    +\infty \qquad \mbox{si l'ensemble est vide.} \\
\end{array}%
\right.
$$
 \end{defn}

\begin{pro}
Sous les hypothèses 2, 
 la fonction   $\rho^+$ est semi-continue inférieurement (s.c.i).
\end{pro}

\noindent\textbf{Preuve:}  De la proposition \ref{recurence}, la fonction $\rho^+$ peut être écrite sous la forme suivante:
 {\small \begin{eqnarray*}
\rho^+(i,x) &=&   \ESS{\mu \in \A, \tau_0^\mu >0} \left[ \sum_{l \geq 0} F_1(l,i,x ) -\sum_{l> 0}F_3(l,i,x)  + \mathbb{E}_{i,x}\left[\int_ {\tau_{N(\omega)}}^{+\infty}  e^{-\beta s} f(\zeta_{N(\omega)},Y_s)ds \right]\right]\nonumber\\
&=& \ESS{\mu \in \A, \tau_0^\mu >0}  \left[ \sum_{l \geq 0}  \left[F_1(l,i,x ) -F_3(l+1,i,x) \right] + \mathbb{E}_{i,x}\left[\int_ {\tau_{N(\omega)}}^{+\infty}  e^{-\beta s} f(\zeta_{N(\omega)},Y_s)ds \right]\right]. 
\end{eqnarray*}}
 (i) Tout d'abord,   prouvons la continuité de la fonction  $x \longmapsto \mathbb{E}_{i,x}\big[\int_ {\tau_{N(\omega)}}^{+\infty}  e^{-\beta s} f(\zeta_{N(\omega)},Y_s)ds \big].$ D'une part, cette fonction peut être écrite sous la forme suivante:
$$\mathbb{E}_{i,x}\left[\int_ {\tau_{N(\omega)}}^{+\infty}  e^{-\beta s} f(\zeta_{N(\omega)},Y_s)ds \right] = \sum_{k \geq 0} 
\mathbb{E}_{i,x} \left[\int_ {\tau_k}^{+\infty}  e^{-\beta s} f(\zeta_{N(\omega)},Y_s)ds ~\1_{\{N=k\}} \right].$$
La fonction $f$ étant bornée positive (hypothèse 2 (ii)), la dernière expression peut être majorée par:
$$
\mathbb{E}_{i,x} \left[\int_ {\tau_k}^{+\infty}  e^{-\beta s} f(\zeta_{N(\omega)},Y_s)ds ~\1_{\{N=k\}} \right] \leq \mathbb{E}_{i,x} \left[\int_ {\tau_k}^{+\infty}  e^{-\beta s}  \|f\|_{\infty} ds ~\1_{\{N=k\}} \right] \nonumber\\
$$
Sur l'ensemble $\{N =k\}$, $\tau_{j-1}<\tau_{j} ~ \forall j\leq k,$  nous obtenons ainsi $\tau_k \geq k.$ D'où, 
$$
\mathbb{E}_{i,x} \left[\int_ {\tau_k}^{+\infty}  e^{-\beta s} f(\zeta_{N(\omega)},Y_s)ds ~\1_{\{N=k\}} \right] \leq  \|f\|_{\infty}e^{-\beta k},
$$
où $\|f\|_{\infty}e^{-\beta k}$ est une série à termes positives convergente.  \\

D'autre part, rappelons que le processus $Y$ est un processus de Markov homogène qui peut être écrit sous la forme suivante:
$$Y_t = Y_{\tau_k} + \int_{\tau_k}^{t}( b(\xi_u,Y_u)du + \sigma \left(\xi_u,Y_u) dW_u \right).$$
La fonction $x \longmapsto \mathbb{E}_{i,x} \left[\int_ {\tau_k}^{+\infty}  e^{-\beta s} f(\zeta_{N(\omega)},Y_s)ds ~\1_{\{N=k\}} \right]$ est égale à:
$$
\mathbb{E} \left[\mathbb{E} \left[\int_ {\tau_k}^{+\infty}  e^{-\beta s} f(\zeta_{N(\omega)},Y_{\tau_k} + \int_{\tau_k}^{s}( b(\xi_u,Y_u)du + \sigma \left(\xi_u,Y_u) dW_u \right))ds ~\1_{\{N=k\}}   \bigg| \F_{\tau_n}\right] \right].
$$
Notons que la loi de la fonction $$x \longmapsto \mathbb{E} \left[\int_ {\tau_k}^{+\infty}  e^{-\beta s} f(\zeta_{N(\omega)},Y_{\tau_k} + \int_{\tau_k}^{s}( b(\xi_u,Y_u)du + \sigma \left(\xi_u,Y_u) dW_u \right))ds ~\1_{\{N=k\}}   \bigg| \F_{\tau_n}\right]$$ est une fonction mesurable en $\xi_k$ et $Y_{\tau_k-}.$
Ensuite, en utilisant l'assertion (\ref{rn}) et la définition \ref{strategadm}.2, nous constatons que la loi de $(\xi_s,Y_s )\big | (\xi_0=i, Y_0 =x)$ est égale à: 
$$\prod_{i=1}^k \mbox{ loi de}~(\xi_i,Y_{\tau_i} \big | \F_{\tau_{i-1}} ) \times \mbox{ loi de} ~ (\xi_0,Y_{\tau_0}) \big |  (i,x),$$
elle même égale à $  \prod_{i=1}^k \mbox{ loi de}~(\xi_i,Y_{\tau_i}) \big | (\xi_{i^-},Y_{\tau_i^-} ) \times \mbox{ loi de} ~ (\xi_0,Y_{\tau_0} )\big |  (i,x).$\\
Sous l'hypothèse 2 (i) et la définition \ref{strategadm}, la loi  produit ci-dessus  est continue en $x$.\\
Enfin, la fonction  $x \longmapsto \mathbb{E}_{i,x}\big[\int_ {\tau_{N(\omega)}}^{+\infty}  e^{-\beta s} f(\zeta_{N(\omega)},Y_s)ds \big]$ est continue. \\
\\
(ii) Ensuite, il suffit  de prouver que $F_1$ et $F_3$ sont continues $\forall l$ et  qu'il existe   une série réelle $(\eps_l)$ à termes positifs convergente telle que $\forall l, \forall x,$
$$\big| F_1(l,i,x) - F_3(l+1,i,x) \big| \leq \eps_l.$$
 Procédons par récurrence pour montrer que  la fonction $F_1(l,i,x)$ est majorée par:
\begin{equation}\label{recmajor}
\big| F_1(l,i,x)\big| \leq \|f\|_{\infty}\left( \mathbb{E}_{\{i,x\}}\left[\1_{\{\tau_0>0\}}  e^{-\beta\,\tau_0}  \right]\right)^l.
\end{equation}
Pour $l=0,$ nous avons:
$$
F_1(0,i,x)= \mathbb{E}_{i,x}\left[\1_{\{\tau_0>0\}} \int_{0}^{\tau_0}e^{-\beta
 t}f\big(i,
Y^x_t)\,dt\right] .$$
De l'hypothèse 2 (ii), la fonction $f$ étant bornée positive, $x\mapsto F_1(0,i,x)$ est majorée par:
$$\big| F_1(l,i,x)\big|\leq \mathbb{P}_{i,x}(\tau_0>0)\|f\|_{\infty} \leq \|f\|_{\infty}.$$
Supposons que l'hypothèse de récurrence (\ref{recmajor}) est vraie  jusqu'à $l-1$ et montrons qu'elle est vraie à l'ordre $l$.  Pour $l>0,$
\begin{eqnarray*}
\big|F_1(l,i,x)\big|&=&\big|\mathbb{E}_{\{i,x\}}\left[\1_{\{\tau_0>0\}}  e^{-\beta\,\tau_0}\, \int_{{U}\times \mathbb{R}} F_1(l-1,j,y)
r (i,Y^x_{\tau_0^-},j,dy)\right] \big|\nonumber\\
&\leq& \mathbb{E}_{\{i,x\}}\left[\1_{\{\tau_0>0\}}  e^{-\beta\,\tau_0}  \right] \|f\|_{\infty}\left( \mathbb{E}_{\{i,x\}}\left[\1_{\{\tau_0>0\}}  e^{-\beta\,\tau_0}  \right]\right)^{l-1}\nonumber\\
&=&\|f\|_{\infty}\left( \mathbb{E}_{\{i,x\}}\left[\1_{\{\tau_0>0\}}  e^{-\beta\,\tau_0}  \right]\right)^{l}. 
\end{eqnarray*}
Ainsi l'hypothèse de récurrence est vérifiée $\forall l.$\\
\\
Procédons d'une manière analogue pour prouver que $F_3(l+1,i,x)$ est majorée par 
\begin{equation}\label{recmajor2}
\big| F_3(l+1,i,x)\big| \leq \|c\|_{\infty}\left(   \mathbb{E}_{\{i,x\}}\left[\1_{\{\tau_0>0\}}  e^{-\beta\,\tau_0}  \right]\right)^{l+1}.
\end{equation}
Pour $l =0,$ du corollaire \ref{relationrec}, nous avons:
$$F_3(1,i,x)= \mathbb{E}_{i,x}\left[\1_{\{\tau_0>0\}} e^{-\beta \tau_0}\int_{{U}\times \mathbb{R}}  \mathbb{P}_{j,y} (\tau_0 >0) c(i, Y_{\tau_0^-}^x,j,y) 
r(i,Y_{\tau_0^-}^x,j,dy) \right]$$
De l'hypothèse 2 (ii), la fonction $c$ est  bornée, la fonction $F_3(1,i,x)$ est donc majorée par
$$\big| F_3(1,i,x)\big| \leq \|c\|_{\infty}\mathbb{E}_{\{i,x\}}\left[\1_{\{\tau_0>0\}}  e^{-\beta\,\tau_0}  \right].$$
Supposons que l'hypothèse de récurrence (\ref{recmajor2}) est vraie  jusqu'à $l$ et montrons qu'elle est vraie à l'ordre $l+1$.  Pour $l>1,$
\begin{eqnarray*}
\big|F_3(l+1,i,x)\big|&=&\big| \mathbb{E}_{i,x}\left[\1_{\{\tau_0>0\}} e^{-\beta \tau_0}\int_{{U}\times \mathbb{R}}  F_3(l,j,y)
r(i,Y_{\tau_0^-}^x,j,dy) \right] \big|\nonumber\\
&\leq& \mathbb{E}_{\{i,x\}}\left[\1_{\{\tau_0>0\}}  e^{-\beta\,\tau_0}  \right] \|c\|_{\infty}\left( \mathbb{E}_{\{i,x\}}\left[\1_{\{\tau_0>0\}}  e^{-\beta\,\tau_0}  \right]\right)^{l}\nonumber\\
&=&\|c\|_{\infty}\left( \mathbb{E}_{\{i,x\}}\left[\1_{\{\tau_0>0\}}  e^{-\beta\,\tau_0}  \right]\right)^{l+1}. 
\end{eqnarray*}
Ainsi l'hypothèse de récurrence est vérifiée $\forall l.$ 
Par suite, nous obtenons, $\forall l,$
$$\big| F_1(l,i,x) -F_3(l,i,x)  \big| \leq (\|f\|_{\infty}+\|c\|_{\infty})   \left( \mathbb{E}_{\{i,x\}}\left[\1_{\{\tau_0>0\}}  e^{-\beta\,\tau_0}  \right]\right)^l,$$
où $\eps_l = (\|f\|_{\infty}+\|c\|_{\infty})   \left( \mathbb{E}_{\{i,x\}}\left[\1_{\{\tau_0>0\}}  e^{-\beta\,\tau_0}  \right]\right)^l$ 
est une série à termes positifs convergente car  $ \mathbb{E}_{\{i,x\}}\left[\1_{\{\tau_0>0\}}  e^{-\beta\,\tau_0}  \right]<1.$ En effet, si $ \mathbb{E}_{\{i,x\}}\left[\1_{\{\tau_0>0\}}  e^{-\beta\,\tau_0}  \right]=1,$ nous obtenons $\tau_0 =0$ et 
$ \1_{\tau_0 > 0}=1$, ce qui est contradictoire.\\
\\
(iii)  Enfin, il nous reste à prouver que $F_1$ et $F_3$ sont continues $\forall l:$ \\
Pour $l =0,$ nous avons: 
$$
F_1(0,i,x)= \mathbb{E}_{i,x}\left[\1_{\{\tau_0>0\}} \int_{0}^{\tau_0}e^{-\beta
 t}f\big(i,x+b(i)t+\sigma(i)W_t)\,dt\right] .$$
De l'hypothèse 2 (ii), la fonction $f$ étant bornée positive, $x\mapsto F_1(0,i,x)$ est continue. \\
Pour $l>0,$ 
la fonction $F_1(l,i,x)$ peut être écrite sous la forme suivante:
$$F_1(l,i,x)=\int_0^{+\infty} e^{-\beta s}\left(  \int_{{U}\times \mathbb{R}} F_1(l-1,j,y)
r (i,x + b(i) s+ \sigma(i) W_s;j,dy)\right)\mathbb{P}_{\tau_0}(ds\big|i,x ).$$
L'hypothèse 2 (i)  (iii) montre  la continuité de la fonction $x\mapsto F_1(l,i,x)$, $\forall l.$\\
\\
De même, pour $l =1,$  la fonction $x\mapsto F_3(1,i,x)$ est donnée par: 
 {\footnotesize
$$ \int_0^{+\infty} e^{-\beta s}\left( \int_{{U}\times \mathbb{R}}  \mathbb{P}_{j,y} (\tau_0 >0) c(i,x + b(i) s+ \sigma(i) W_s,j,y) 
r(i,x + b(i) s+ \sigma(i) W_s,j,dy)\right)\mathbb{P}_{\tau_0}(ds\big|i,x ),$$}
elle est  est continue   car, d'après l'hypothèse 2 (ii),  $x\mapsto c(i,x,..)$ est bornée continue et   $x \mapsto  r(i,x;.,.)$ est continue.\\
Pour $l>1,$
$$F_3(l,i,x)=\int_0^{+\infty} e^{-\beta s}\left( \int_{{U}\times \mathbb{R}}  F_3(l-1,j,y)
r(i,x + b(i) s+ \sigma(i) W_s,j,dy) \right)\mathbb{P}_{\tau_0}(ds\big|i,x ).$$
Ainsi,  l'hypothèse 2(i)  (iii)  induit la continuité de la fonction $x\mapsto F_3(l,i,x)$, $\forall l.$

\endproof
 
\begin{rem}
On peut montrer par une preuve analogue que la fonction $\rho$ est elle aussi 
semi-continue inférieurement. Rappelant que $m\rho^+$ est continue (proposition \ref{opcont})
$\rho-m\rho^+$ est semi-continue inférieurement et l'ensemble $E$ est un fermé de $U\times\R.$
La remarque \ref{ensE} dit qu'il contient l'ensemble 
$\{  (i,x) \in U \times \R:   ~   \exists ~  \alpha   ~  \exists ~ (\omega,t), Y_t (\omega) = x  \}$:
par continuité, il contient aussi tout $(i,x)$ de $U \times \R,$ d'où l'égalité partout
$\rho=m\rho^+$. 
\end{rem}
 
\begin{cor}
Sous les hypothèses 1, 2 et 3, l'ensemble $I$ est fermé et $C$ est ouvert. 
\end{cor}

\noindent\textbf{Preuve:} 
  La proposition \ref{opcont} donne que  la fonction $m \rho^+$ est continue donc  la fonction  $ \rho$ est continue.
 Ensuite,  la fonction $m^\ast \rho^+$  est semi-continue inférieurement comme supremum de fonction semi-continue inférieurement. D'où, $(\rho-m^\ast\rho^+)$ est semi-continue inférieurement
 ce qui montre  que $I$ est fermé et $C$ est ouvert.  
\endproof

Cette propriété topologique des ensembles $I$ et $C$ permet de montrer le lemme suivant:

 \begin{lem}\label{CI}
 $\forall \omega \in \Omega,$ $ T^\ast ((i,x),\omega)  > 0 $
si et seulement si $(i,x) \in C.$  C'est à dire en d'autres termes que 
$(i,x) \in I$  si et seulement si $T^\ast ((i,x),\omega)  = 0$.
 \end{lem}

 \noindent\textbf{Preuve:} 
Supposons que $ T^\ast ((i,x),\omega)  > 0.$ $T^\ast ((i,x),\omega) $ est le temps de sortie d'un processus  continu (donc continu à gauche)  et adapté de l'ouvert $C$.  Par conséquent, $(\xi_t(\omega),Y_t(\omega)) \in C , ~\forall t \in [0,T^\ast ((i,x),\omega)[$, en particulier $(i,x) \in C.$\\
 \\
 Inversement, soit $(i,x) \in C.$ Il existe une suite de temps d'arrêt $(T_n((i,x),.))_n$  qui décroît vers $ T^\ast ((i,x),.)$ telle que, $\forall \omega \in \Omega,$ 
 $(\xi_{T_n ((i,x),\omega)},Y_{T_n ((i,x),\omega)}) \in I.$ L'ensemble $I$ étant fermé, $(\xi_{T^\ast ((i,x),\omega)},Y_{T^\ast ((i,x),\omega)}) \in I.$
 Or $I \cap C = \emptyset,$ d'où $ T^\ast ((i,x),\omega)  > 0.$
 
  \endproof

\begin{lem}\label{propT}
1. $T^\ast$ v\'erifie pour toute loi $\nu$ sur $({U}\times
\mathbb{R}\times
\Omega,\mathcal{P}({U})\times\mathcal{B}_{\mathbb{R}}\times
\F)$:
\begin{equation}\label{rho1}
 (\xi_{T^{\ast-}((i,x),.)},Y^x_{T^{\ast-}((i,x),.)})\in I \quad \nu
\;\mbox{p.s.}
\end{equation}

2. $T^\ast$ est une application mesurable  de $({U}\times
\mathbb{R}\times
\Omega,\mathcal{P}({U})\times\mathcal{B}_{\mathbb{R}}\times
\F)$ dans
$(\mathbb{R}_+,\mathcal{B}_{\mathbb{R}_+})$
telle que pour tout $(i,x)$, $T^\ast((i,x),.)$  est  un
$\G\mbox{-temps d'arr\^{e}t}$.\\
\\
3.   Pour tout $(i,x)\in C$:
\begin{equation}\label{T}
\rho^+(i,x) = \mathbb{E}_{\{i,x\}}\left(\int_0^{T^\ast((i,x),.)}
e^{-\beta s}f(i,Y_s)\, ds + e^{-\beta T^\ast((i,x),.)}
m\rho^+(i,Y_{T^{\ast-}((i,x),.)})\right).
\end{equation}
\end{lem}

\noindent\textbf{Preuve} \\
1. Le temps $T^\ast((i,x),.)$ est le temps d'entrée du processus continu $(i,Y^{i,x}_{t})$ dans le borélien fermé $I$. Ainsi, par continuité et définition du temps d'entrée, 
$(i,Y^{i,x}_{T^\ast((i,x),.)})\in I~\nu$ presque sûrement. Ce temps sera un temps de saut, mais 
avant le saut le processus continu coincide avec sa limite à gauche, d'où
  pour toute
loi $\nu$ sur$({U}\times \mathbb{R}\times
\Omega,\mathcal{P}({U})\times\mathcal{B}_{\mathbb{R}}\times
\F)$:
   $$(\xi_{T^{\ast-}((i,x),.)},Y^x_{T^{\ast-}((i,x),.)})\in I \qquad \nu \;\mbox{p.s.}$$
D'o\`{u} l'\'egalit\'e (\ref{rho1}) est vérifiée.\\
\\
2. Soient $\tilde{\Omega} = {U}\times
\mathbb{R}\times \Omega$ et $\tilde{\G}_t = \mathcal{P}({U})\times\mathcal{B}_{\mathbb{R}}\times
\G_t$. Pour tout $\tilde{\omega} \in \tilde{\Omega}, $ notons $X_t(\tilde{\omega}) = (\xi_{t^-}(\omega),Y^{i,x}_{t^-}(\omega)).$ D'après 1,  le temps $T^\ast$ peut être écrit sous la forme suivante:
$$T^\ast (\tilde{\omega}) = \inf\{t, X_t(\tilde{\omega}) \in I\}.$$
Le processus $X$ étant  continu  et adapté, il est  optionnel pour la filtration $\tilde{\G}$ (cf.  proposition 1.24 de \cite[p. 6]{jacod}). De plus, le temps $T^\ast$ est le temps d'entrée de $X_t$ dans l'ensemble borélien $I$. D'où, c'est un temps optionnel et par suite, c'est un  $\tilde{\G}$-temps d'arr\^{e}t 
(cf.  corollaire 2.4 de \cite[p. 6]{karatzas} ). 
Par conséquent, 
 $$\left\{\tilde{\omega}\in \tilde{\Omega}:~~ \{T^\ast(\tilde{\omega}) >
t\}\right\}\in\mathcal{P}({U})\times\mathcal{B}_{\mathbb{R}}\times
\G_t.$$ 
Ce qui entra\^{\i}ne que $T^\ast$ est une
application mesurable sur $({U}\times
\mathbb{R}\times\Omega,\mathcal{P}({U})\times\mathcal{B}_{\mathbb{R}}\times
\G_t)$.\\
Il reste à prouver que   pour tout $(i,x)$ fixé, $T^\ast((i,x),.)$ est un
$\G$-temps d'arr\^{e}t.\\ D'après 1, le temps $T^\ast$ peut être écrit sous la forme suivante:
$$T^\ast ((i,x), \omega) = \inf\{t, X_t((i,x), \omega) \in I\}.$$
Le processus $X$ étant continu  et adapté, il est  optionnel pour la filtration $\tilde{\G}$.  De plus, le temps $T^\ast$ est le temps d'entrée de $X_t$ dans l'ensemble borélien $I$. D'où, c'est un temps optionnel et par suite, c'est un  $\G$-temps d'arr\^{e}t. 
Par conséquent,  pour $(i,x)$ fixé, nous obtenons:
 $$\left\{\omega\in \Omega:~T^\ast((i,x), \omega) >
t\right\}\in
\G_t.~p.s.$$ 
3.  Soit $(i,x) \in C$,  alors $T^\ast((i,x),.)$ est un $\G$-temps d'arrêt strictement positif appartenant à  $\underline{R}_{-1}$.
D'o\`{u}, d'apr\`es la proposition  \ref{egalrho},
$$
\rho^+(i,x) \geq \mathbb{E}_{\{i,x\}}\left(\int_0^{T^\ast((i,x),.)}
e^{-\beta s}f(i,Y_s)\, ds + e^{-\beta T^\ast((i,x),.)}
m\rho^+(i,Y_{T^{\ast-}((i,x),.)})\right).
$$
Pour  toute strat\'egie $\mu \in \A$
 telle que $\tau_0^\mu=T^\ast((i,x),.)$, utilisant l'\'egalit\'e (\ref{promrhw}) il vient:
$$
\rho^+(i,x) \geq \mathbb{E}_{\{i,x\}}\left(\int_0^{T^\ast((i,x),.)}
e^{-\beta s}f(i,Y_s)\, ds  + W^\mu_{T^\ast((i,x),.)} \right).$$ Par
d\'efinition du gain maximal conditionnel, nous avons, pour toute
strat\'egie $\mu$ qui d\'emarre avec $\xi_0=i$ et
$\tau_0^\mu=T^\ast((i,x),.)> 0,$

\begin{eqnarray*}
\mathbb{E}_{\{i,x\}}\left(\int_0^{T^\ast((i,x),.)} e^{-\beta
s}f(i,Y_s)\, ds + W^\mu_{T^\ast((i,x),.)} \right)&\geq&
\mathbb{E}_{\{i,x\}} \big(k (\mu)-
k_{T^\ast((i,x),.)}(\mu)\nonumber\\
&+& \mathbb{E} [k_{T^\ast((i,x),.)}(\mu)|\;\G_{T^\ast((i,x),.)}
]\big).
\end{eqnarray*}
 L'expression  $k (\mu)-
k_{T^\ast((i,x),.)}(\mu)$ \'etant
$\G_{T^\ast((i,x),.)}\mbox{-mesurable}$,  nous avons
$$
\mathbb{E}_{\{i,x\}} \left(k (\mu)- k_{T^\ast((i,x),.)}(\mu)+
\mathbb{E} [k_{T^\ast((i,x),.)}(\mu)|\;\G_{T^\ast((i,x),.)}
]\right)= \mathbb{E}_{\{i,x\}}(k (\mu)). $$ D'o\`{u}, pour toute
strat\'egie $\mu$ qui d\'emarre avec $\xi_0=i$ et
$\tau_0=T^\ast((i,x),.)> 0,$ nous avons
$$
\mathbb{E}_{\{i,x\}}\left(\int_0^{T^\ast((i,x),.)} e^{-\beta
s}f(i,Y_s)\, ds + W^\mu_{T^\ast((i,x),.)}
\right)\geq\mathbb{E}_{\{i,x\}}(k (\mu)).
$$
Par d\'efinition de l'ess sup (cf. \cite{neveu1}), il en suit:
$$
\mathbb{E}_{\{i,x\}}\left(\int_0^{T^\ast((i,x),.)} e^{-\beta
s}f(i,Y_s)\, ds + W^\mu_{T^\ast((i,x),.)} \right)\geq \ESS{\mu
\in \A, \tau_0>0}\mathbb{E}_{\{i,x\}}(k
(\mu)) = \rho^+(i,x),
$$ d'o\`{u}
l'\'egalit\'e. Ce qui entraîne (\ref{T}). 
\endproof

%\begin{lem}
%Sous les hypothèses  1 et  2 (iii),  le noyau borélien  $r^\ast  $ introduit dans  la proposition \ref{prnoyaubor} vérifie,  pour tout $(i,x)$ de I
%\begin{equation}
%\label{lemme}
%  m\rho^+(i,x)=   m^\ast\rho^+(i,x) =
%  \int_{{U}\times\mathbb{R}}r^\ast
%(i,x;j,dy)\,\left(-c(i,x,j,y) + \rho^+(j,y)\right)\qquad \nu
%\;\mbox{p.s.}
%\end{equation}
%{\color{red} non, c'est faux}
%\end{lem}
%
% \noindent \textbf{Preuve } 
% Les hypothèses  1  et  2 (iii)  permettent d'appliquer la proposition \ref{prnoyaubor}: le noyau borélien $r^\ast  $ vérifie:
% $$m \rho^+ (i,x) = \int_{{U}\times \mathbb{ R}} r^* (i,x;j,dy)\left(-c(i,x;j,y) +\rho^+(j,y)\right).$$
%Pour tout $(i,x) \in I,$ {\color{red} non, c'est faux, seulement dans $E$ cf. mes remarques
%ci-dessus,} nous avons
%$\rho(i,x)=m^\ast\rho ^+(i,x).$ De plus, par la proposition
%\ref{mrhorho}, nous avons $\rho(i,x)=m\rho^+(i,x)$. D'où
%$$\forall (i,x) \in I,~~  m \rho^+ (i,x) = m^\ast\rho ^+(i,x).$$
%\endproof

 \begin{defn}\label{stroptimal}
Posons $\mathbb{P}_{\{i,x\}}$ p.s:
$$
\widehat{\tau}_0:=
 \left\{%
\begin{array}{ll}
    T^\ast((i,x),. ) & \hbox{si} \quad  T^\ast((i,x),. )>0
\\
+\infty &\mbox{ si } T^\ast((i,x),. )=0.
\end{array}%
\right.  $$
En $\widehat{\tau}_0$, la loi du   couple  $(\xi_1, Y_{\widehat{\tau}_0}) $ sachant $\G_{\widehat{\tau}_0}$ est donnée par 
 $r^\ast(\xi_0, Y_{\widehat{\tau}_0^-}, ., . ).$ 
\end{defn}

On note $N(\omega) = \inf\{n, T^*((\zeta_{n+1}, Y_{\widehat{\tau}_n}),. )=0\}$. 
Sur l'ensemble $\{\omega, N(\omega) < \infty\}$, nous avons  $\widehat{\tau}_{k} = \widehat{\tau}_{N(\omega)}$, $\forall k \geq N(\omega)$
 et à partir de cet instant la politique proposée est de continuer sans plus
d'impulsion, car on est en un point où $T^\ast(\zeta_k,Y_{\tau_k^-})=0.$
\\
Rappelons que la suite de la stratégie après $\tau_0$ est entièrement définie
par la construction explicitée par   la définition \ref{strategadm}.

\begin{thm} 
  Sous les hypothèses 1, 2 et 3,  la famille  
$\widehat{\alpha} =  (\widehat{\tau}_{0},r^\ast)$ définie par la définition \ref{stroptimal} est une stratégie  {\bf  admissible optimale } vérifiant     
  la propriété suivante:
\begin{equation}\label{propadmis}
\mathbb{P}(\tau_n^{\widehat{\alpha}} < \tau^{\widehat{\alpha}}  < +\infty,\; \forall n) = 0,
\end{equation}
 c'est à dire que
 soit $\widehat{\tau} = + \infty$ et dans ce cas la suite $(\widehat{\tau}_n)$ est strictement croissante, 
  soit il existe une v.a.  à valeurs entières $N$ telle que
$\widehat{\tau}_{k} = \widehat{\tau}_{N(\omega)}$, $\forall k \geq N(\omega).$
\end{thm}

\noindent \textbf{Preuve}\\
1. Remarquons que $\Omega = H  \cup S $ où 
 $$H  = \displaystyle
 \cap_{n }  \{T^{\ast}((\zeta_{n+1} , Y_{\widehat{\tau}_{n}}),.) >0\}\quad\mbox{et} \quad S  =\displaystyle \cup_{n }\, \{T^{\ast}((\zeta_n , Y_{\widehat{\tau}_{n}}),.) =0\}.$$
Sur l'ensemble $H $,  la suite $(\widehat{\tau}_{n})$ est  strictement croissante.  
  Sinon, sur l'ensemble complémentaire $S ,$ nous avons l'existence de 
 $N(\omega)$ tel que  $ \widehat{\tau}_{k} = \widehat{\tau}_{k+1} ,   ~  \forall k \geq N(\omega),$ et à partir de cet instant, la suite $(\widehat{\tau}_{n})$  s'arrête et nous continuons sans plus d'impulsion,  c'est à dire exactement (\ref{propadmis}). Ainsi, 
grâce à la construction explicitée par    la définition \ref{strategadm} de  la
famille $\widehat{\alpha}$, nous avons bien   une stratégie admissible .\\
\\
2. Du th\'eor\`eme \ref{stratoptimal}, nous sommes
ramen\'es \`a \'etablir  que la strat\'egie $\widehat{\alpha}$
satisfait les \'egalit\'es (\ref{egalite5}) et
(\ref{egalite6})  pour tout $(i,x) \in {U}\times \R$:\\

 i/ Sous les hypothèses 1 et 2(iii),   l'égalité
 (\ref{noyaubor}) de la proposition \ref{prnoyaubor}
 prise au point $(\zeta_{n},
Y_{\widehat{\tau}_n^-})$ implique:
$$m\rho^+(\zeta_{n},
Y_{\widehat{\tau}_n^-}) =
\int_{{U}\times\mathbb{R}}r^\ast
(\zeta_{n}, Y_{\widehat{\tau}_n^-};j,dy)\,\left(-c(\zeta_{n},
Y_{\widehat{\tau}_n^-},j,y) + \rho^+(j,y)\right).$$ 
D'où l'égalité
(\ref{egalite5}) est obtenue. \\

ii/ La partie droite de l'inégalité  (\ref{egalite6})  se décompose  selon la somme  suivante:
$$ \mathbb{E}_{\{\zeta_{n+1},
Y_{\widehat{\tau}_n}\}}\left[ \int_{\widehat{\tau}_n}^{\widehat{\tau}_{n+1}}
e^{-\beta s}f(\zeta_{n+1},Y_s)\, ds + e^{-\beta \widehat{\tau}_{n+1}}
m\rho^+(\zeta_{n+1},Y_{\widehat{\tau}_n^-})\right]$$
\begin{eqnarray*}
 & = &  \mathbb{E}_{\{\zeta_{n+1},
Y_{\widehat{\tau}_n}\}}\left[\left( \int_{\widehat{\tau}_n}^{\widehat{\tau}_{n+1}}
e^{-\beta s}f(\zeta_{n+1},Y_s)\, ds + e^{-\beta \widehat{\tau}_{n+1}}
m\rho^+(\zeta_{n+1},Y_{\widehat{\tau}_n^-})\right) \1_{\{\widehat{\tau}_{n+1}> \widehat{\tau}_{n}}\}\right] \nonumber\\
&+&\mathbb{E}_{\{\zeta_{n+1},
Y_{\widehat{\tau}_n}\}}\left[\left(\int_{\widehat{\tau}_n}^{+\infty}
e^{-\beta s}f(\zeta_{n+1},Y_s)\, ds \right)\1_{\{\widehat{\tau}_{n+1} = +\infty\}}\right].
\end{eqnarray*}
D'après  la définition \ref{strategadm}, la loi du couple $(\widehat{\tau}_{n+1} - \widehat{\tau}_{n}, (\zeta_{n+1},Y_.)\1_{[ \widehat{\tau}_{n}, \widehat{\tau}_{n+1}[})$ sachant $( \zeta_{n+1}, Y_{\widehat{\tau}_n}) $ est égale à  la $\mathbb{P}_{\{i,x\}}$-loi du couple $(T^\ast((i,x),.), (i,Y_.)\1_{[0, T^\ast((i,x),.)[})$ prise en $i = \zeta_{n+1}, x =Y_{\widehat{\tau}_n}$. Par suite, l'expression précédente  est égale à
$$ \mathbb{E}_{\{\zeta_{n+1},
Y_{\widehat{\tau}_n}\}}\left[ \int_{\widehat{\tau}_n}^{\widehat{\tau}_{n+1}}
e^{-\beta s}f(\zeta_{n+1},Y_s)\, ds + e^{-\beta \widehat{\tau}_{n+1}}
m\rho^+(\zeta_{n+1},Y_{\widehat{\tau}_n^-})\right]=$$
{\footnotesize{
\begin{eqnarray*}
&
\mathbb{E}_{\{\zeta_{n+1},
Y_{\widehat{\tau}_n}\}}\big[
\1_{(T^\ast(\zeta_{n+1},Y_{\widehat{\tau}_n}))>0}(\int_0^{T^\ast((\zeta_{n+1},Y_{\widehat{\tau}_n}))}
e^{-\beta  (s -\tau_n)}f(\zeta_{n+1}, Y_{\widehat{\tau_n}+s} ) ds
+ e^{-\beta (T^\ast((\zeta_{n+1},Y_{\widehat{\tau}_n}))-\tau_n) }
m\rho^+(\zeta_{n+1},Y_{T^\ast((\zeta_{n+1},Y_{\widehat{\tau}_n}))^-}))\big] 
\\
&+ \mathbb{E}_{\{\zeta_{n+1},
Y_{\widehat{\tau}_n}\}}\left[\1_{\{T^\ast(\zeta_{n+1},Y_{\widehat{\tau}_n},.)= 0\}}
\left(\int_{\tau_n}^{+\infty}
e^{-\beta (s -\tau_n)}f(i,Y_s)\, ds \right)\right].
\end{eqnarray*}
}}

D'après le lemme \ref{CI},  $T^\ast(\zeta_{n+1},Y_{\widehat{\tau}_n})> 0$  équivaut à $(\zeta_{n+1},Y_{\widehat{\tau}_n})\in C.$ 
 Sur cet événement,
 d'après l'égalité   (\ref{T}) du lemme \ref{propT},   le premier terme  ci-dessus 
  est exactement $\1_{(T^\ast(\zeta_{n+1},Y_{\widehat{\tau}_n})>0)}e^{-\beta \tau_n}\rho^+(\zeta_{n+1},Y_{\widehat{\tau}_n})$   puisque  $(\zeta_{n+1},Y_{\widehat{\tau}_n})\in C$
  et l'événement $(T^\ast(\zeta_{n+1},Y_{\widehat{\tau}_n})>0)\in \F_{\widehat\tau_n}.$
    
    Concernant le second terme, 
    $T^\ast(\zeta_{n+1},Y_{\widehat{\tau}_n})= 0$  entraine $\tau=\widehat{\tau_n}$:  l'intégrand 
  est exactement $e^{-\beta \tau}k_\tau^+$ d'où le second terme est 
  $\1_{(T^\ast(\zeta_{n+1},Y_{\widehat{\tau}_n})=0)}e^{-\beta \widehat{\tau_n}}  
  \rho^+(\zeta_{n+1},Y_{\widehat{\tau}_n}).$  
  
  La somme des deux termes est ainsi l'espérance 
  $$\mathbb{E}_{\{\zeta_{n+1},Y_{\widehat{\tau}_n}\}}[e^{-\beta \tau_n}\rho^+(\zeta_{n+1},Y_{\widehat{\tau}_n})]$$
  c'est à dire    l'égalité (\ref{egalite6}).  
\\
Ainsi, 
 les \'egalités  (\ref{egalite5}) et (\ref{egalite6})
sont  vérifiées  pour tout $(i,x) $ ce qui conclut l'optimalité de la statégie proposée. 
\endproof

\begin{exmp}
Nous donnons un exemple qui satisfait les hypothèses suffisantes de  l'existence d'une solution optimale.
  L'ensemble $U$ des technologies possibles est restreint à $\{0,1\}$ où $0$ est l'ancienne technologie et $1$ est la nouvelle technologie. Le processus $Y$ représentant le log de la valeur de la firme est un modèle de  Black-Scholes  satisfaisant:
 $$dY_t = b(\xi_t)dt + \sigma(\xi_t) dW_t,$$
 où $\xi_t$ indique la technologie au temps $t$. \\
  Nous supposons  que
 $b(1) > b(0)$ et $\sigma(1)> \sigma(0)$. Le profit net de la firme est donné   par   $f(i,x) = Arctg ~ x +\frac{\pi}{2}$ 
et le coût de changement de technologie  est égale à $$c(i,x,1-i,y) = 1-\frac{1}{1+ (y -x)^2}.$$
 \\
 Les lois conditionnelles des sauts $(\zeta_{n+1}, \Delta_n)$ sachant $\F_{\tau_n}$ sont $r(\zeta_n,Y_{\tau_n}^-;.,.)$ où $r$ appartient à l'ensemble des probabilités de transition suivant:
$$
   M_{(i,x)}= \left\{ r(i,x;1-i,dy) =  p_{i,j }\otimes \mathcal{N}(x+m,1) dy; [p_{i,j}]  = \left(
\begin{array}{ccc}
0 & 1 \\
1 & 0
\end{array}
\right)
,  m \in [\underline{m},\overline{m}] 
    \right\} \cup \{\delta_{i,x}\}.$$
Le paramètre  $m$ appartenant à un compact, les hypothèses 1 et 3 sont vérifiées  dans cet exemple.  De plus,  les lois conditionnelles r
 introduites ci-dessus sont  continues en x. Elles satisfont donc    l'hypothèse 2 (i) et (iii).  \\
 Les fonctions f et c étant bornées continues, l'hypothèse 2 (ii) est vérifiée. 

\end{exmp}

\section{Comparaison d'hypothèses}
\label{section5}
Dans cette section,   nous donnons une comparaison entre les hypothèses et définitions
utilisées tout au long de ce chapitre et celles utilisées par 
J.P.  Lepeltier et B. Marchal (\cite{lepl, lepeltier}). \\
   L'hypothèse 1 \cite[p. 48]{lepeltier} permet 
à J.P.  Lepeltier et B. Marchal de restreindre leur   modèle   à  la famille des lois markoviennes stationnaires.\\
 Nous nous plaçons dans le même cadre    en définissant les stratégies admissibles
avec des  lois de passage $r$  indépendantes de $n$.\\
\\
  L'hypothèse 2 \cite[p. 655]{lepl} impose que $\forall ~(i,x) M_{(i,x)}$ est faiblement relativement compact. 
 Nous renforçons  cette hypothèse en supposant que  l'ensemble des  probabilités de transition  $M_{(i,x)}$ est  fermé et  compact pour la topologie faible  (notre hypothèse 1). De plus, nous supposons que l'ensemble  $A= \{r(.,-c+\rho^+), r \in M\}$ est fermé et compact pour la topologie définie  par (\ref{topHyp3}) (notre hypothèse 3). \\
\\
L'hypothèse 3 \cite[p. 662]{lepl}  impose que  l'application $m^\ast\rho^+$ soit  semi-continue supérieurement sur les noyaux d'arrêt. 
 Dans notre cas,  notre hypothèse 2(iii) ( $x \To \nu(i,x;f)$ est  continue $\forall i \in U, \forall \nu \in M, ~ \forall f \in \mathcal{B}_b$)  permet de montrer que 
 l'application $m^\ast\rho^+$ est semi-continue inférieurement comme supremum de fonctions s.c.i.\\
\\
L'hypothèse 4 \cite[p. 662]{lepl}   suppose  que    $(i,x) \in I,$ l'application 
$$r(i,x,.,.)\longmapsto \int_{U \times \R}  r(i,x,j,dy)(-c(i,x,j,y) + \rho^+(j,y))$$
 atteint son minimum sur $M^\ast_{(i,x)} $. 
  Quant à nous,  cette propriété résulte  de  nos hypothèses 1 et 2  $\forall (i,x)$ et nous avons tenu
à en préciser la preuve en faisant appel à des résultats anciens sur les sélections mesurables
dans notre proposition \ref{prnoyaubor}. De plus, notre hypothèse 3 assure la continuité de l'application précédente. \\
\\
 J.P.  Lepeltier et B. Marchal imposent  les hypothèses 5 et 6 (\cite[p. 72]{lepl}) pour éliminer toute possibilité d'effectuer deux impulsions simultanées et elles s'énoncent comme suit: \\
 
Hypothèse  5 : pour tout noyaux borélien $r_1$ et $r_2$ de $M^\ast$, il existe un noyau    borélien $r\in M$ vérifiant:
$$  r(i_1,x_1; d(i_3,x_3)) = \int_{{U}\times \R} r_1(i_1,x_1; d(i_2,x_2))r_2 (i_2,x_2; d(i_3,x_3)).$$

Hypothèse 6: Pour tout $(i,x) \in {U}\times \R$, tout $\nu \in M^\ast_{(i,x)}$ et   $(j,y) \in {U}\times \R,$ on a l'inégalité suivante:
$$   \int_{{U}\times \R} \nu(d(w,z))( c(i,x,w,z) + c(w,z,j,y)) > c(i,x,j,y).$$
Ces deux hypothèses présentent une restriction de l'ensemble des stratégies. Dans notre modèle, notre définition \ref{strategadm} impose deux propriétés qui  restreignent à leur  tour l'ensemble des stratégies admissibles: la première assure l'intégrabilité des fonctions $f$ et $c$ tandis que la deuxième suppose que la loi entre deux sauts est donnée par la loi du couple initial
et assure des propriétés de continuité qui permettent d'éviter ces hypothèses 5 et 6 de \cite{lepeltier}. \\
   Ces   hypothèses 5 et 6 nous ont paru difficiles à vérifier dans un cadre concret:
 elles restreignent  l'ensemble des probabilités de transition  et 
   placent le modèle dans un cadre étroit. Nous avons jugé  plus réaliste de placer
    notre modèle markovien sous les propriétés  mentionnées dans la définition \ref{strategadm}.
    De plus,  nos hypothèses 1 et 2 
    et le théorème de sélection mesurable \cite[p. 85]{castaing}  ont permis de
     prouver l'existence d'un noyau borélien $r^\ast$ réalisant le supremum  de l'application $m \rho,$  mais aussi de donner les propriétés topologiques des ensembles
     $C$ et $I$ de continuation et d'impulsion. 
     \\
\\
L'hypothèse 7 \cite[p. 663]{lepl} indique que le coût de changement  d'état est minoré. 
  Cette propriété a été utilisée par les deux auteurs pour prouver l'admissibilité de la stratégie optimale. \\
\\
 Enfin, dans la thèse de J.P.  Lepeltier et B. Marchal \cite{lepeltier},  les fonctions  $f$ et $c$ sont boréliennes positives, fonctions de coûts qui amènent donc à une minimisation
 de valeurs positives. Puisque nous cherchons une optimisation d'une différence de valeurs
 positives, nous avons supposé    ces deux fonctions   continues bornées (Hypothèse 2 (ii)).  
  \\
L'hypothèse 2 (i) suppose   que  la loi du couple $(\tau_0,(i,Y_.))$ sachant $(i,x)$ est faiblement continue. Cette hypothèse plus l'hypothèse 2(iii) ont  contribué  à montrer  que les fonctions $\rho$ et $\rho^+$ sont semi-continues inférieurement.  Par conséquent, nous obtenons des propriétés topologiques sur l'ensemble d'impulsion $I$ pour s'assurer qu'aux  instants d'impulsion $ \tau_n$  le processus $(\xi_{\tau_n^-},Y_{\tau_n^-}) \in I$.

\section{Conclusion}
Dans ce chapitre, nous  considérons le problème du contrôle impulsionnel comme une suite de changements d'état et nous  montrons l'existence d'une solution optimale pour ce problème en utilisant une approche markovienne. Nous retrouvons dans un modèle constructif les résultats de J.P.  Lepeltier et B. Marchal (\cite{lepeltier,leplprob}) en utilisant  des outils différents. Au lieu de leur construction canonique, nous   construisons trajectoriellement le processus de l'évolution du système impulsé. 

Dans le prochain chapitre, nous utiliserons des exemples concrets pour présenter une solution optimale. Ainsi, on devrait être en mesure d'appliquer concrètement ces résultats. Dans un cas particulier (instants d'impulsion à intervalles réguliers et déterministes,  changement de technologie de l'ancienne à la nouvelle ou de la nouvelle vers l'ancienne), nous trouverons une loi de saut optimale. Nous constaterons que si l'on inclut les coûts de changement de technologie, il est préférable de sauter avec les précautions nécessaires pour réaliser plus de profit: ainsi, plus le saut est petit,  plus le profit est élevé.

\end{document}